\documentclass[12pt,letterpaper]{article}

\usepackage{renstyles}
\usepackage{comment}
\usepackage{cleveref}

\title{Recovering the initial condition and physical coefficients in a nonlinear PDE model of cell invasion}
\date{}
\author{
    Beiji Chen \thanks{Department of Applied Physics and Applied Mathematics, Columbia University, New York, NY 10027; bc3073@columbia.edu}
    \and 
    Kui Ren \thanks{Department of Applied Physics and Applied Mathematics, Columbia University, New York, NY 10027; kr2002@columbia.edu} 
} 

\begin{document}

\maketitle

\begin{abstract}

This paper investigates an inverse problem for the simultaneous reconstruction of two spatially varying reaction coefficients, the local proliferation rate and the competition (saturation) coefficient, together with the unknown initial condition, in a nonlinear, density-dependent reaction-diffusion model motivated by cell invasion and tumor growth dynamics. Using Carleman estimates, we establish a global uniqueness result together with a Lipschitz-type stability estimate for the reaction coefficients and a weaker, logarithmic stability estimate for the initial condition. For the numerical reconstructions, we develop a two-stage algorithm employing a time-shift strategy to decouple the coefficient and the initial condition. Numerical experiments are presented to illustrate the feasibility, accuracy, and robustness of the proposed inversion method.

\end{abstract}




\section{Introduction}

In recent years, there has been significant interest in calibrating mathematical models using experimental data. Such tasks are ubiquitous in the mathematical modeling of tumor growth~\cite{BaJoSv-IMA23,FaCoCaBa-arXiv23,Fuentes-arXiv25,HeKe-IP13,HuMo-arXiv24,LiSuMaCoBa-arXiv23}, propagation of epidemics~\cite{CaMoLa-BMB09,DuLa-arXiv24,KlLiYa-IPI25,NgYiZh-SIAM20,ZhVeDeKa-EC18}, traffic flow~\cite{KaAgSa-IEEE15,ThAmSaJa-TRC24}, and heat transfer~\cite{BoTsZa-IJHMT09,RoOrCoSzNoBiOs-CAMES17}, to name a few. Model calibration serves two main purposes. The first is to improve the accuracy of the models in describing the physical process behind the phenomenon to be modeled, and the second is to use the calibrated model to predict future dynamics~\cite{LiLoViFa-arXiv23,OdBaFa-ECM17,SwBrMuAl-JNS03,YaQuEv-CR15}.

In this work, we consider inverse problems for a model of cell invasion in brain tumor dynamics~\cite{FeLiZh-arXiv23,Murray2003,SwBrMuAl-JNS03,CoSu-AMC21,TuHaKuKa-MB25}:
\begin{equation}\label{EQ:Diff}
    \begin{cases}
            \partial_t \rho(x,t) -\nabla\cdot (\gamma(\rho(x,t)) \nabla \rho(x,t)) + f(x,\rho) = 0, & \mbox{in}\ \Omega\times (0,T)\\
            \rho(x,t) = b(x,t), & \mbox{on}\ \partial\Omega \times (0,T)\\
            \rho(x, 0) = \rho_0(x), & \mbox{in}\ \Omega
    \end{cases}
\end{equation}
where $\rho(x,t)$ denotes the tumor cell density. The term $\nabla \cdot (\gamma(\rho)\nabla\rho)$ characterizes the spatial spread of tumor cells via a \emph{density-dependent} diffusion coefficient $\gamma(\rho)$. The reaction term is modeled as 
\begin{equation}\label{EQ:Nonlinearity}
f(x,\rho) = \mu(x)\rho - \xi(x)\rho^{2}
\end{equation}
with $\mu(x)$ representing the local proliferation rate and $\xi(x)$ capturing the saturation effects arising from competition for space and resources. Finally, the boundary condition enforces tumor confinement within the brain domain $\Omega$, while the initial condition $\rho(x,0)=\rho_0(x)$ specifies the pre-existing tumor distribution.

Models like~\eqref{EQ:Diff} have been generalized to handle effects such as those due to the co-existence of multispecies and nonlocality and have been used in other areas of biological modeling; see, for instance, ~\cite{CaDeLa-ET23,CaEsMiTa-arXiv24,CoSaByMaLo-PRSA21,EnBuCa-PTRSA14,FaCoCaBa-arXiv23,Fede-DCDS17,GhMaBi-JMB16,HuMo-arXiv24,LiSuMaCoBa-arXiv23,PaWaZh-SIAM26,Smoller-Book12,ViLuGu-arXiv24} and references therein for some random samples of existing research in the area.

We are interested in reconstructing model parameters of~\eqref{EQ:Diff} from snapshot measurements of the density of the cell population $\rho$ at different observation times. More precisely, let $0 < t_0 < t_1 < T$ be two distinct time instants. We assume the available observation data include two full-domain snapshots at $t_0$ and $t_1$. Moreover, for technical convenience, we assume that we also have observations in a small subdomain $\omega\subset\Omega$ for all time in $(0, T)$. That is, the available data we have is:
\begin{equation}\label{EQ:Data}
    \begin{aligned}
        d_0 &:= \rho(x,t_0), && x\in \Omega\\
        d_1 &:= \rho(x,t_1), && x\in \Omega\\
        d_\omega &:= \rho(x,t), && (x,t)\in \omega\times (0,T).
    \end{aligned}
\end{equation}
Our objective is to simultaneously reconstruct the coefficients $\mu(x), \xi(x)$ and the initial distribution $\rho_0(x)$ from these data. Note that such snapshot data acquisition methods are now typical in many applications~\cite{BrBrNa-IJNM22,CaEsMiTa-arXiv24}.

Inverse problems for models like~\eqref{EQ:Diff} have been studied extensively with measured data of different types; see, e.g.,~\cite[Chapter 9]{Isakov-Book06} for an early summary of existing results. Most existing work focuses on reconstruction from boundary measurements, particularly via Dirichlet-to-Neumann maps~\cite{KiUh-ARMA23, WaLi-JCP18, ChLi-AML20,LiLiLo-arXiv25,LiLiLo-arXiv25B}. Klibanov~\cite{Klibanov-IP04} studied an inverse boundary value problem for general nonlinear diffusion equations $\partial_t u = F(u,\nabla u,x,t)$ by Carleman estimates. The reconstruction of nonlinear source terms from boundary data has been investigated in~\cite{CrKaNaRo-IP14,EgEnKl-IP04}, with corresponding numerical methods developed in~\cite{BodeSc-ESAIM21}. Uniqueness and stability for the recovery of nonlinear diffusion coefficients and general reaction terms from the Dirichlet-to-Neumann map have likewise been established in~\cite{EgPiSc-IP17,KiUh-ARMA23}. Meanwhile, inverse problems based on internal measurements have been investigated to a less extend. Such data arise naturally in biological applications and often lead to fundamentally different analytical and numerical frameworks~\cite{Bedr15, Cristofol_2006}. In particular, following Isakov’s seminal work~\cite{Isakov-CPAM91}, the reconstruction of linear parabolic equations from final-time internal observations has been extensively studied; see, for example,~\cite{Triki-JMPA21} and the references therein.

The interest in recovering the reaction term $f(x, \rho)$ in equation~(\ref{EQ:Diff}) originates from early works using time-trace data~\cite{PiRu-NMPDE87,PiRu-CPDE88}. More recently, Kaltenbacher and Rundell investigated several related inverse problems utilizing final-time data. Specifically, in~\cite{KaRu-IP19}, they recovered the spatial coefficient $q(x)$ in $u_t-\Delta u= q(x)g(u)$, assuming $g(u)$ is known. 
Additionally, the simultaneous reconstruction of two coefficients in a semilinear model was studied in~\cite{CrRo-IP13}. Closely related are the simultaneous recovery of the conductivity together with the nonlinear reaction term~\cite{KaRu-IPI20} and the uniqueness and reconstruction of a nonlinear diffusion term in a parabolic equation~\cite{KaRu-JMAA21}, the latter being particularly relevant to the density-dependent diffusion considered here. Lipschitz and H\"older stability for the simultaneous determination of quasilinear terms in parabolic equations was obtained in~\cite{ChKi-arXiv24}, and the stable determination of coefficients in semilinear parabolic systems in~\cite{AiBeHa-IP22}. Most closely related to the present work, Martinez and Vancostenoble~\cite{MaVa-DCDSS21} established Lipschitz stability for the growth-rate coefficient of a nonlinear Fisher--KPP equation from internal measurements.

The reconstruction of initial conditions in reaction-diffusion equations is likewise a fundamental problem, closely related to backward parabolic equations and data assimilation. On the theoretical front, uniqueness and stability results for recovering initial data from final-time or overdetermined observations have been established using Carleman estimates~\cite{ChYa-NA08}. Computationally, stable recovery from partial observations has been achieved through novel methods such as those based on Carleman estimates~\cite{BodeSc-ESAIM21,WaZhDuLi-DCDS18,WaLi-CMA12}. More recently, attention has extended to nonlinear and singularly perturbed reaction–diffusion-advection models~\cite{Lukyanenko-Math21}, as well as nonlinear parabolic problems utilizing Carleman-based numerical schemes, such as the Carleman–Newton methodology~\cite{Abhishek-JCAM24}. Further recent advances include the recovery of initial states in semilinear parabolic problems from time-averages~\cite{ScWa-arXiv24} and a Rothe-type, Carleman-based numerical framework for ill-posed initial-data problems in parabolic equations~\cite{Klibanov-arXiv24}. Despite these extensive studies, the simultaneous reconstruction of both coefficients and the initial condition remains less explored. Most existing results focus on identifying either the parameters, assuming a known initial state, or vice versa. Recent contributions toward joint recovery include the reconstruction of the initial condition for parabolic equations with Log-Lipschitz coefficients~\cite{DeSaPr-AMPA25} and the simultaneous identification of a piecewise-constant reaction coefficient together with the initial condition in a reaction-diffusion equation~\cite{WeZhHaZh-JMAA26}. We should also point out the interesting work in~\cite{ChNaWa-IPI25} for the linear diffusion model.

Recent works on inverse problems for reaction-convection-diffusion types of equations with nonlocal effect have also attracted a lot of attention in the community~\cite{CaEsMiTa-arXiv24,KoWa-RMS23,LiLiUh-arXiv22,LaLa-IPI22,LiZi-arXiv23}.

In this paper, we address the simultaneous reconstruction of the spatially varying reaction coefficients $\mu(x)$ and $\xi(x)$, alongside the unknown initial condition $\rho_0(x)$, in a nonlinear reaction–diffusion model. This simultaneous identification is inherently challenging due to the strong coupling between the coefficients and the initial state, which typically leads to severe ill-posedness. The closest prior result is that of Martinez and Vancostenoble~\cite{MaVa-DCDSS21}, who established unconditional Lipschitz stability for the growth-rate coefficient of a Fisher--KPP equation from internal measurements, but in the setting of a \emph{constant} diffusion coefficient and a \emph{known} initial state. Our analysis departs from theirs in three essential respects: the diffusion is genuinely density-dependent, $\gamma=\gamma(\rho)$; both reaction coefficients $\mu$ and $\xi$ are recovered simultaneously rather than a single growth rate; and the initial condition $\rho_0$ is itself unknown and reconstructed jointly with the coefficients. Our main contributions are threefold. First, we derive global uniqueness and stability results using Carleman estimates, demonstrating that the unknown parameters can be uniquely identified utilizing observations in snapshots at two distinct time instances, supplemented by a subdomain. Second, we develop a novel 'time-shift' strategy to overcome the coefficient-state coupling. By treating the snapshot at $t_0$ as a pseudo-initial condition, this approach decouples the inverse problem, enabling the recovery of reaction coefficients without a priori knowledge of the true initial state. Finally, we formulate a PDE-constrained optimization method based on this strategy, with numerical experiments confirming the algorithm's accuracy and stability.

The remaining sections are organized as follows. In Section~\ref{SEC:Pre}, we introduce the notation and functional setting and establish the well-posedness of the forward problem. Section~\ref{SEC:Unique} is devoted to the main uniqueness and stability results. In Section~\ref{SEC:algo}, we propose a two-stage reconstruction algorithm. Numerical experiments illustrating the theoretical analysis are presented in Section~\ref{SEC:Numer}.

\section{Well-posedness of forward problem}
\label{SEC:Pre}

\subsection{Notations and functional spaces}

We first fix some of the standard notations we will use throughout the paper.

We use $C$ to indicate the generic positive constants, which may depend on the domain or dimension. For any $\alpha \in (0,1)$ and non-negative integers $m, n$, we denote by  
\begin{equation*}
    C^{m+\alpha}(\overline{\Omega}), \quad C^{m+\alpha, \frac{m+\alpha}{2}}(\overline{\Omega} \times [0,T]), \quad C^{m+\alpha, \frac{m+\alpha}{2}}(\partial \Omega \times [0,T])
\end{equation*}
the standard Hölder spaces, as introduced in~\cite{Evans-Book10}. Likewise, we write  
\begin{equation*}
    H^{m}(\Omega), \quad H^{m}(0, T; H^n(\Omega))
\end{equation*}
to denote the Sobolev spaces and the time-dependent Sobolev--Bochner spaces, respectively, following the notation in~\cite{Adams-Book75}.



\subsection{Fréchet derivative and general mean-value property}

Let $F : X \to Y$ be a mapping between Banach spaces $X$ and $Y$. The Fréchet derivative of $F$ at a point $x_0 \in X$, denoted by $F'(x_0)$, is a bounded linear operator $F'(x_0): X \to Y$ satisfying
\[
\lim_{\|h\|\to 0}\frac{\|F(x_0+h)-F(x_0)-F'(x_0)(h)\|}{\|h\|}=0.
\]

Here, $h$ is a perturbation in $X$, and $F'(x_0)(h)$ represents the linear approximation of the change in $F$ at $x_0$. We denote by $C^1(X,Y)$ the space of continuously differentiable mappings $F : X \to Y$, meaning $F$ is continuous, and its Fréchet derivative $F'$ is bounded and continuous. If the Banach spaces $X$ and $Y$ coincide, we simplify the notation to $C^1(X)$.

We now introduce a general mean-value property valid for the space $C^1(X,Y)$ (see \cite{Deimling-Book13}): For any $F\in C^1(X,Y)$ and $x,y \in X$, there exists a constant $\tilde{\lambda}\in(0,1)$ such that
\begin{equation}\label{mvp}
    F(x)-F(y)=F'(\tilde{\lambda} x+(1-\tilde{\lambda})y)(x-y).
\end{equation}

\subsection{Well-posedness of forward problem}
\label{well-posedness}
In this paper, we assume that $\Omega$ is a smooth and bounded domain with $\partial \Omega$ in class $C^{4
}$ and the boundary function $b(x,t)$ defined in (\ref{EQ:Diff}) satisfies
\begin{equation*}
    b(x,t)\in C^{4
    +\alpha,\frac{4+\alpha}{2}}(\partial \Omega\times [0,T]), \quad b(x,t)>b_0 \quad \mbox{on}\ \partial\Omega \times [0,T]
\end{equation*}
with some constant $b_0>0$. For fixed constant $M>0$ and $r_0>0,$ we set
\begin{align*}
\mathcal{U} = \left\{
\begin{aligned}
&(\mu,\xi,\rho_0) \in C^{2+\alpha}(\overline{\Omega}) \times C^{2+\alpha}(\overline{\Omega}) \times C^{4+\alpha}(\overline{\Omega}) \; \text{such that} \\
&{\rho_0}|_{\partial \Omega} = b(\cdot, 0), \quad
\left.\nabla\cdot (\gamma(\rho_0) \nabla \rho_0) + \mu \rho_0-\xi \rho_0^2\right|_{\partial \Omega} = b_t(\cdot, 0),\\
&\|\mu\|_{C^{2+\alpha}(\overline{\Omega})} \leq M, \quad \mu > r_0 \; \text{on } \overline{\Omega},\\
&\|\xi\|_{C^{2+\alpha}(\overline{\Omega})} \leq M, \quad \xi > r_0 \; \text{on } \overline{\Omega},\\
&\|\rho_0\|_{C^{4+\alpha}(\overline{\Omega})} \leq M, \quad \rho_0 \geq r_0 > 0 \; \text{in } \overline{\Omega}.
\end{aligned}
\right\}
\end{align*}

Let us first recall some necessary results on the well-posedness of the forward model~\eqref{EQ:Diff}.

Local well-posedness results are easy to get by fixed point arguments; see, e.g,~\cite[Proposition 2.1]{FeKiUh-NA22}. For the degeneracy case, i.e., assuming that $\gamma(\rho)=\beta(\bx) \rho^\kappa$, well-posedness for such porous media equations is now standard~\cite{Aronson-NDP86,Vazquez-Book07, CaGhNa-arXiv21}. For related studies on inverse problems, we refer the reader to~\cite{KaSaTr-arXiv24,LiYu-SIAM19,CaMaVa-SIAM08,WaZhDuLi-DCDS18, CaMaVa-Book16}. In our study, we assume that $\gamma\in C^1(H^1(0,T;H^2(\Omega))).$ Furthermore, we assume that there exists $\underline\alpha$, $\overline\alpha$ such that $0<\underline \alpha \le \gamma(\rho)\le \overline \alpha<+\infty$. This assumption avoids the possibility of having degeneracy in the diffusion coefficient. Therefore, by the classical parabolic theory (see, e.g., \cite{Amann-Book95-1}), we have the following well-posedness theorem.

\begin{theorem} For any $(\mu, \xi, \rho_0)\in \mathcal{U},$ there exists a solution $\rho \in C^{4+\alpha,\frac{4+\alpha}{2}}(\overline{\Omega}\times [0,T])$ to problem (\ref{EQ:Diff}) satisfying the estimate
\begin{equation*}
 \|\rho\|_{C^{4+\alpha,2+\alpha/2}(\overline{Q})}
\leq C \left(
\|f\|_{C^{2+\alpha, 1+\alpha/2}(\overline{Q})}
+ \|b\|_{C^{4+\alpha,2+\alpha/2}(\partial\Omega \times [0,T])}
+ \|\rho_0\|_{C^{4+\alpha}(\overline{\Omega})}\right)\leq CM,
\end{equation*}
where $C$ depends only on $\Omega, T, \alpha, \underline\alpha$, and $\overline\alpha$.
\end{theorem}
By the maximum principle (see \cite{Friedman-Book08}), we have for any $x\in \overline{\Omega}$ and $t\in (0,T),$
\begin{equation}\label{maxprin}
    \rho(\mu,\xi,\rho_0)(x,t)\geq e^{-MT}\min\{\min_{x\in \overline{\Omega}} \rho_0(x), \min_{x\in \partial \Omega,\, t\geq0} b(x,t)\}\geq e^{-MT} \min\{r_0, b_0\}.
\end{equation}
Thus, $\rho(\mu,\xi,\rho_0)$ is strictly positive. Beyond this positivity, the analysis in Section~\ref{SEC:Unique} also requires a non-degeneracy condition on the time derivative of the reference solution. To this end, we assume that there exists a constant $c_0>0$ such that
\begin{equation}\label{ndeg}
    |\partial_t \rho(\mu,\xi,\rho_0)(x,t)| \geq c_0 \quad \text{on}\ \Omega\times (0,T).
\end{equation}
This condition cannot be deduced from~\eqref{maxprin}, since $\partial_t \rho$ may vanish even when $\rho$ is uniformly positive; it is needed for the change of variables $v = u/\partial_t \tilde\rho$ used in the proof of Theorem~\ref{stabi}. We note that this condition is a genuine restriction: it requires the reference solution to be strictly monotone in time and is most naturally satisfied during the growth phase, before the density saturates. If it holds only on a compact subset $\Omega'\Subset\Omega$, the stability estimates below remain valid with $\Omega$ replaced by $\Omega'$.

\section{Uniqueness and stability results}
\label{SEC:Unique}

In this section, we employ Carleman estimates to reconstruct the reaction coefficients $\mu$ and $\xi$, together with the initial condition $\rho_0$, using the data $d_0$, $d_1$, and $d_\omega$ defined in \eqref{EQ:Data}. In particular, the uniqueness of $\mu$ and $\xi$ can be established using only $d_0$, $d_1$, and the partial measurement $\rho|_{\omega \times (t_0 - \delta, t_1 + \delta)}$ for an arbitrarily small $\delta > 0$. However, the measurement $d_\omega$ remains indispensable for the recovery of the initial condition $\rho_0$.

We now explicitly state our main result:

\begin{theorem}\label{stabi}
Let $\tilde{\rho}$ be the solution with respect to coefficients $(\tilde{\mu},\tilde{\xi},\tilde\rho_0)\in \mathcal{U}$, i.e., $\tilde{\rho}$ satisfies
\begin{equation}\label{EQ:tirho}
    \begin{cases}
            \partial_t \tilde{\rho} -\nabla\cdot (\gamma(\tilde{\rho}) \nabla \tilde{\rho}) + \tilde{\rho} (\tilde{\mu}-\tilde{\xi}\tilde{\rho}) = 0, & \mbox{in}\ Q,\\
            \tilde{\rho}(x,t) = b(x,t), & \mbox{on}\ \Sigma,\\
            \tilde{\rho}(x, 0) = \tilde\rho_0(x), & \mbox{in}\ \Omega.
    \end{cases}
\end{equation}

We set
\begin{equation*}
    G(\rho,\tilde{\rho}) = \|\rho - \tilde{\rho}\|^2_{H^2(0, T, L^2(\omega))} + \|(\rho - \tilde{\rho})(t_0, \cdot)\|^2_{H^2(\Omega)} + \|(\rho-\tilde{\rho}) (t_1,\cdot)\|^2_{H^2(\Omega)}.
\end{equation*}
There exists a constant $C$ such that
\begin{align}
\label{estmu}
        \|\mu - \tilde{\mu}\|^2_{L^2(\Omega)}+\|\xi - \tilde{\xi}\|^2_{L^2(\Omega)} &\leq C G(\rho,\tilde{\rho})\\
\label{estnu}
        \|\rho_0 - \tilde\rho_0\|_{L^2(\Omega)} &\leq C \,\bigl|\log \left(G(\rho,\tilde{\rho})\right)\bigr|^{-1}.
\end{align}
\end{theorem}

\subsection{Carleman Estimate}

We first introduce a function $\psi \in C^2(\Omega)$ that satisfies
\[
\psi > 0 \text{ in } \Omega,\quad \psi = 0 \text{ on } \partial\Omega,\quad |\nabla \psi| > 0 \text{ on } \overline{\Omega\setminus\omega_0},
\]
where $\omega_0$ is a subdomain of $\omega$. The existence of such a function is established in~\cite{FuIm-Notes96}. Without loss of generality, we can assume that $t_0 =\frac{1}{3}T$ and $t_1 =\frac{2}{3}T$. Then we set $Q_i =(t_i-\delta, t_i+\delta)\times \Omega$ and $Q_{\omega, i} =(t_i-\delta, t_i+\delta)\times \omega$ with $\delta=\frac{1}{3}T$ for $i=0,1$. Fixing a sufficiently large constant $\lambda>0$, we define the weight functions
\[
\varphi_i(t,x) := \frac{e^{\lambda \psi(x)}}{(t-t_i+\delta)(t_i-t+\delta)}, \quad
\eta_i(t,x) := \frac{e^{\lambda \psi(x)} - e^{2\lambda \|\psi\|_{C(\bar{\Omega})}}}{(t-t_i+\delta)(t_i-t+\delta)}, \quad (x,t)\in Q,
\]
for $i=0,1.$ We note that $\eta_0(t_0,\cdot) = \eta_1(t_1,\cdot)$ in $\overline{\Omega}.$

With these definitions, we recall the following Carleman estimate, proved in~\cite[Lemma 1.2]{FuIm-Notes96}.

\begin{lemma}
\label{carle}
Let $a \in L^\infty(Q_i)$ satisfy $0 < a_0 \le a(x, t) \le a_1 < \infty$
for a.e. $(x, t) \in Q_i,$
for some constants $a_0, a_1$. Then there exist constants $s_0 \ge 1$
and $C > 0$, depending only on $\Omega$, $\omega_0$, $T$, $\lambda$,
$a_0$, and $a_1$, such that for every $s \ge s_0$ and $U \in H^1\bigl(t_i - \delta,\; t_i + \delta;\; H^2(\Omega) \cap H^1_0(\Omega)\bigr),$
the following estimate holds:
\begin{equation}\label{carleman}
\begin{aligned}
\int_{Q_i} \!\Bigl(
  \frac{1}{s\varphi_i}\bigl(|\partial_t U|^2 &+ |\Delta U|^2\bigr)
  + s \varphi_i\, |\nabla U|^2
  + s^3 \varphi_i^3\, |U|^2
\Bigr) e^{2 s \eta_i}\, dx\, dt \\
&\le\;
C \int_{Q_i} \bigl|\partial_t U - \nabla \!\cdot\! (a \nabla U)\bigr|^2 e^{2 s \eta_i}\, dx\, dt
\;+\;
C \int_{Q_{\omega, i}} s^3 \varphi_i^3\, |U|^2 e^{2 s \eta_i}\, dx\, dt.
\end{aligned}
\end{equation}
\end{lemma}  

\begin{remark}\label{rmk1}
Set $L_1^{(i)} V := \Delta V + s^2 |\nabla\varphi_i|^2 V - s (\partial_t \eta_i) V$ and $L_2^{(i)} V := \partial_t V + 2s \nabla\varphi_i \cdot \nabla V$ where $V = e^{s\eta_i} U,$ it follows from \cite{FuIm-Notes96} that the right-hand side of (\ref{carleman}) also gives an upper bound of $\|L_1^{(i)} V\|^2_{L^2(Q_i)} + \|L_2^{(i)} V\|^2_{L^2(Q_i)}$.
\end{remark}
\begin{remark}
    For simplicity, we denote 
\begin{equation*}
    I_i(U) := \int_{Q_i} \left( \frac{1}{s \varphi_i} \left| \partial_t U \right|^{2} + \frac{1}{s \varphi_i} |\Delta U|^{2} + s \varphi_i |\nabla U|^{2} + s^{3} \varphi_i^{3} |U|^{2} \right) e^{2s\eta_i} \, dx \, dt
\end{equation*}
\end{remark}

\subsection{Proof of Theorem \ref{stabi}}
We first denote $y := \rho - \tilde{\rho}$, $\alpha= \mu-\tilde{\mu}$, $\beta=\xi-\tilde{\xi}$, $h=\rho_0 - \tilde\rho_0$, $p=\rho +\tilde{\rho}$, $q = \gamma(\rho)$ and $R= \tilde{\rho}$. Since \( \gamma \in C^1(H^1(0,T, H^2(\Omega))) \), we recall the mean value property (\ref{mvp}), which states that there exists a constant \( \tilde{\lambda} \in (0,1) \) such that
\begin{equation*}
    \gamma(\rho) - \gamma(\tilde{\rho}) = \gamma'(\zeta)(\rho - \tilde{\rho})
\end{equation*}
with \( \zeta = \tilde{\lambda} \rho + (1 - \tilde{\lambda}) \tilde{\rho} \). Then the equation satisfied by $y$ can be derived by taking the difference of the equations (\ref{EQ:Diff}) and (\ref{EQ:tirho}) given as
\begin{equation}\label{yfpp}
    \begin{cases}
            \partial_t y -\nabla\cdot (q \nabla y) - \gamma'(\zeta)\nabla R\cdot \nabla y +(-\gamma'(\zeta)\Delta R+\mu-\xi p)y +\alpha R -\beta R^2= 0, & \mbox{in}\ Q,\\  
            y(x,t) = 0, & \mbox{on}\ \Sigma,\\
            y(x, 0) = h(x), & \mbox{in}\ \Omega.
    \end{cases}
\end{equation}
Here, the zeroth- and first-order coefficients of $y$ in~\eqref{yfpp} are displayed only schematically: the term $\nabla\!\cdot\!\big(\gamma'(\zeta)\,y\,\nabla R\big)$ contributes additional lower-order terms whose precise form is immaterial, since all such terms carry bounded coefficients on $\overline Q$ and are absorbed into the left-hand side of the Carleman estimate once $s$ is large. 

Applying Lemma \ref{carle} to the solution $y$ in (\ref{yfpp}), we can estimate $y$ by
\begin{equation}\label{yest}
\begin{aligned}
& \int_{Q_i}\left(\frac{1}{s \varphi_i}\left(\left|\frac{\partial y}{\partial t}\right|^{2}+|\Delta y|^{2}\right)+s \varphi_i|\nabla y|^{2}+s^{3} \varphi_i^{3} |y|^{2}\right)e^{2 s \eta_i} \, dx \, dt\\
& \leq C\left(\int_{Q_i}|\alpha - \beta R|^{2}e^{2 s \eta_i} \, dx \, dt+\int_{Q_{\omega, i}} s^{3} \varphi_i^{3} |y|^{2}e^{2 s \eta_i} \, dx \, dt\right).
\end{aligned}
\end{equation}
In particular, the term involving $y$ and $\nabla y$ can be absorbed into the left-hand side of \eqref{yest}, provided that the parameter $s$ is sufficiently large. This absorption technique will be used repeatedly throughout the proof, relying on the magnitude of $s$.

By introducing the transformation $z = \frac{y}{R}$, and defining
\begin{equation*}
    A_1 = 2q\frac{\nabla R}{R} + \gamma'(\zeta)\nabla R, \quad
A_2 = \frac{\partial_t R}{R} - q \frac{\Delta R}{R} - \nabla q \cdot \frac{\nabla R}{R} - \gamma'(\zeta)\frac{|\nabla R|^2}{R} - \gamma'(\zeta) \Delta R + \mu - \xi p,
\end{equation*}
we reformulate equation~\eqref{yfpp} as
\begin{equation}\label{zfpp}
    \begin{cases}
        \partial_t z - \nabla \cdot (q \nabla z) -A_1 \cdot \nabla z + A_2 z + \alpha - \beta R = 0, & \text{in } Q, \\
        z(x,t) = 0, & \text{on } \Sigma, \\
        z(x,0) = \dfrac{h(x)}{R(x)}, & \text{in } \Omega.
    \end{cases}
\end{equation}
Applying Lemma \ref{carle} to the solution $z$ in (\ref{zfpp}), we derive the estimate for $z$ as
\begin{equation}\label{zest}
\begin{aligned}
& \int_{Q_i}\left(\frac{1}{s \varphi_i}\left(\left|\frac{\partial z}{\partial t}\right|^{2}+|\Delta z|^{2}\right)+s \varphi_i|\nabla z|^{2}+s^{3} \varphi_i^{3} |z|^{2}\right)e^{2 s \eta_i} \, dx \, dt\\
& \leq C\left(\int_{Q_i}|\alpha - \beta R|^{2}e^{2 s \eta_i} \, dx \, dt+\int_{Q_{\omega, i}} s^{3} \varphi_i^{3} |z|^{2}e^{2 s \eta_i} \, dx \, dt\right).
\end{aligned}
\end{equation}
Set $u=\partial_t z$ and differentiate the equation (\ref{zfpp}) with respect to $t$. We then obtain
\begin{equation}\label{ufpp}
    \begin{cases}
            \partial_t u -\nabla\cdot (q \nabla u) - A_1\cdot \nabla u + A_2 u -\nabla\cdot (\gamma_t(\rho)\nabla z)\\-(\partial_t A_1)\cdot \nabla z +(\partial_t A_2)z -\beta \partial_t R= 0, & \mbox{in}\ Q,\\
            u(x,t) = 0, & \mbox{on}\ \Sigma,\\
            u(x, 0) = (\nabla\cdot (q\nabla z)+A_1\cdot \nabla z -A_2 z -\alpha +\beta R)(x,0), & \mbox{in}\ \Omega.
    \end{cases}
\end{equation}
To estimate $u$, we recall~\cite[Lemma 2.1]{Klibanov-IP04}.
\begin{lemma}\label{kli}
There exists a constant $C>0$ depending only on $t_i$ and a constant $C_s>0$ depending on $s$ such that 
\begin{equation}
    \int_{Q_i} e^{2s\eta_i} \left|\int_{t_i}^t g(\tau,x) d\tau\right|^2 \, dt \, dx \leq \frac{C}{s} \|e^{s\eta_i}g\|^2_{L^2(Q_i)} + C_s \|\tilde{g}(t_i,\cdot)\|^2_{L^2(\Omega)}.
\end{equation}
for all $g\in L^2(Q_i), s>0$ and $\partial_t \tilde{g} =g$.
\end{lemma}
Using Lemma~\ref{carle} and Lemma~\ref{kli}, we derive the following estimate for \( u \):
\begin{equation}\label{uest}
\begin{aligned}
& \int_{Q_i}\left(\frac{1}{s \varphi_i}\left(\left|\frac{\partial u}{\partial t}\right|^{2}+|\Delta u|^{2}\right)+s \varphi_i|\nabla u|^{2}+s^{3} \varphi_i^{3} |u|^{2}\right)e^{2 s \eta_i} \, dx \, dt\\
& \leq C_s\left(\int_{Q_{\omega, i}} s^{3} \varphi_i^{3} |z|^{2}e^{2 s \eta_i} \, dx \, dt +\|z(t_i,\cdot)\|^2_{H^2(\Omega)}\right)+C \int_{Q_i} |\beta|^2 e^{ 2 s \eta_i} \, dx \, dt.
\end{aligned} 
\end{equation}
We denote
\begin{equation*}
\begin{aligned}
        &A_3 = 2q \frac{\nabla \partial_t R}{\partial_t R} + A_1,\\
        &A_4 = q \frac{\Delta \partial_t R}{\partial_t R} +(\nabla q + A_1) \cdot \frac{\nabla\partial_t R}{\partial_t R} + A_2,
\end{aligned}
\end{equation*}
and  introduce the linear operators $P_1$ and $P_2$ by
\begin{equation*}
    \begin{aligned}
        &P_1(v) := - \nabla\cdot (\partial_t q)\nabla v +\partial_t A_3 \cdot \nabla v +(\partial_t A_4)v,\\
        &P_2(z) := \partial_t \left(\frac{\nabla\cdot (\partial_t q)\nabla z+(\partial_t A _1)\cdot \nabla z + (\partial_t A_2)z}{\partial_t R}\right).
    \end{aligned}
\end{equation*}
Setting $v = \frac{u}{\partial_t R}$ and $\omega = \partial_t v$, we find that $\omega$ satisfies \begin{equation}\label{omfpp}
    \begin{cases}
            \partial_t \omega -\nabla\cdot (q \nabla \omega) + A_3\cdot \nabla \omega +A_4 \omega + P_1(v) + P_2(z)= 0, & \mbox{in}\ \Omega\times (0,T),\\
            \omega(x,t) = 0, & \mbox{on}\ \partial\Omega \times (0,T).
    \end{cases}
\end{equation}
From Lemma~\ref{carle}, it follows that $\omega$ satisfies the estimate
\begin{equation}\label{pom}
\begin{aligned}
& \int_{Q_i}\left(\frac{1}{s \varphi_i}\left(\left|\frac{\partial \omega}{\partial t}\right|^{2}+|\Delta \omega|^{2}\right)+s \varphi_i|\nabla \omega|^{2}+s^{3} \varphi_i^{3} |\omega|^{2}\right)e^{2 s \eta_i} \, dx \, dt\\
& \leq C_s\left(\int_{Q_{\omega, i}} s^{3} \varphi_i^{3} |\omega|^{2}e^{2 s \eta_i} \, dx \, dt +\int_{Q_i} \left(|P_1(v)|^2 + |P_2(z)|^2\right)e^{ 2 s \eta_i} \, dx \, dt\right).
\end{aligned} 
\end{equation}
Furthermore, Lemma \ref{kli} leads to
\begin{equation}\label{p2z}
    \int_{Q_i} |P_2(z)|^2 e^{ 2 s \eta_i} \, dx \, dt \leq \frac{C}{s} \int_{Q_i} (|u|^2+|\nabla u|^2 + |\Delta u|^2) e^{2 s \eta_i} \, dx \, dt + C_s \|z(t_i, \cdot)\|^2_{H^2(\Omega)}.
\end{equation}
By combining equations~\eqref{uest}, \eqref{pom}, and \eqref{p2z}, we derive an estimate for $\omega$ by
\begin{equation}\label{omest}
\begin{aligned}
& \int_{Q_i}\left(\frac{1}{s \varphi_i}\left(\left|\frac{\partial \omega}{\partial t}\right|^{2}+|\Delta \omega|^{2}\right)+s \varphi_i|\nabla \omega|^{2}+s^{3} \varphi_i^{3} |\omega|^{2}\right)e^{2 s \eta_i} \, dx \, dt\\
& \leq C_s\int_{Q_{\omega, i}} s^{3} \varphi_i^{3} \left(|\omega|^{2} + |u|^2\right) e^{2 s \eta_i} \, dx \, dt + \frac{C}{s} \int_{Q_i} |\beta|^2 e^{ 2 s \eta_i} \, dx \, dt +C_s\|z(t_i, \cdot)\|^2_{H^2(\Omega)}.
\end{aligned}     
\end{equation}

We now proceed to prove estimate~\eqref{estmu}, aiming to bound $\|\mu - \tilde{\mu}\|^2_{L^2(\Omega)} + \|\xi - \tilde{\xi}\|^2_{L^2(\Omega)}$. As a first step, we estimate $\|\alpha e^{ s \eta_0 (t_0, \cdot)}\|^2_{L^2(\Omega)}$. To this end, we recall Lemma 4.6 in \cite{CrKaNaRo-IP14}.
\begin{lemma}\label{crka}
    There exists a constant $s_1$ depending only on $t_0$, such that, for all $w\in H^1((t_0 - \delta, t_0 + \delta),L^2(\Omega))$ and $s\geq s_1,$ we have
    \begin{equation*}
        \int_{\Omega} |w(t_0, x)|^2 \, dx \leq 2\left(s \int_{Q_0} |w(t, x)|^2 \, dx \, dt + s^{-1}\int_{Q_0} |\partial_t w(t,x)|^2 \, dx \, dt \right).
    \end{equation*}
\end{lemma}

By choosing $w(t,x) = u(x,t) e^{s\eta_0(t_0,x)}$ in Lemma \ref{crka}, we obtain the following estimate:
\begin{equation*}
    \int_{\Omega} |u(t_0, x)|^2 e^{2s \eta_0(t_0, x)} \, dx \leq 2s \int_{Q_0} |u(t,x)|^2 e^{ 2 s \eta_0 (t_0,x)} \, dx \, dt + \frac{2}{s} \int_{Q_0} |\partial_t u(t,x)|^2 e^{ 2 s \eta_0(t_0, x )} \, dx \, dt
\end{equation*}
for $s$ sufficiently large. Evaluating equation~\eqref{zfpp} at $t = t_0$ yields the identity
\begin{equation*}
    \alpha - \beta R(t_0,\cdot) = - u(t_0,\cdot) + \nabla\cdot(q\nabla z)(t_0,\cdot) + A_1\cdot\nabla z(t_0,\cdot) - A_2\, z(t_0,\cdot).
\end{equation*}
Combining this with Lemma~\ref{crka} applied to $u$, we obtain
\begin{equation*}
\begin{aligned}
    \|(\alpha - \beta R(t_0,\cdot))\, e^{s \eta_0(t_0, \cdot)}\|^2_{L^2(\Omega)}
        &\leq C\left(\|u(t_0, \cdot) e^{s \eta_0(t_0, \cdot)}\|^2_{L^2(\Omega)}+\|z(t_0,\cdot) e^{s \eta_0(t_0, \cdot)}\|^2_{H^2(\Omega)}\right)\\
        & \leq C\left(s^{-1}\|\omega e^{s \eta_0}\|^2_{L^2(t_0-\delta, t_0+\delta , L^2(\Omega))} + s \|u e^{s \eta_0}\|^2_{L^2(t_0-\delta, t_0+\delta, L^2(\Omega))} \right.\\
        &\left.+\|z(t_0,\cdot) e^{s \eta_0(t_0, \cdot)}\|^2_{H^2(\Omega)}\right).
\end{aligned}
\end{equation*}
Then it follows from the estimate of $u$ in (\ref{uest}) that
\begin{equation*}
\|(\alpha - \beta R(t_0,\cdot))\, e^{s \eta_0(t_0, \cdot)}\|^2_{L^2(\Omega)} \leq C \left(\frac{I_0(\omega)}{s^4} + \frac{I_0(u)}{s^2} + \|z(t_0, \cdot)e^{s\eta_0(t_0,\cdot)}\|^2_{H^2(\Omega)}\right).
\end{equation*}
From the Carleman estimate in Lemma \ref{carle}, we conclude that there exist constants $s_2>0$ and $C>0$ such that for all $s>s_2$, we have
\begin{equation}\label{aest1}
\begin{aligned}
        \|(\alpha - \beta R(t_0,\cdot))\, e^{s \eta_0(t_0, \cdot)}\|^2_{L^2(\Omega)} &\leq C_s\left( \int_{Q_{\omega_0}} s^3\varphi_0^3 (|\omega|^2+|u|^2) e^{ 2s \eta_0} \, dx \, dt + \|z(t_0,\cdot)\|^2_{H^2(\Omega)}\right)
    \\&+ \frac{C}{s^2} \|\beta e^{s \eta_0(t_0,\cdot)}\|^2_{L^2(\Omega)}
\end{aligned}
\end{equation}
where constant $C_s>0$ depends on $s$.

After establishing the first estimate, we now focus on bounding the second term, $\|(\alpha - \beta R(t_1, \cdot))\, e^{s \eta_1(t_1, \cdot)}\|^2_{L^2(\Omega)}.$
To this end, we introduce the function $E = e^{s\eta_1} u$. Noting that $ E(t_1 - \delta, \cdot) = 0 $ and using the identity $\nabla \varphi_1 = \lambda \varphi_1 \nabla \psi$, we obtain
\begin{equation*}
\begin{aligned}
     \|e^{s \eta_1(t_1, \cdot)}u(t_1,\cdot)\|^2_{L^2(\Omega)}    &= \int_{t_1-\delta}^{t_1} \int_{\Omega} \partial_t |E(t,x)|^2 \, dx \, dt \\
     &=2 \left(\int_{t_1-\delta}^{t_1} \int_{\Omega}  (L_2 E) E \, dx \, dt
      + s \int_{t_1-\delta}^{t_1} \int_{\Omega}  \Delta \varphi_1 E^2 \, dx \, dt \right),
\end{aligned}
\end{equation*}
where $L_2$ is defined in Remark \ref{rmk1}. This yields
\begin{equation*}
    \int_{\Omega} |u(t_1,x)|^2e^{2 s \eta_1(t_1, x)} dx \leq 2\left|\int_{t_1-\delta}^{t_1} \int_{\Omega} (L_2 E) E \, dx\, dt\right| + C s \int_{t_1-\delta}^{t_1} \int_\Omega e^{2s\eta_1} \varphi_1 |u|^2 \, dx \, dt.
\end{equation*}
On the other hand, by Young's inequality, we obtain
\begin{equation*}
    \left|\int_{t_1-\delta}^{t_1} \int_{\Omega} (L_2 E) E \, dx\, dt\right| \leq \frac{1}{2} s^{-\frac{3}{2}}\left(\|L_2E\|_{L^2(\Omega\times(t_1-\delta,t_1))}^2 + s^3\int_{t_1-\delta}^{t_1}\int_{\Omega} e^{2s\eta_1} \varphi_1^3  |u|^2 \, dx\, dt\right).
\end{equation*}
Then it follows from Remark \ref{rmk1} that
\begin{equation*}
    \left|\int_{t_1-\delta}^{t_1} \int_{\Omega} (L_2 E) E \, dx\, dt\right| \leq C\left(s^{\frac{3}{2}}\int_{t_1-\delta}^{t_1} \int_{\omega}e^{2s\eta_1} \varphi_1^3 |u|^2 \, dx \,  dt +s^{-\frac{3}{2}}\int_{t_1-\delta}^{t_1} \int_{\Omega} e^{2s\eta_1} |\beta|^2 \, dx\, dt\right). 
\end{equation*}
Therefore, we obtain
\begin{equation}\label{esttu}
    \|e^{s \eta_1}u(t_1,\cdot)\|^2_{L^2(\Omega)}   \leq C_s
    \int_{Q_{\omega_1}} e^{2s\eta_1} \varphi_1^3 |u|^2 \, dx \, dt + Cs^{-\frac{3}{2}}\int_{t_1-\delta}^{t_1}\int_{\Omega} e^{2s\eta_1}|\beta|^2\, dx \, dt.
\end{equation}
Now return to equation~(\ref{zfpp}) evaluated at time $t=t_1$, i.e., we consider
\begin{equation*}
  u(t_1,\cdot) = \left(\nabla\cdot (q\nabla z) + A_1 \cdot \nabla z + A_2 z +\alpha -\beta R\right)(t_1,\cdot). 
\end{equation*}
Multiplying both sides by $e^{s\eta_1(t_1,\cdot)}$ on both sides and combining this identity with estimate~\eqref{esttu}, we deduce that there exist constants $s_3>0$ and $C >0$ such that for all $s>s_3,$ the following inequality holds:
\begin{equation}\label{abest.}
\begin{aligned}
        \|(\alpha - \beta R(t_1,\cdot))e^{s \eta_1(t_1,\cdot)}\|^2_{L^2(\Omega)} 
        &\leq C_s \left(\int_{Q_{\omega_1}} \varphi_1^3 e^{ 2 s \eta_1 } |u|^2 \, dx \,  dt+\|z(t_1,\cdot)\|^2_{H^2(\Omega)}\right) \\&+s^{-\frac{3}{2}} \int_{t_1 - \delta}^{t_1} \int_{\Omega} e^{2s\eta_1} |\beta|^2 \, dx \, dt.
\end{aligned}
\end{equation}
where constant $C_s>0$ depends on $s$.

We have already established the estimates~\eqref{aest1} for $\|(\alpha - \beta R(t_0,\cdot))\,e^{s\eta_0(t_0,\cdot)}\|^2_{L^2(\Omega)}$ and~\eqref{abest.} for $\|(\alpha - \beta R(t_1,\cdot))e^{s \eta_1(t_1,\cdot)}\|^2_{L^2(\Omega)}$ as a first step toward bounding $\|\mu - \tilde{\mu}\|^2_{L^2(\Omega)} + \|\xi - \tilde{\xi}\|^2_{L^2(\Omega)}.$

By the triangle inequality and the identity $\eta_0(t_0,\cdot) = \eta_1(t_1,\cdot)$, we have
\begin{equation*}
\begin{aligned}
    \|\beta(R(t_1,\cdot) - R(t_0,\cdot))\,e^{s\eta_0(t_0,\cdot)}\|^2_{L^2(\Omega)}
    &\leq 2\,\|(\alpha - \beta R(t_0,\cdot))\,e^{s\eta_0(t_0,\cdot)}\|^2_{L^2(\Omega)}\\
    &\quad + 2\,\|(\alpha - \beta R(t_1,\cdot))\,e^{s\eta_1(t_1,\cdot)}\|^2_{L^2(\Omega)}.
\end{aligned}
\end{equation*}
The non-degeneracy assumption~\eqref{ndeg} implies that $\partial_t R$ has constant sign on $\Omega\times(0,T)$. Hence
\begin{equation*}
    |R(t_1,x) - R(t_0,x)| = \left|\int_{t_0}^{t_1}\partial_t R(x,\tau)\,d\tau\right| \geq c_0(t_1 - t_0) =: c_1 > 0 \qquad \text{on}\ \Omega,
\end{equation*}
which gives the lower bound
\begin{equation*}
    c_1^2\,\|\beta\,e^{s\eta_0(t_0,\cdot)}\|^2_{L^2(\Omega)} \leq \|\beta(R(t_1,\cdot) - R(t_0,\cdot))\,e^{s\eta_0(t_0,\cdot)}\|^2_{L^2(\Omega)}.
\end{equation*}
Combining the above with~\eqref{aest1} and~\eqref{abest.}, and noting that $\eta_1(t,\cdot) \leq \eta_1(t_1,\cdot) = \eta_0(t_0,\cdot)$ on $(t_1-\delta,t_1)$ implies $\int_{t_1-\delta}^{t_1}\!\int_\Omega |\beta|^2 e^{2s\eta_1}\,dx\,dt \leq \delta\,\|\beta\,e^{s\eta_0(t_0,\cdot)}\|^2_{L^2(\Omega)}$, we obtain
\begin{equation*}
\begin{aligned}
    c_1^2\,\|\beta\,e^{s\eta_0(t_0,\cdot)}\|^2_{L^2(\Omega)}
    &\leq C_s\biggl(\int_{Q_{\omega_0}} s^3\varphi_0^3(|\omega|^2 + |u|^2)\,e^{2s\eta_0}\,dx\,dt + \int_{Q_{\omega_1}}\varphi_1^3 e^{2s\eta_1}|u|^2\,dx\,dt\biggr)\\
    &\quad + C_s\bigl(\|z(t_0,\cdot)\|^2_{H^2(\Omega)} + \|z(t_1,\cdot)\|^2_{H^2(\Omega)}\bigr)
    + \frac{C}{s^{3/2}}\,\|\beta\,e^{s\eta_0(t_0,\cdot)}\|^2_{L^2(\Omega)}.
\end{aligned}
\end{equation*}
For $s$ sufficiently large, the last term on the right is absorbed into the left-hand side, yielding a bound on $\|\beta\,e^{s\eta_0(t_0,\cdot)}\|^2_{L^2(\Omega)}$ in terms of the data. Combining this with~\eqref{aest1} and the elementary inequality
\begin{equation*}
\|\alpha\,e^{s\eta_0(t_0,\cdot)}\|^2_{L^2(\Omega)} \leq 2\,\|(\alpha-\beta R(t_0,\cdot))\,e^{s\eta_0(t_0,\cdot)}\|^2_{L^2(\Omega)} + 2\|R(t_0,\cdot)\|^2_{L^\infty(\Omega)}\,\|\beta\,e^{s\eta_0(t_0,\cdot)}\|^2_{L^2(\Omega)}
\end{equation*}
yields a corresponding bound on $\|\alpha\,e^{s\eta_0(t_0,\cdot)}\|^2_{L^2(\Omega)}$.

Since $e^{2s\eta_0(t_0,\cdot)}$ is bounded above and below by positive constants on $\overline{\Omega}$ (with constants depending on $s$ and $\lambda$), we obtain
\begin{equation}\label{ineq:alpha_and_beta}
\begin{aligned}
        \|\alpha\|^2_{L^2(\Omega)} + \|\beta\|^2_{L^2(\Omega)}
        &\leq C\left(\int_{0}^{T} \int_{\omega} (|u|^2+|\omega|^2) \, dx \, dt +\|z(t_0,\cdot)\|^2_{H^2(\Omega)}+\|z(t_1,\cdot)\|^2_{H^2(\Omega)}\right)\\
        &\leq C \left(\|\rho-\tilde{\rho}\|^2_{H^2(0,T,L^2(\Omega))} +\|(\rho-\tilde{\rho})(t_0,\cdot)\|^2_{H^2(\Omega)}+\|(\rho-\tilde{\rho})(t_1,\cdot)\|^2_{H^2(\Omega)}\right)
\end{aligned}
\end{equation}

Next, we turn to the proof of (\ref{estnu}). To this end, we denote $\tilde{y} = \partial_t y$, $p_t = \partial_t p$, $q_t = \partial_t q$ and $R_t = \partial_t R$. Differentiating equation (\ref{yfpp}) with respect to $t$, we obtain
\begin{equation}\label{zlz}
    \begin{cases}
            \partial_t \tilde{y} -\nabla\cdot (q \nabla \tilde{y}) + \gamma'(\zeta)\nabla R\cdot \nabla \tilde{y} +(\gamma'(\zeta)\Delta R+\mu-\xi p)\tilde{y}\\
            -\nabla\cdot(q_t\nabla y) + \gamma'(\zeta)\nabla R_t \cdot \nabla y +(\gamma'(\zeta)\Delta R_t-\xi p_t) y + \alpha R_t - 2 \beta R R_t= 0, & \mbox{in}\ Q,\\
            \tilde{y}(x,t) = 0, & \mbox{on}\ \Sigma,\\
            \tilde{y}(x, t_0) = r(x), & \mbox{in}\ \Omega,
    \end{cases}
\end{equation}
where $r(x)=(\nabla\cdot (q \nabla y) - \gamma'(\zeta)\nabla R\cdot \nabla y -(\gamma'(\zeta)\Delta R+\mu-\xi p)y -\alpha R +\beta R^2)(x, t_0).$ Applying Lemma~\ref{carle} to equation~\eqref{zlz}, we obtain
\begin{equation}\label{est2}
\begin{aligned}
\int_{Q_0} \left( \frac{1}{s \varphi_0} |\partial_t \tilde{y}|^2 + s \varphi_0 |\nabla \tilde{y}|^2 + s^3 \varphi_0^3 |\tilde{y}|^2 \right) e^{2s\eta_0} \, dx \, dt 
\leq  C \int_{Q_{\omega_0}} s^3 \varphi_0^3 |\tilde{y}|^2 e^{2s\eta_0} \, dx \, dt \\
+C \int_{Q_0} \left( |\alpha|^2 + |\beta|^2 + |\Delta y|^2+ |\partial_t y|^2 + |\nabla \partial_t y|^2 \right) e^{2s\eta_0} \, dx \, dt.
\end{aligned}
\end{equation}
Noting that $q(x,t) \geq \underline\alpha$ on $Q$, it follows that
\[
\int_{Q_0} |\Delta y|^2 e^{2s\eta_0} \, dx \, dt 
\leq C \int_{Q_0} |q \Delta y|^2 e^{2s\eta_0} \, dx \, dt.
\]
Moreover, we substitute $q\Delta y$ in equation~\eqref{yfpp} and then get the following estimate:
\begin{equation}\label{Laplace}
\begin{aligned}
\int_{Q_0} |\Delta y|^2 e^{2s\eta_0} \, dx \, dt 
\leq C \int_{Q_0} \left( |\alpha|^2 + |\beta|^2 + |\partial_t y|^2 + |\nabla y|^2 + |y|^2  \right) e^{2s\eta_0} \, dx \, dt.
\end{aligned}
\end{equation}
Based on the estimates in (\ref{yest}), (\ref{est2}), and (\ref{Laplace}), together with the absorption of terms involving \( y \), \( \nabla y \), \( \partial_t y \), and \( \nabla \partial_t y \), we obtain the following lemma:
\begin{lemma}\label{ttctrl}
Let $y \in H^{1}(t_0-\delta,t_0+\delta,H^2(\Omega))$ satisfy equation (\ref{yfpp}). Then, there exist constants $s_4 > 0$ and $C > 0$ such that for all $s > s_4$,
\[
\int_{Q_0} \left(\frac{1}{s \varphi_0} \left|\partial^2_{t} y\right|^{2} + s^{3} \varphi_0^{3} |\partial_t y|^{2}\right) e^{2 s \eta_0}  \, dx \, dt \leq C\int_{Q_0} (|\alpha|^2 + |\beta|^2) e^{2 s \eta_0}  \, dx \, dt + C H ,
\]
where
\[
H = \int_{Q_{\omega_0}} s^{3} \varphi_0^{3} (|y|^{2} + |\partial_t y|^2) e^{2 s \eta_0}  \, dx \, dt.
\]
\end{lemma}
Applying Lemma~\ref{ttctrl} to estimate~(\ref{Laplace}) yields
\begin{equation}\label{Hy}
    \int_{Q_0} (|\Delta y|^2 +|\nabla y|^2 + |y|^2) e^{2s\eta_0} \, dx \, dt 
\leq C\int_{Q_0} (|\alpha|^2 + |\beta|^2) e^{2 s \eta_0}  \, dx \, dt + C H.
\end{equation}

We now proceed to decompose (\ref{zlz}) as follows:
\begin{equation}\label{dev}
\begin{cases}
    \partial_t \tilde{v} -\nabla\cdot (q \nabla \tilde{v}) + \gamma'(\zeta)\nabla R\cdot \nabla \tilde{v} +(\gamma'(\zeta)\Delta R+\mu-\xi p)\tilde{v}\\
    -\nabla\cdot(q_t\nabla y) + \gamma'(\zeta)\nabla R_t \cdot \nabla y +(\gamma'(\zeta)\Delta R_t-\xi p_t) y+ \alpha R_t - 2 \beta R R_t = 0 &\mbox{in}\ Q,\\
    \tilde{v}(x,t) = 0 & \mbox{on}\ \Sigma,\\
    \tilde{v}(x,0) = 0 & \mbox{in}\ \Omega.
\end{cases} 
\end{equation}
and
\begin{equation}\label{deu}
\begin{cases}
    \partial_t \tilde{u} -\nabla\cdot (q \nabla \tilde{u}) + \gamma'(\zeta)\nabla R\cdot \nabla \tilde{u} +(\gamma'(\zeta)\Delta R+\mu-\xi p)\tilde{u} =0 &\mbox{in}\ Q,\\
    \tilde{u}(x,t) = 0 &\mbox{on}\ \Sigma,\\
    \tilde{u}(x, t_0) = r(x) - \tilde{v}(x, t_0) & \mbox{in}\ \Omega.
\end{cases}
\end{equation}
It is easy to verify that
\begin{equation}\label{zuv}
    \tilde{y} = \tilde{u} + \tilde{v}; \quad \tilde{y}(x, 0) = \tilde{u}(x,0),\quad x\in \Omega.
\end{equation}
On one hand, in view of (\ref{deu}), we note that $\|\tilde{y}(\cdot,0)\|_{L^\infty(\Omega)}$ is bounded which indicates $\|\tilde{u}(\cdot,0)\|_{L^\infty(\Omega)} \leq M_0$ for some constant $M_0$. Moreover, since $\|\tilde{u}(\cdot,t)\|_{L^2(\Omega)}$ is logarithmically convex with respect to $t\in (0,T)$ (see, e.g.,~\cite{ChYa-NA08}), one has
\begin{equation}\label{estu}
    \|\tilde{u}(\cdot,t)\|_{L^2(\Omega)} \leq M_0^{\frac{t_0-t}{t_0}} \|\tilde{u}(\cdot,t_0)\|_{L^2(\Omega)}^{\frac{t}{t_0}} , \quad 0<t<t_0.
\end{equation}
On the other hand, in view of (\ref{dev}), it follows from the regularity of the parabolic equation (see e.g.,~\cite{Evans-Book10}) that
\begin{equation*}
    \max_{0\leq t\leq t_0} \|\tilde{v}(\cdot, t)\|_{L^2(\Omega)} \leq C \left(\|\alpha\|_{L^2(\Omega)}+\|\beta\|_{L^2(\Omega)} + \|y\|_{L^2(0, \frac{2}{3}T, H^2(\Omega))}\right).
\end{equation*}
By (\ref{Hy}), we obtain
\begin{equation}\label{estv}
    \max_{0\leq t\leq t_0} \|\tilde{v}(\cdot, t)\|_{L^2(\Omega)} \leq C \left(\|\alpha\|_{L^2(\Omega)}+\|\beta\|_{L^2(\Omega)} + \|y\|_{H^1(0, \frac{2}{3}T, L^2(\omega))}\right)
\end{equation}
By combining estimates (\ref{zuv}), (\ref{estu}), and (\ref{estv}), we arrive at
\begin{equation*}
    \begin{aligned}
        \|\tilde{y}(\cdot, t)\|_{L^2(\Omega)}&\leq \|\tilde{u}(\cdot, t)\|_{L^2(\Omega)} +\|\tilde{v}(\cdot, t)\|_{L^2(\Omega)}\\
        &\leq C \left(\|\tilde{u}(\cdot,t_0)\|_{L^2(\Omega)}^{\frac{t}{t_0}}+\|\alpha\|_{L^2(\Omega)}+\|\beta\|_{L^2(\Omega)} + \|y\|_{H^1(0, \frac{2}{3}T, L^2(\omega))}\right)\\
        &\leq C\left(\left(\|r\|_{L^2(\Omega)}+\|\tilde{v}(\cdot, t_0)\|_{L^2(\Omega)}\right)^{\frac{t}{t_0}} + \|\alpha\|_{L^2(\Omega)}+\|\beta\|_{L^2(\Omega)} + \|y\|_{H^1(0, \frac{2}{3}T, L^2(\omega))}\right)\\
        &\leq C\left(\left(\|y(t_0,\cdot)\|_{H^2(\Omega)}+\|\alpha\|_{L^2(\Omega)}+\|\beta\|_{L^2(\Omega)}\right)^{\frac{t}{t_0}} + \|\alpha\|_{L^2(\Omega)}+\|\beta\|_{L^2(\Omega)} + \|y\|_{H^1(0, \frac{2}{3}T, L^2(\omega))}\right)
    \end{aligned}
\end{equation*}
for $0<t<t_0$. 

For simplicity, we set
\begin{equation*}
    \mathfrak{a} := \|y(t_0,\cdot)\|_{H^2(\Omega)} + \|\alpha\|_{L^2(\Omega)} + \|\beta\|_{L^2(\Omega)}.
\end{equation*}
We recall $\mathfrak{a} \leq C\sqrt{G(\rho,\tilde\rho)}$ in~\eqref{ineq:alpha_and_beta}. In addition, in the small-data regime relevant for stability, we may assume $0 < \mathfrak{a} < 1$. Integrating the bound for $\tilde y$ obtained from~\eqref{estu}--\eqref{estv} over $(0,t_0)$ gives \begin{equation}\label{nufin}
\begin{aligned}
        \|\rho_0 - \tilde\rho_0\|_{L^2(\Omega)} &= \|y(\cdot, 0)\|_{L^2(\Omega)} = \left\|-\int_0^{t_0} \tilde{y}(\cdot, s) \, ds + y(\cdot, t_0)\right\|_{L^2(\Omega)}\\
        &\leq C \int_0^{t_0} \mathfrak{a}^{\frac{t}{t_0}} \, dt + Ct_0\left(\|\alpha\|_{L^2(\Omega)}+\|\beta\|_{L^2(\Omega)} + \|y\|_{H^1(0, T, L^2(\omega))}\right) + \|y(t_0,\cdot)\|_{L^2(\Omega)}.
\end{aligned}
\end{equation}
We now bound each term on the right-hand side of~\eqref{nufin}.

A direct computation gives
\begin{equation*}
    \int_0^{t_0} \mathfrak{a}^{\frac{t}{t_0}}\,dt \;=\; \frac{t_0\,(\mathfrak{a}-1)}{\log \mathfrak{a}} \;\leq\; \frac{Ct_0}{\bigl|\log \mathfrak{a}\bigr|}.
\end{equation*}
Applying the elementary inequality $-x\log x \leq 1/e$ on $(0,1)$ to $x = \mathfrak{a}$ yields $\mathfrak{a} \leq C/\bigl|\log \mathfrak{a}\bigr|$. Since each of $\|\alpha\|_{L^2(\Omega)}$, $\|\beta\|_{L^2(\Omega)}$, and $\|y(t_0,\cdot)\|_{L^2(\Omega)}$ is at most $\mathfrak{a}$, all three are likewise bounded by $C/\bigl|\log \mathfrak{a}\bigr|$.

The remaining data term $\|y\|_{H^1(0,T;L^2(\omega))}$ is not controlled by $\mathfrak{a}$, so we estimate it directly in terms of $G(\rho,\tilde\rho)$. Since the $H^2$-norm in time dominates the $H^1$-norm, the definition of $G(\rho,\tilde\rho)$ gives
\begin{equation*}
    \|y\|_{H^1(0,T;L^2(\omega))} \;\leq\; \|y\|_{H^2(0,T;L^2(\omega))} \;\leq\; \sqrt{G(\rho,\tilde\rho)}.
\end{equation*}
Applying the same elementary inequality to $x = \sqrt{G(\rho,\tilde\rho)}$ then yields $\sqrt{G(\rho,\tilde\rho)} \leq C/\bigl|\log G(\rho,\tilde\rho)\bigr|$.

Substituting these bounds into~\eqref{nufin}, we obtain
\begin{equation}\label{nufin2}
    \|\rho_0 - \tilde\rho_0\|_{L^2(\Omega)} \;\leq\; \frac{C}{\bigl|\log \mathfrak{a}\bigr|} \,+\, \frac{C}{\bigl|\log G(\rho,\tilde\rho)\bigr|}.
\end{equation}
Finally, since $\mathfrak{a} \leq C\sqrt{G(\rho,\tilde\rho)}$, we have $\bigl|\log \mathfrak{a}\bigr| \geq \tfrac{1}{2}\bigl|\log G(\rho,\tilde\rho)\bigr|$ whenever $G(\rho,\tilde\rho)$ is sufficiently small, so both terms on the right-hand side of~\eqref{nufin2} are bounded by $C/\bigl|\log G(\rho,\tilde\rho)\bigr|$ after adjusting the constant $C$.

Combining this with~\eqref{ineq:alpha_and_beta}, we conclude
\begin{equation*}
  \|\rho_0 - \tilde\rho_0\|_{L^2(\Omega)} \leq C \,\bigl|\log (\|\rho(\cdot, t_{0})-\tilde{\rho}(\cdot, t_{0})\|_{H^{2}(\Omega)}+\|\rho(\cdot, t_{1})-\tilde{\rho}(\cdot, t_{1})\|_{H^{2}(\Omega)}
    + \|\rho-\tilde{\rho}\|_{H^{1}(0, T ; L^{2}(\omega))})\bigr|^{-1}
\end{equation*}
for some constant $C$. This completes the proof of Theorem \ref{stabi}.

\section{A two-stage reconstruction algorithm}\label{SEC:algo}

We now develop a two-stage numerical reconstruction algorithm for simultaneously recovering the coefficients in the model~\eqref{EQ:Diff} with reaction nonlinearity~\eqref{EQ:Nonlinearity}, that is,
\begin{equation}\label{EQ:Diff KPP}
    \begin{cases}
            \partial_t \rho -\nabla\cdot (\gamma(\rho) \nabla \rho) + \rho (\mu-\xi\rho) = 0, & \mbox{in}\ \Omega\times (0,T),\\
            \rho(x,t) = b(x,t), & \mbox{on}\ \partial\Omega \times (0,T),\\
            \rho(x, 0) = \rho_0(x), & \mbox{in}\ \Omega.
    \end{cases}
\end{equation}
In contrast to the analysis of Section~\ref{SEC:Unique}, the numerical study does not enforce the strict positivity ($\rho_0,\mu,\xi>r_0$ and $b>b_0$) or the non-degeneracy condition~\eqref{ndeg}. These were required only to establish the stability estimates. 

For the sake of generality, we assume that we have snapshot data at $J$ time instances $0<t_0<t_1<\cdots<t_{J-1}\le T$. The theory of the previous section requires $J\ge 2$. In the experiments below, we take $J\ge 3$. The measurement information is therefore:
\begin{equation}\label{EQ:Data Numer}
    \begin{aligned}
        &d_j:=\rho(x,t_j),\quad (x, t)\in\Omega\times\{t_j\}_{j=0}^{J-1},\\
        &d_\omega:=\rho(x,t),\quad (x,t)\in \omega\times (0,T).
    \end{aligned}
\end{equation}
For the numerical simulations in Section~\ref{SEC:Numer}, we observe that one can have reasonable reconstructions without the data in $\omega$. However, if one can indeed have data for large $\omega$, the data would significantly improve the quality of the reconstructions. We include this datum in the numerical algorithm. Specifically, $d_\omega$ is used to construct the initial guess for $(\mu,\xi)$ described below, while the fitting functionals~\eqref{optimize} and~\eqref{optimize_rho0} use only the snapshot data $\{d_j\}_{j=0}^{J-1}$.

\subsection{Constructing initial guess of \texorpdfstring{$(\mu,\xi)$}{}}

First, we exploit the interior data $d_\omega$ to construct an inexpensive initial estimate of $\mu$ and $\xi$ on $\omega$. Since $d_\omega$ is available on the whole time interval, $\partial_t\rho$ and $\nabla\!\cdot(\gamma(\rho)\nabla\rho)$ can be evaluated on $\omega$ by numerical differentiation, so we define 
\begin{equation*}
    \sigma(x,t) := \frac{-\partial_t\rho + \nabla\!\cdot(\gamma(\rho)\nabla\rho)}{\rho}.
\end{equation*}
The governing PDE~\eqref{EQ:Diff} then implies the pointwise linear relation $\sigma(x,t) = \mu(x) - \xi(x)\,\rho(x,t)$. At each fixed spatial point $x\in \omega$, we fit this linear relation in $\rho$ by least squares
\[
    \min_{\mu(x),\xi(x)}\; \mathbb{E}_t[(\sigma - \mu + \xi\rho)^2]\,,
\]
which leads to the closed-form estimates of the unknown parameters $\mu$ and $\xi$ in $\omega$:
\begin{equation*}
     \tilde{\xi}(x) = -\frac{\mathrm{Cov}_t(\rho,\sigma)}{\mathrm{Var}_t(\rho)}, 
    \qquad
    \tilde{\mu}(x) = \mathbb{E}_t[\sigma] + \tilde{\xi}(x)\,\mathbb{E}_t[\rho],\quad x\in\omega
\end{equation*}
where $\bbE_t$, ${\rm Var}_t$, and ${\rm Cov}_t$ denote the empirical mean, variance, and covariance taken over the sampled observation time, and $\sigma$ is evaluated from $d_\omega$ by numerical differentiation.

We then extend $(\tilde\mu,\tilde\xi)$ from the subdomain $\omega$ to the entire domain $\Omega$ to serve as the initial guess.
Specifically, we employ polynomial extensions defined by
\begin{equation}\label{EQ:Extension}
\begin{aligned}
\mu^{(0)}|_{\omega} &= \tilde{\mu}, \qquad
\mu^{(0)} = \arg\min_{P}
\int_\omega \bigl| P(x)-\tilde{\mu}(x) \bigr|^2 dx
+ \lambda_P \|P\|_2^2,\\
\xi^{(0)}|_{\omega} &= \tilde{\xi}, \qquad
\xi^{(0)} = \arg\min_{P}
\int_\omega \bigl| P(x)-\tilde{\xi}(x) \bigr|^2 dx
+ \lambda_P \|P\|_2^2,
\end{aligned}
\end{equation}
where $P$ ranges over $d$-variate polynomials of fixed degree and $\lambda_P>0$ is a regularization parameter, and $\|\cdot\|_2$ denotes the Euclidean norm of the polynomial's coefficient vector. We remark that the specific choice of polynomial extension is not mandatory; any other smooth continuation method that extends the coefficients from $\omega$ to the entire domain $\Omega$ can be employed as a valid initial guess. Since forming $\sigma$ requires differentiating the data $d_\omega$, an operation that amplifies measurement noise, this closed-form estimate is used only to initialize the optimization; the regularized fit~\eqref{EQ:Extension} provides additional smoothing.

\subsection{First stage: reconstruction of the coefficients}

In the first stage of the reconstruction, we aim to recover the reaction coefficients $\mu$ and $\xi$. Since the initial condition of $\rho$ is unknown, we shift the temporal origin to the first observation time $t_0$. Specifically, we define:
\begin{equation*}
    \tilde{\rho}(x,t) = \rho(x, t + t_0).
\end{equation*}
Then $\tilde \rho$ has known initial condition
$\tilde{\rho}(x,0) = \rho(x,t_0)=d_0$ and solves
\begin{equation*}
    \begin{cases}
            \partial_t \tilde\rho -\nabla\cdot (\gamma(\tilde \rho) \nabla \tilde\rho) + \tilde\rho (\mu-\xi\tilde\rho) = 0, & \mbox{in}\ \Omega\times (0,T-t_0),\\
            \tilde \rho(x,t) = 0, & \mbox{on}\ \partial\Omega \times (0,T-t_0),\\
            \tilde \rho(x, 0) = d_0, & \mbox{in}\ \Omega.
    \end{cases}
\end{equation*}
This approach enables the independent reconstruction of the reaction coefficients $\mu$ and $\xi$ without requiring knowledge of the original, unknown initial state of $\rho$. 
We solve for $(\mu, \xi)$ by minimizing the least-squares functional
\begin{equation}\label{optimize}
    \Phi(\mu,\xi) :=\frac12\sum_{j=1}^{J-1}\!\int_\Omega(\tilde \rho(\cdot,t_{j}-t_0)-d_j)^2\,dx
+\;\frac{\alpha_\mu}{2}\|\nabla\mu\|_{L^2(\Omega)}^2+\frac{\alpha_\xi}{2}\|\nabla\xi\|_{L^2(\Omega)}^2.
\end{equation} 
Here the parameters $\alpha_\mu,\alpha_\xi>0$ are $H^1$ semi-norm regularization weights. 

To differentiate $\Phi$, we regard the parameter-to-state map $(\mu,\xi)\mapsto\tilde\rho(\mu,\xi)$ as a mapping from $C^{2+\alpha}(\overline\Omega)\times C^{2+\alpha}(\overline\Omega)$ into $C^{4+\alpha,2+\alpha/2}(\overline\Omega\times[0,T-t_0])$, which is Fr\'echet differentiable by the well-posedness theory of Section~\ref{well-posedness}. For admissible perturbations $(\delta\mu,\delta\xi)$, the directional derivative $\tilde\rho'=D\tilde\rho(\mu,\xi)[\delta\mu,\delta\xi]$ is the unique solution of the linearized (sensitivity) problem obtained by differentiating the forward model,
\begin{equation}\label{EQ:Sensitivity}
    \begin{cases}
        \partial_t\tilde\rho' -\nabla\cdot(\gamma(\tilde\rho)\nabla\tilde\rho') -\nabla\cdot\big(\gamma'(\tilde\rho)\tilde\rho'\,\nabla\tilde\rho\big) +(\mu-2\xi\tilde\rho)\tilde\rho' = -(\delta\mu-\tilde\rho\,\delta\xi)\,\tilde\rho, & \mbox{in}\ \Omega\times(0,T-t_0),\\
        \tilde\rho'(x,t)=0, & \mbox{on}\ \partial\Omega\times(0,T-t_0),\\
        \tilde\rho'(x,0)=0, & \mbox{in}\ \Omega.
    \end{cases}
\end{equation}
By the chain rule, $\Phi$ is then Fr\'echet differentiable as the composition of this differentiable map with the smooth quadratic misfit and the $H^1$ semi-norm penalty, and its derivative is the bounded linear functional
\begin{equation}\label{EQ:Frechet Raw}
    \Phi'(\mu,\xi)[\delta\mu,\delta\xi]
    = \sum_{j=1}^{J-1}\int_\Omega \big(\tilde\rho(\cdot,t_j-t_0)-d_j\big)\,\tilde\rho'(\cdot,t_j-t_0)\,dx
    + \alpha_\mu\!\int_\Omega \nabla\mu\cdot\nabla\delta\mu\,dx
    + \alpha_\xi\!\int_\Omega \nabla\xi\cdot\nabla\delta\xi\,dx.
\end{equation}
To remove the implicit dependence on the sensitivity $\tilde\rho'$, we introduce the adjoint state $p$ solving
\begin{equation}\label{EQ:Diff KPP Adjoint}
    \begin{cases}
            -\partial_t p -\nabla\cdot (\gamma(\tilde \rho) \nabla p)+\gamma'(\tilde \rho)\nabla{\tilde\rho}\cdot\nabla p + (\mu-2\xi\tilde\rho)p = S, & \mbox{in}\ \Omega\times (0,T-t_0),\\
            p(x,t) = 0, & \mbox{on}\ \partial\Omega \times (0,T-t_0),\\
            p(x, T-t_0) =0, & \mbox{in}\ \Omega,
    \end{cases}
\end{equation}
where the source for the adjoint equation is
\[
S(x,t)=
- \sum_{j=1}^{J-1}\delta(t - (t_{j}-t_0))\,\big(\tilde \rho(x,t_j-t_0) - d_j\big)\,.
\]
This equation runs backward in time, from $t=T-t_0$ down to $t=0$. Multiplying the sensitivity equation~\eqref{EQ:Sensitivity} by $p$, integrating over $\Omega\times(0,T-t_0)$, and integrating by parts in space and time, 
using the homogeneous boundary and terminal conditions on $p$ together with $\tilde\rho'(\cdot,0)=0$, eliminates $\tilde\rho'$ from~\eqref{EQ:Frechet Raw} and yields the representation
\begin{equation}\label{EQ:Frechet Stage1}
    \Phi'(\mu,\xi)[\delta\mu,\delta\xi] = \int_\Omega g_\mu\,\delta\mu\,dx + \int_\Omega g_\xi\,\delta\xi\,dx.
\end{equation}
Here, assuming the natural homogeneous Neumann condition $\partial_n\mu=\partial_n\xi=0$ on $\partial\Omega$ associated with the semi-norm penalty (so that the boundary terms from the integration by parts of the regularization vanish), the $L^2(\Omega)$ gradients are
\begin{equation}\label{EQ:Grad Stage1}
    g_\mu = \int_{0}^{T-t_0}\!\tilde\rho(x,t)\,p(x,t)\,dt - \alpha_\mu\Delta\mu,
    \qquad
    g_\xi = -\int_{0}^{T-t_0}\!\tilde\rho(x,t)^2\,p(x,t)\,dt - \alpha_\xi\Delta\xi.
\end{equation}
By construction, $g_\mu$ and $g_\xi$ are the Riesz representatives of the Fr\'echet derivative $\Phi'(\mu,\xi)$ in $L^2(\Omega)$, i.e.\ the gradients used in the optimization. Therefore, to evaluate the Fr\'echet derivative, we need one forward PDE solve for the state variable $\tilde\rho$, followed by one backward adjoint PDE solve for $p$.

\subsection{Second stage: reconstruction of the initial condition}

To reconstruct the unknown initial condition, we undo the time shift introduced in the first stage and treat the snapshots collected after $t=0$ as the data for the inversion. Specifically, after obtaining $(\mu,\xi)$ from the first stage, we solve the forward model on the full interval $[0,T]$ with the coefficients held fixed and treat $\rho(\cdot,0)=\rho_0$ as the control variable. The earliest snapshot $d_0=\rho(\cdot,t_0)$ is no longer imposed as a pseudo-initial condition. Instead, all measured snapshots $d_j$ at times $t_j$ ($0\le j\le J-1$) are used as fitting data. This time-shift-back strategy enables the identification of $\rho_0$ from information available strictly after $t=0$, while keeping the coefficient estimates fixed.

With $(\mu,\xi)$ fixed at their first stage reconstructions, we estimate the initial condition $\rho_0$ by solving the PDE-constrained optimization problem with the objective function:
\begin{equation}\label{optimize_rho0}
    \Psi(\rho_0):=
\frac{1}{2}\sum_{j=0}^{J-1}\int_\Omega \big(\rho(x,t_{j})-d_j\big)^2\,dx
+\frac{\alpha_\rho}{2} \|\nabla\rho_0\|_{L^2(\Omega)}^2\,,
\end{equation}
where $\alpha_\rho>0$ is the $H^1$ semi-norm regularization parameter, and $\rho=\rho(\rho_0)$ solves the forward model~\eqref{EQ:Diff KPP} on $[0,T]$ with the fixed coefficients $(\mu,\xi)$ and initial value $\rho(\cdot,0)=\rho_0$.

Due to the fact that the forward PDE~\eqref{EQ:Diff} is nonlinear, the inverse initial condition problem is also nonlinear. The objective function is therefore nonconvex. Differentiating the forward model with respect to $\rho_0$ in a direction $\delta\rho_0$, the sensitivity $\rho'=D\rho(\rho_0)[\delta\rho_0]$ solves the linearized problem
\begin{equation}\label{EQ:Sensitivity-2}
    \begin{cases}
        \partial_t\rho' -\nabla\cdot(\gamma(\rho)\nabla\rho') -\nabla\cdot\big(\gamma'(\rho)\rho'\,\nabla\rho\big) +(\mu-2\xi\rho)\rho' = 0, & \mbox{in}\ \Omega\times(0,T),\\
        \rho'(x,t)=0, & \mbox{on}\ \partial\Omega\times(0,T),\\
        \rho'(x,0)=\delta\rho_0, & \mbox{in}\ \Omega.
    \end{cases}
\end{equation}
Let $q$ be the solution to the adjoint problem
\begin{equation}\label{EQ:Diff KPP Adjoint-2}
    \begin{cases}
            -\partial_t q -\nabla\cdot (\gamma(\rho) \nabla q)+\gamma'(\rho)\nabla\rho\cdot\nabla q + (\mu-2\xi\rho)q = Q, & \mbox{in}\ \Omega\times (0,T),\\
            q(x,t) = 0, & \mbox{on}\ \partial\Omega \times (0,T),\\
            q(x, T) =0, & \mbox{in}\ \Omega,
    \end{cases}
\end{equation}
where the source for the adjoint equation is
\[
Q(x,t)=
- \sum_{j=0}^{J-1}\delta(t - t_{j})\,\big(\rho(x,t_j) - d_j\big)\,.
\]
Pairing the sensitivity equation~\eqref{EQ:Sensitivity-2} with $q$, integrating over $\Omega\times(0,T)$, and integrating by parts, using $q(\cdot,T)=0$, the homogeneous boundary conditions, and $\rho'(\cdot,0)=\delta\rho_0$, eliminates $\rho'$ and identifies the Fr\'echet derivative of $\Psi$ with respect to $\rho_0$ as the bounded linear functional
\begin{equation}\label{EQ:Frechet-2}
    \Psi'(\rho_0)[\delta\rho_0] = \int_\Omega \big(-q(x,0)-\alpha_\rho\,\Delta\rho_0\big)\,\delta\rho_0\,dx,
\end{equation}
where we have again imposed the natural homogeneous Neumann condition $\partial_n\rho_0=0$ on $\partial\Omega$ for the semi-norm penalty. Equivalently, the $L^2(\Omega)$ gradient is $\nabla_{\rho_0}\Psi=-q(\cdot,0)-\alpha_\rho\,\Delta\rho_0$, whose evaluation requires one forward solve for $\rho$ followed by one backward adjoint solve for $q$.
 
\section{Numerical experiments}\label{SEC:Numer}

We now present some numerical simulations based on the algorithms in Section \ref{SEC:algo}. For simplicity, we restrict attention to two-dimensional simulations. The extension to three dimensions is straightforward, subject only to an increase in computational cost. The computational domain is defined as the square $\Omega=(-1,1)^2$ over the time interval $(0,0.4)$, discretized using a uniform $201\times201$ grid. We implemented a quasi-Newton method with the BFGS updating rule for the Hessian matrix, following the implementation in~\cite{ReBaHi-SIAM06}.

We fix the diffusion coefficient to the constant value $\gamma\equiv1$; the algorithm applies without modification to a genuinely density-dependent $\gamma(\rho)$, at the cost of additional forward-solve expense. Measurements are acquired within the interior subdomain $\omega = (-0.5, 0.5)^2$ at the time snapshots $t \in \{0.02, 0.04, 0.06\}$, which play the roles of $t_0,t_1,t_2$. Throughout this work, spatial coordinates in $\mathbb{R}^2$ are denoted interchangeably by $x$ or $(x,y)$.

The functionals~\eqref{optimize} and~\eqref{optimize_rho0}  are minimized by the quasi-Newton iteration described above, terminating at step $k$ once
\begin{equation*}
\frac{\Phi^{k}-\Phi^{k+1}}{\max\{|\Phi^{k}|,|\Phi^{k+1}|,1\}}\le 10^{-20}
\qquad\text{or}\qquad
\|\nabla\Phi^{k}\|_{\infty}\le 10^{-8},
\end{equation*}
and analogously for $\Psi$, with the iteration count capped at $200$ for~\eqref{optimize} and $80$ for~\eqref{optimize_rho0}; in practice the latter converges in about $30$ iterations. We take $\alpha_\mu=\alpha_\xi=\alpha_\rho=10^{-7}$, so that the data misfit dominates and the $H^1$ penalties act only as mild smoothers, and we extend the interior estimate of $(\mu,\xi)$ to $\Omega$ through~\eqref{EQ:Extension} using degree-five tensor-product Legendre polynomials with $\lambda_P=10^{-7}$.

To quantify the quality of the reconstructions, we employ relative $L^2$ errors. Specifically, 
let $f$ denote a quantity to be reconstructed, with $f_t$ and $f_r$ representing its true and reconstructed values, then the relative $L^2$ reconstruction error is defined as
\begin{equation}
\mathcal{E}_{L^2}(f)
= \frac{\| f_t - f_r \|_{L^2(\Omega)}}{\| f_t \|_{L^2(\Omega)}}.
\end{equation}
Noisy data are generated by adding independent Gaussian noise to each measurement, with standard deviation chosen so that the noise-to-signal ratio (NSR) in the $L^2$ norm, $\|f_{\rm noisy}-f\|_{L^2}/\|f\|_{L^2}$, equals the prescribed level $\eta$.

In what follows, we present a selection of representative numerical results.

\begin{figure}[!htb]
\centering
\begin{subfigure}{0.24\textwidth}
  \centering
  \includegraphics[width=\linewidth]{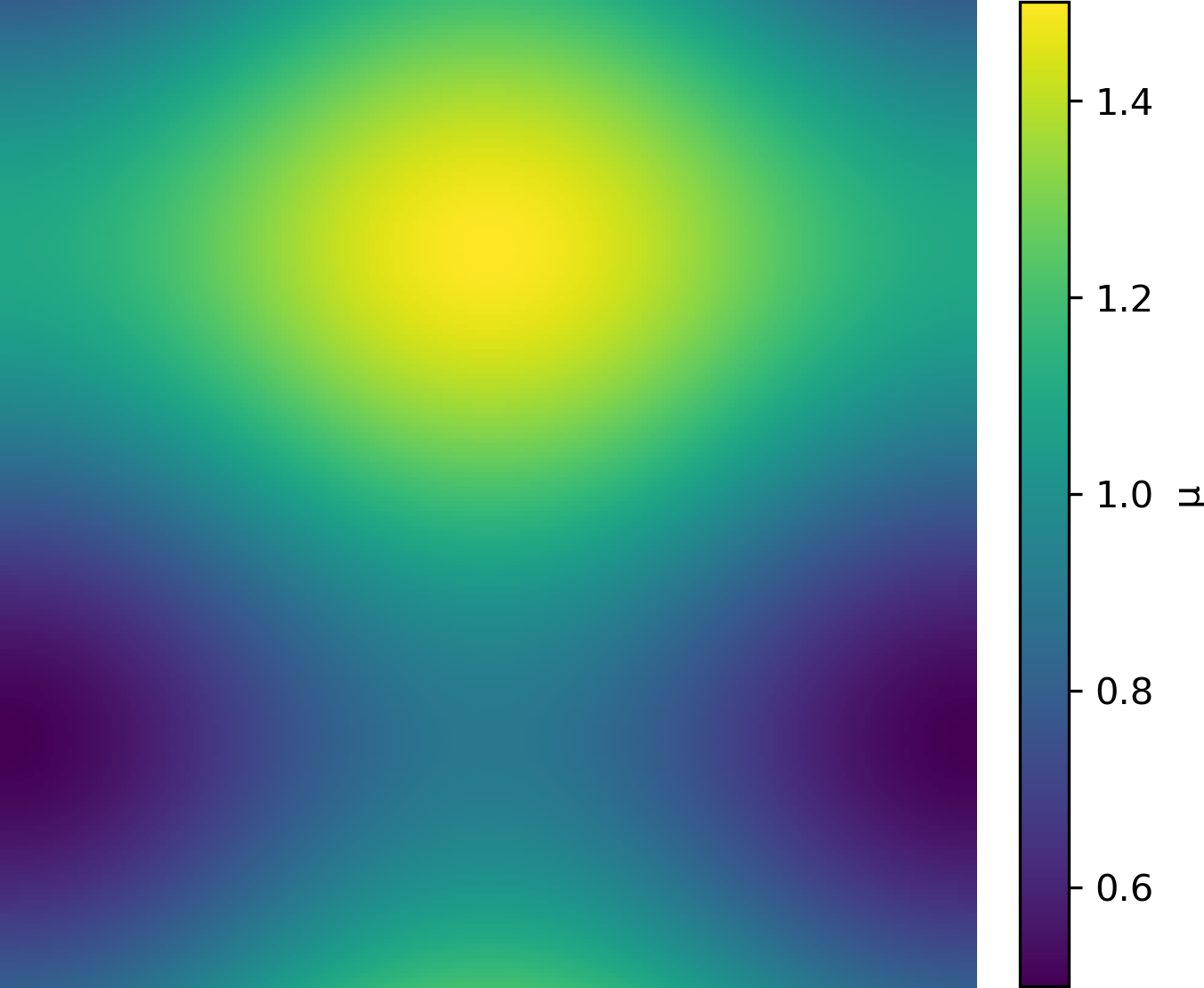}
  \caption{$\mu_{\mathrm{true}}$}
\end{subfigure}\hfill
\begin{subfigure}{0.24\textwidth}
  \centering
  \includegraphics[width=\linewidth]{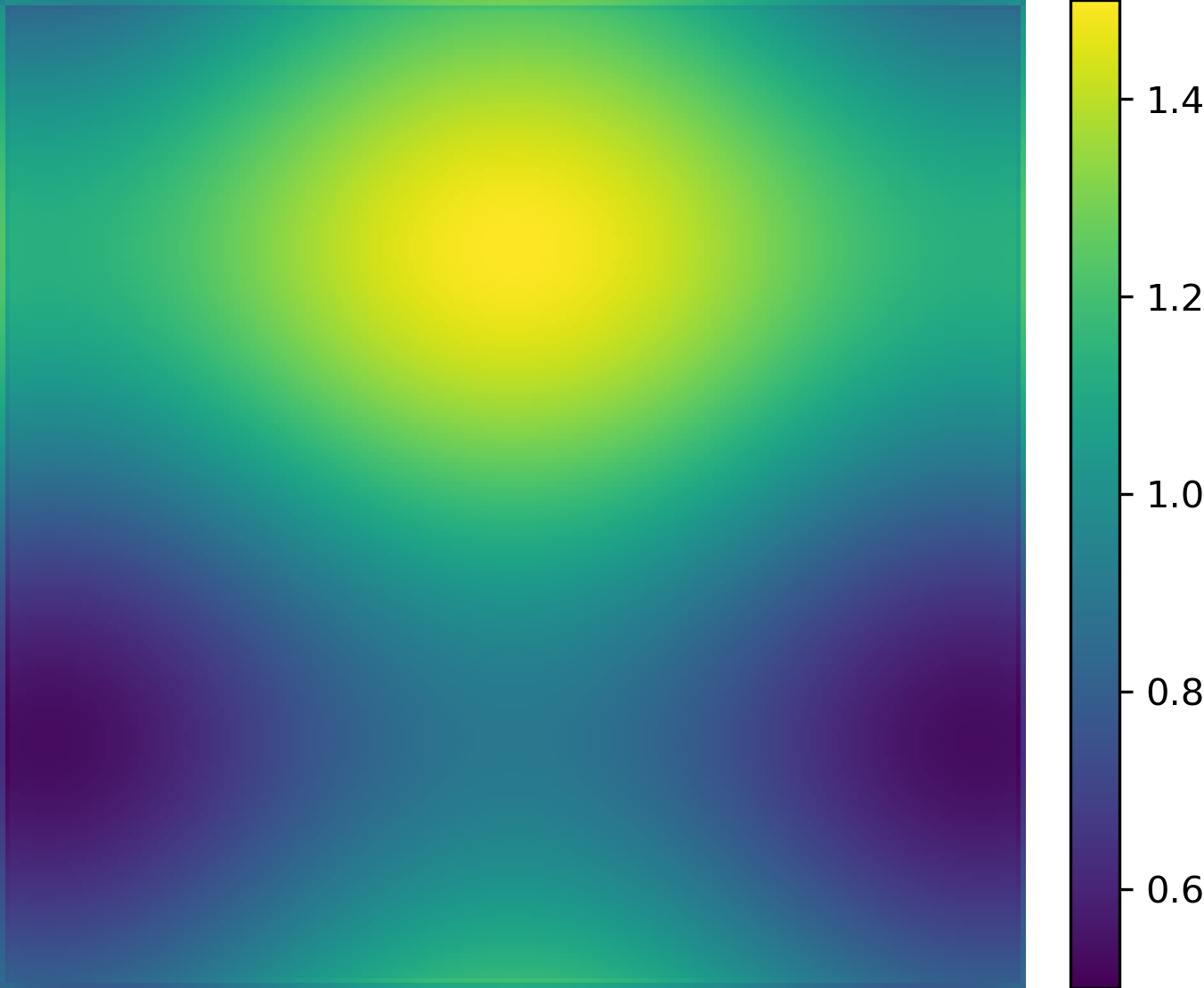}
  \caption{$\mu_{\mathrm{rec}}$}
\end{subfigure}\hfill
\begin{subfigure}{0.24\textwidth}
  \centering
  \includegraphics[width=\linewidth]{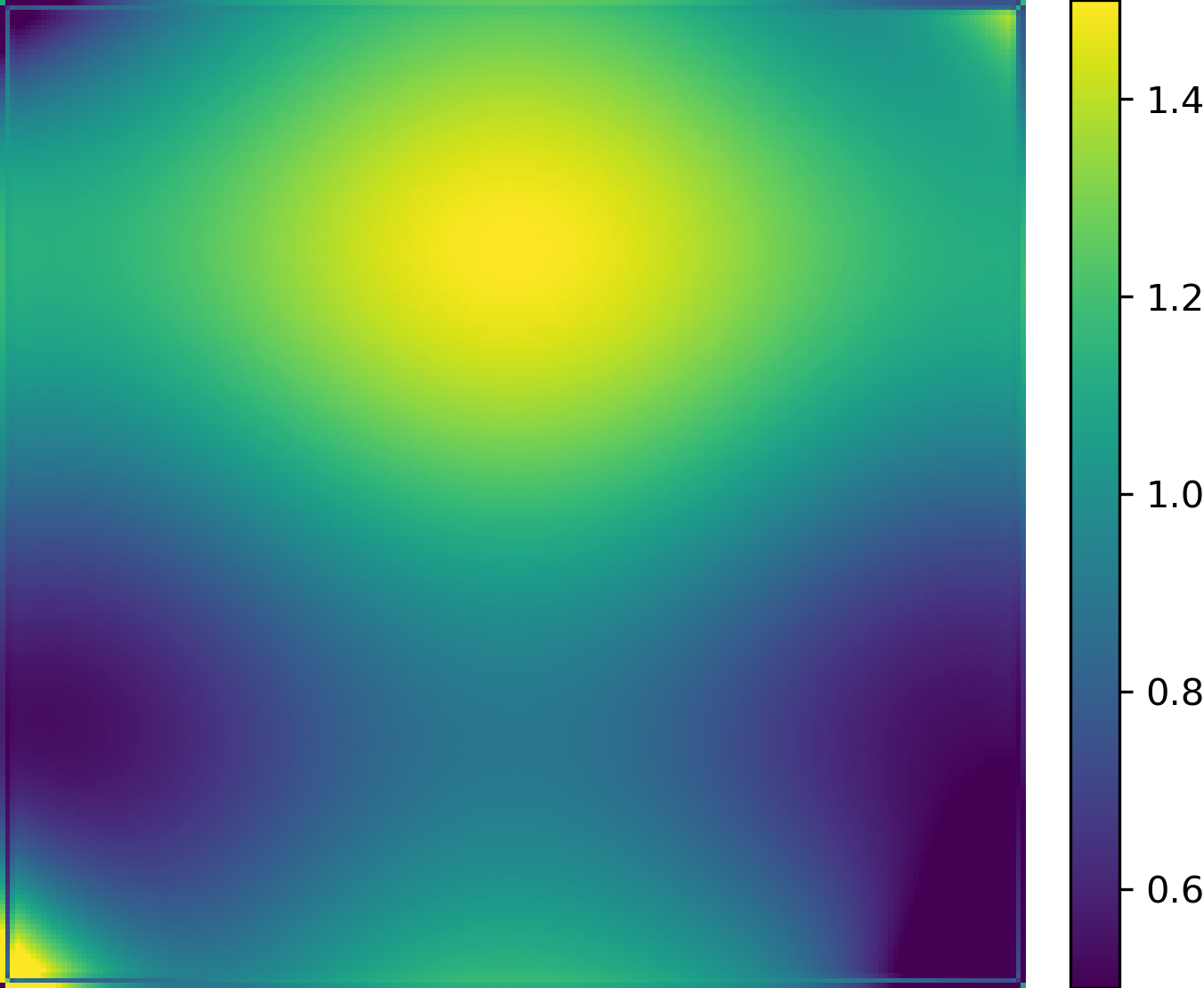}
  \caption{$\mu_{\mathrm{rec}}$ (1\% noise)}
\end{subfigure}\hfill
\begin{subfigure}{0.24\textwidth}
  \centering
  \includegraphics[width=\linewidth]{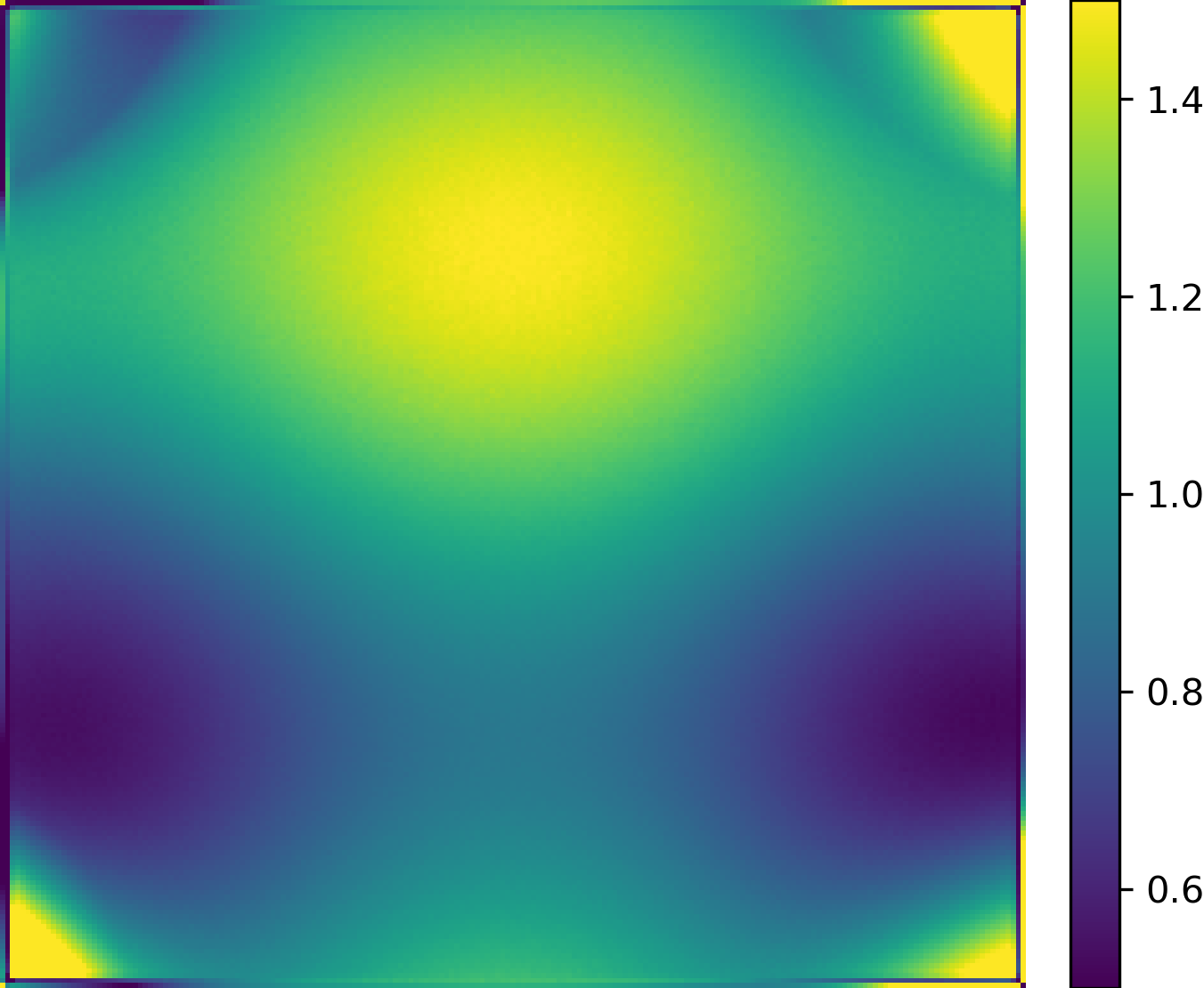}
  \caption{$\mu_{\mathrm{rec}}$ (5\% noise)}
\end{subfigure}

\vspace{2mm}

\begin{subfigure}{0.24\textwidth}
  \centering
  \includegraphics[width=\linewidth]{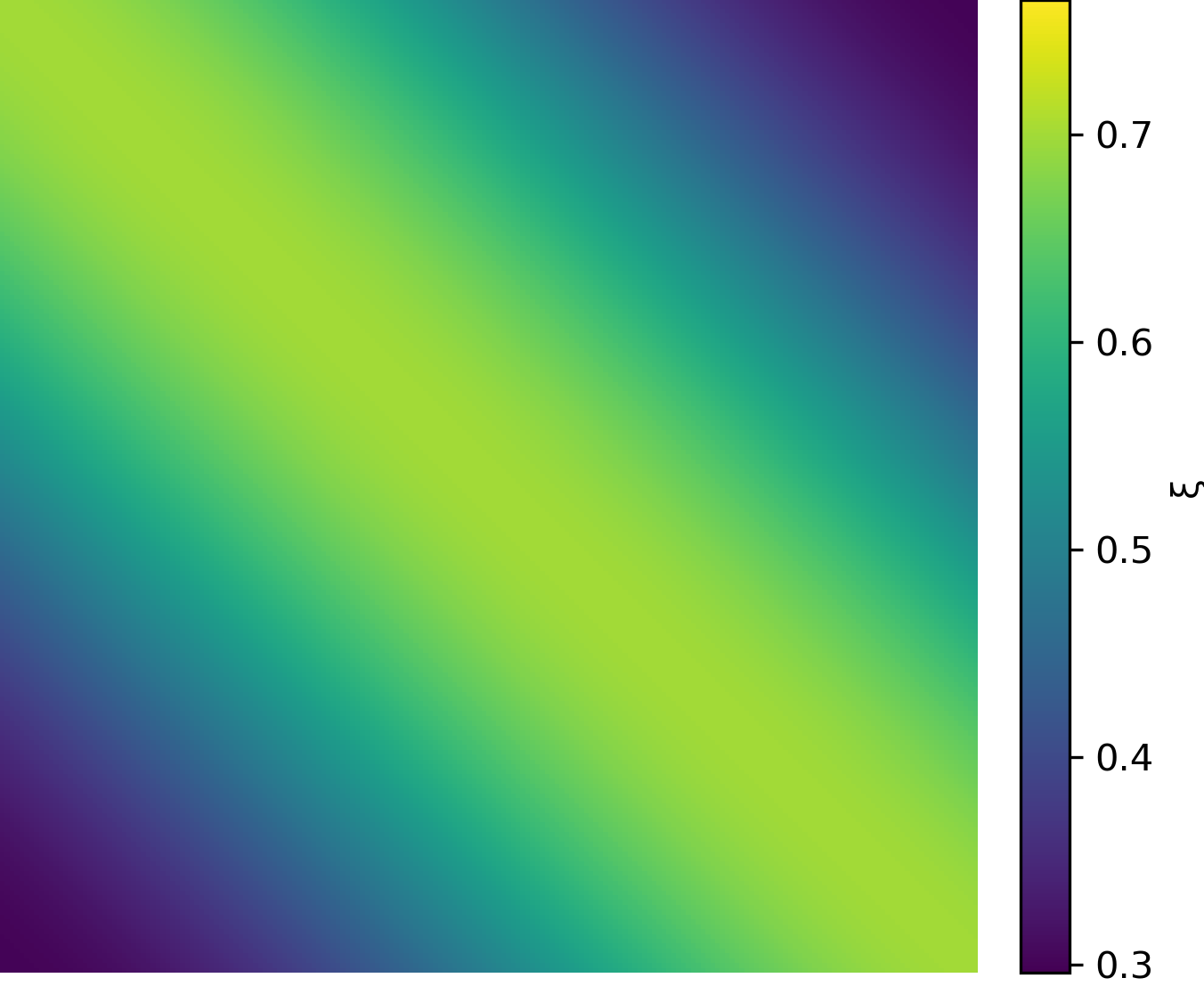}
  \caption{$\xi_{\mathrm{true}}$}
\end{subfigure}\hfill
\begin{subfigure}{0.24\textwidth}
  \centering
  \includegraphics[width=\linewidth]{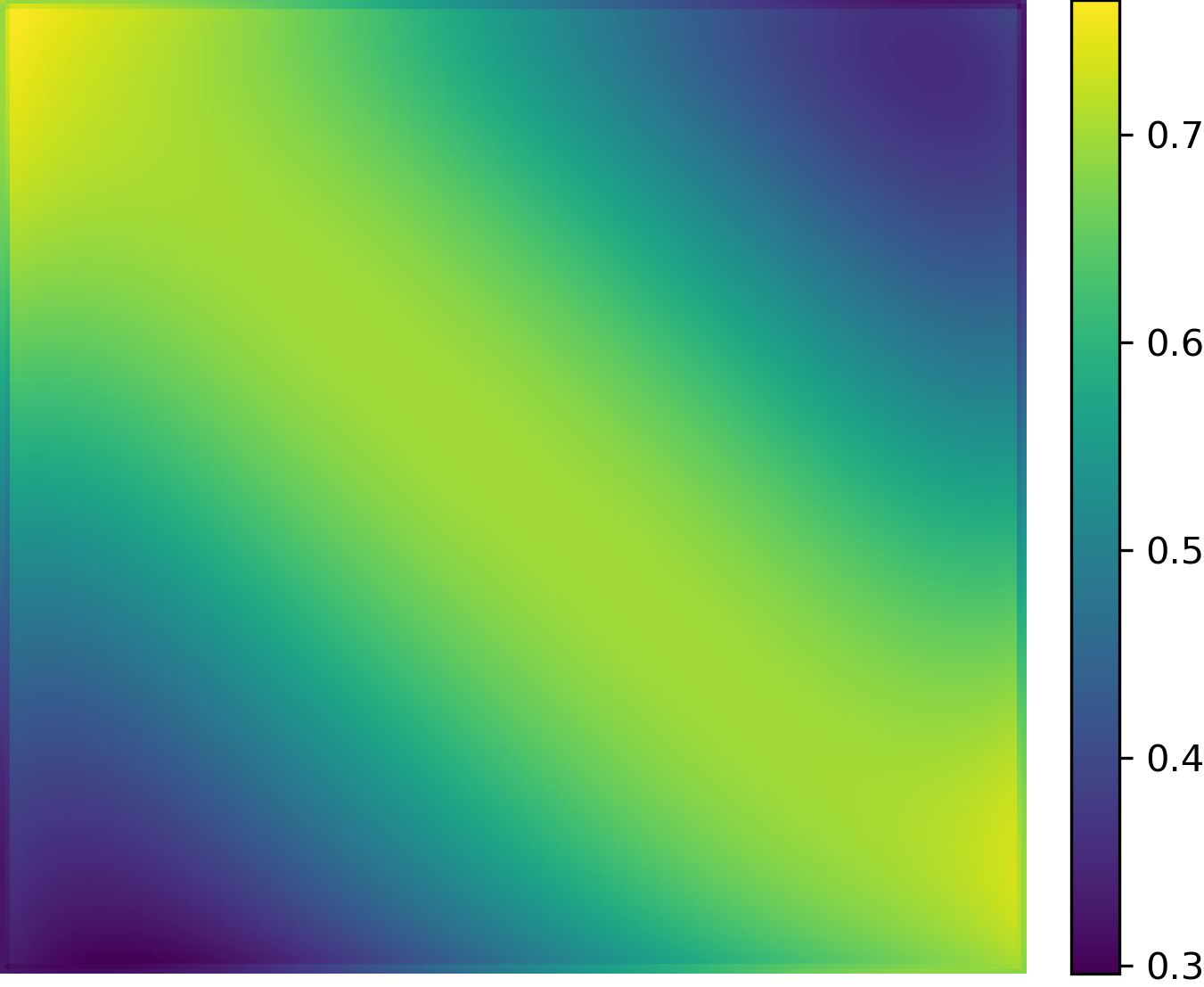}
  \caption{$\xi_{\mathrm{rec}}$}
\end{subfigure}\hfill
\begin{subfigure}{0.24\textwidth}
  \centering
  \includegraphics[width=\linewidth]{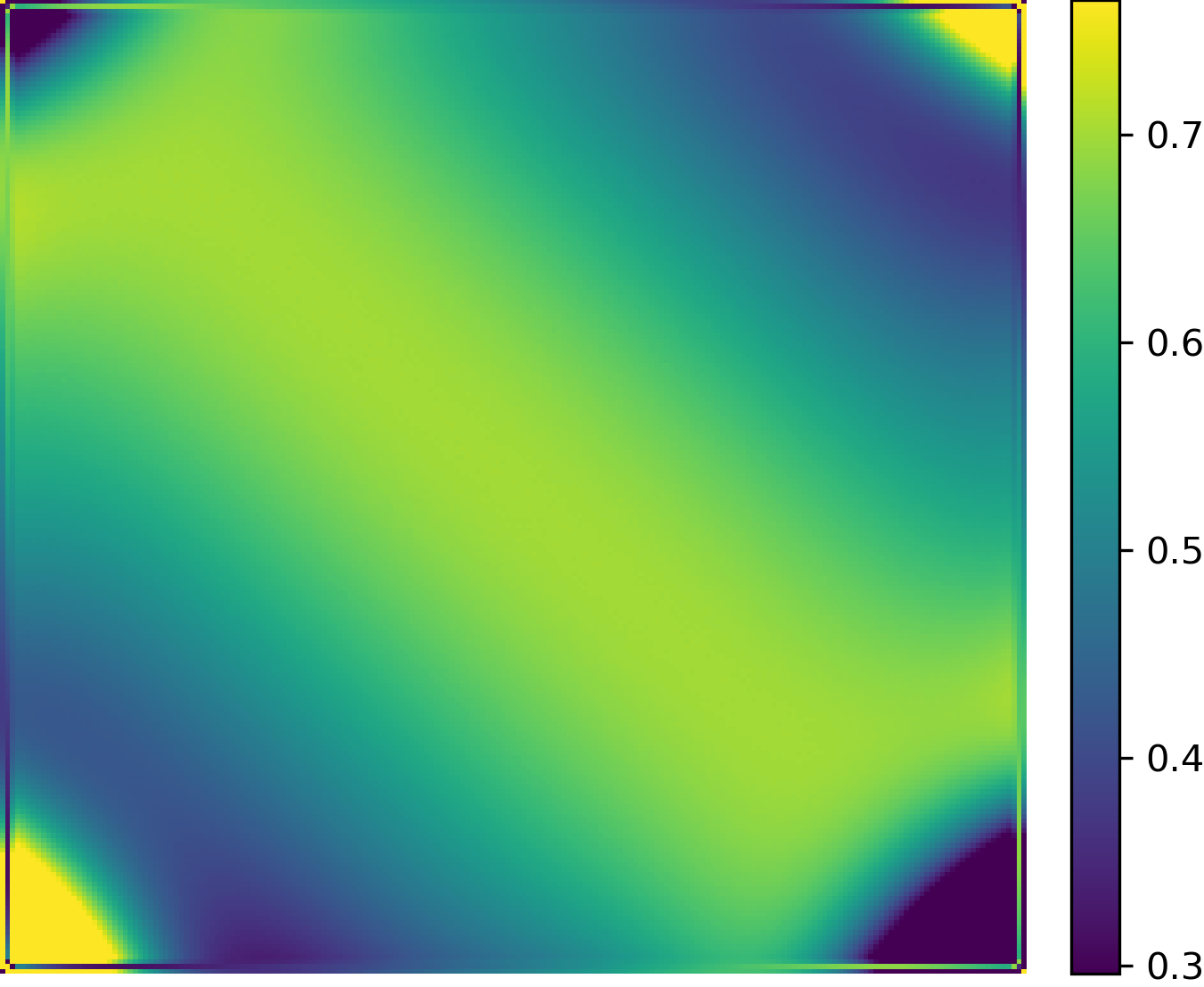}
  \caption{$\xi_{\mathrm{rec}}$ (1\% noise)}
\end{subfigure}\hfill
\begin{subfigure}{0.24\textwidth}
  \centering
  \includegraphics[width=\linewidth]{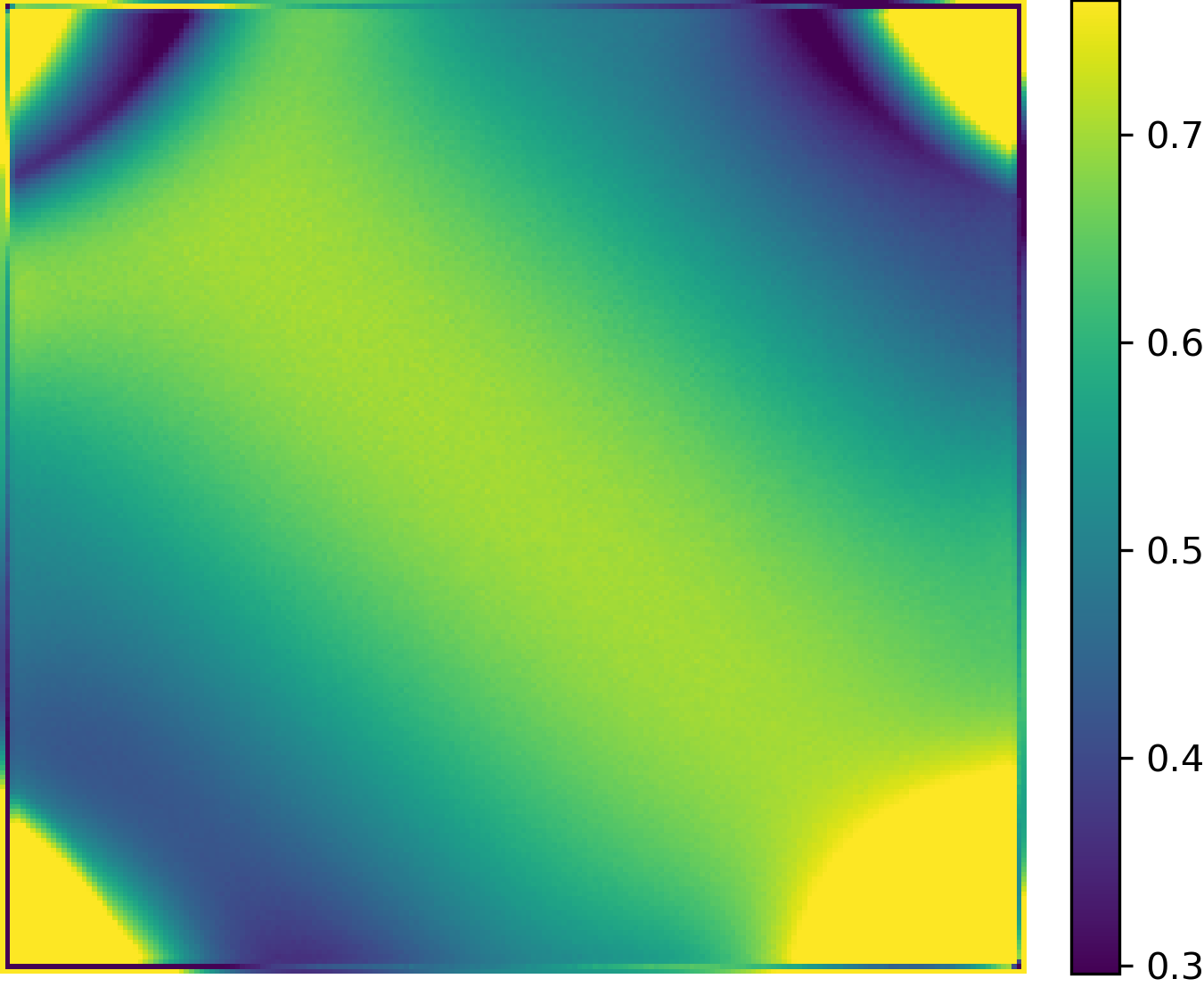}
  \caption{$\xi_{\mathrm{rec}}$ (5\% noise)}
\end{subfigure}

\caption{
Reconstruction results for the coefficients.
Top row: $\mu$. Bottom row: $\xi$.
From left to right: true coefficient, noiseless reconstruction,
reconstruction with 1\% noise, and reconstruction with 5\% noise.
}
\label{fig:coeff_recon_1}
\end{figure}

\paragraph{Experiment I.}
In the first numerical experiment, we focus on the simultaneous reconstruction of the reaction coefficients $(\mu,\xi)$ given by:
\begin{equation*}
\begin{aligned}
\mu(x,y) &= 1.0 + 0.3 \sin(\pi x) + 0.2 \cos(\pi y), \\
\xi(x,y) &= 0.5 + 0.2 \cos\bigl(0.5\pi(x+y)\bigr).
\end{aligned}
\end{equation*}
The reconstruction is performed with the initial condition set to $\rho_0(x,y) = 0.5 + 0.2 \cos(0.5\pi x) + 0.15 \cos(0.5\pi y).$ The recovered coefficients are visualized in Figure~\ref{fig:coeff_recon_1}, where the panels (from left to right) correspond to data with noise-to-signal ratios (NSR) of $\eta = 0$, $0.01$, and $0.05$. The method demonstrates robust performance, yielding relative $L^2$ errors for the pair $(\mu,\xi)$ of $(1.89,\,3.00)\times 10^{-2}$, $(7.79,\,24.15)\times 10^{-2}$, and $(2.21,\,4.55)\times 10^{-1}$, respectively.

\begin{figure}[!htb]
\centering
\begin{subfigure}{0.24\textwidth}
  \centering
  \includegraphics[width=\linewidth]{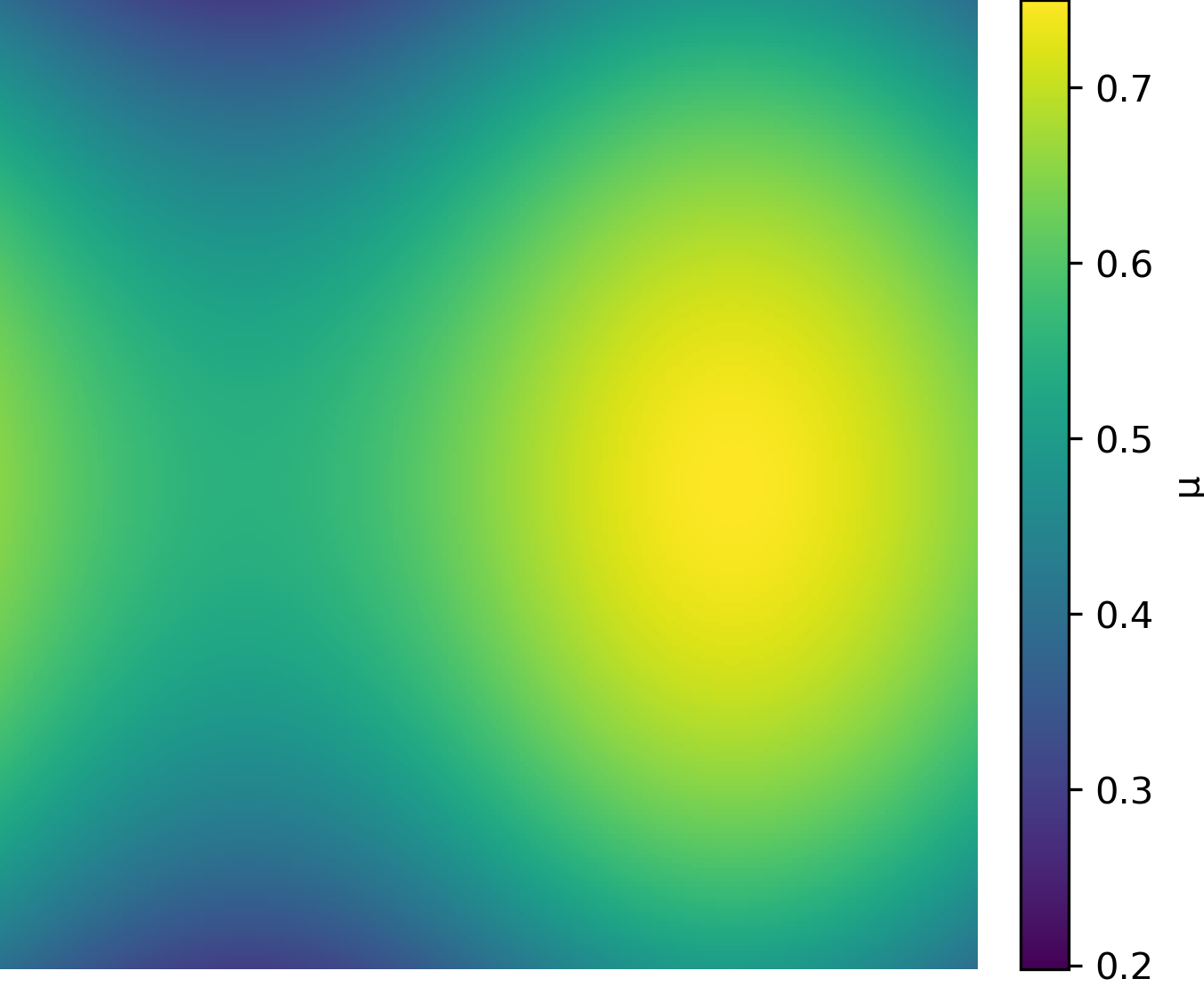}
  \caption{$\mu_{\mathrm{true}}$}
\end{subfigure}\hfill
\begin{subfigure}{0.24\textwidth}
  \centering
  \includegraphics[width=\linewidth]{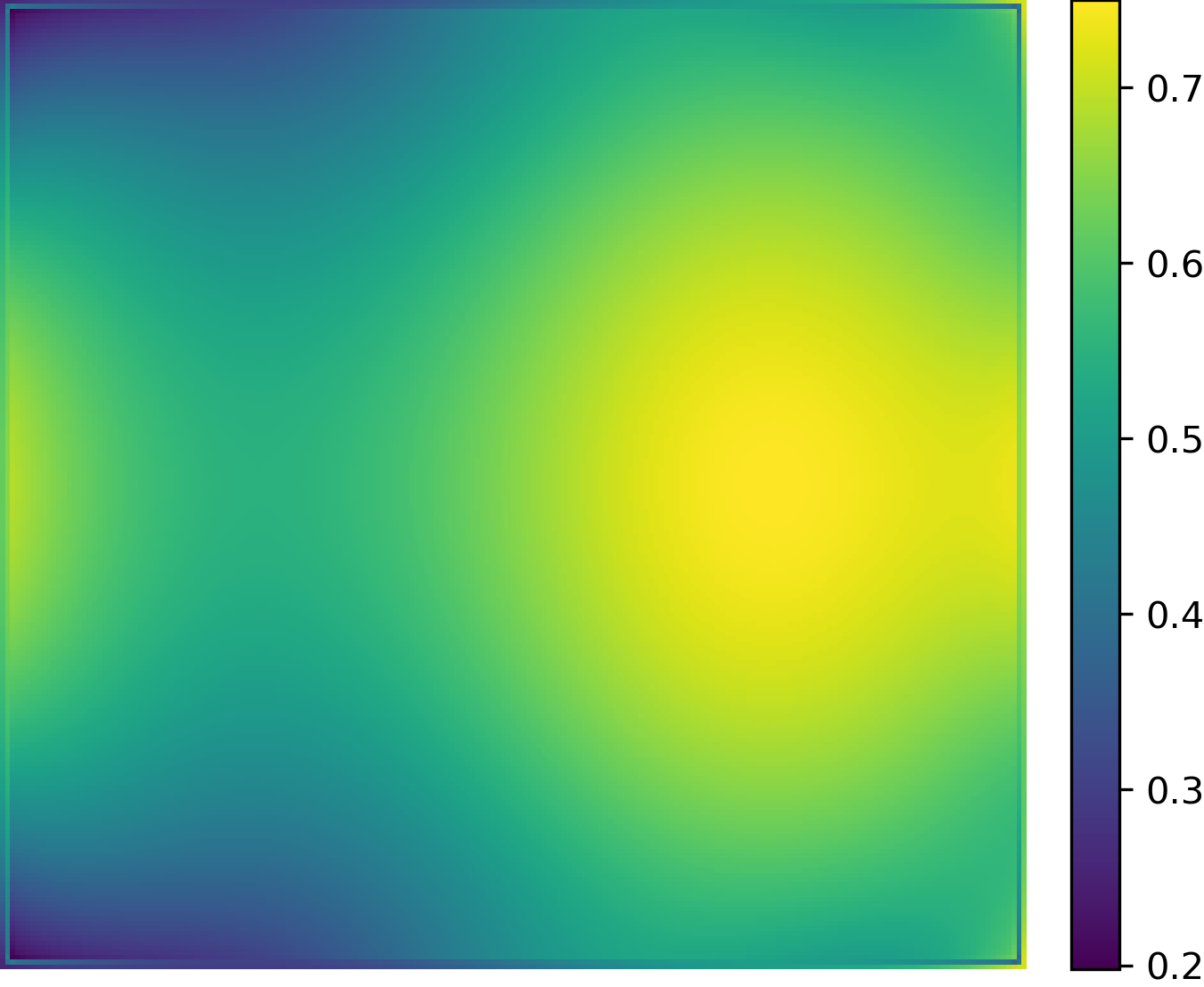}
  \caption{$\mu_{\mathrm{rec}}$}
\end{subfigure}\hfill
\begin{subfigure}{0.24\textwidth}
  \centering
  \includegraphics[width=\linewidth]{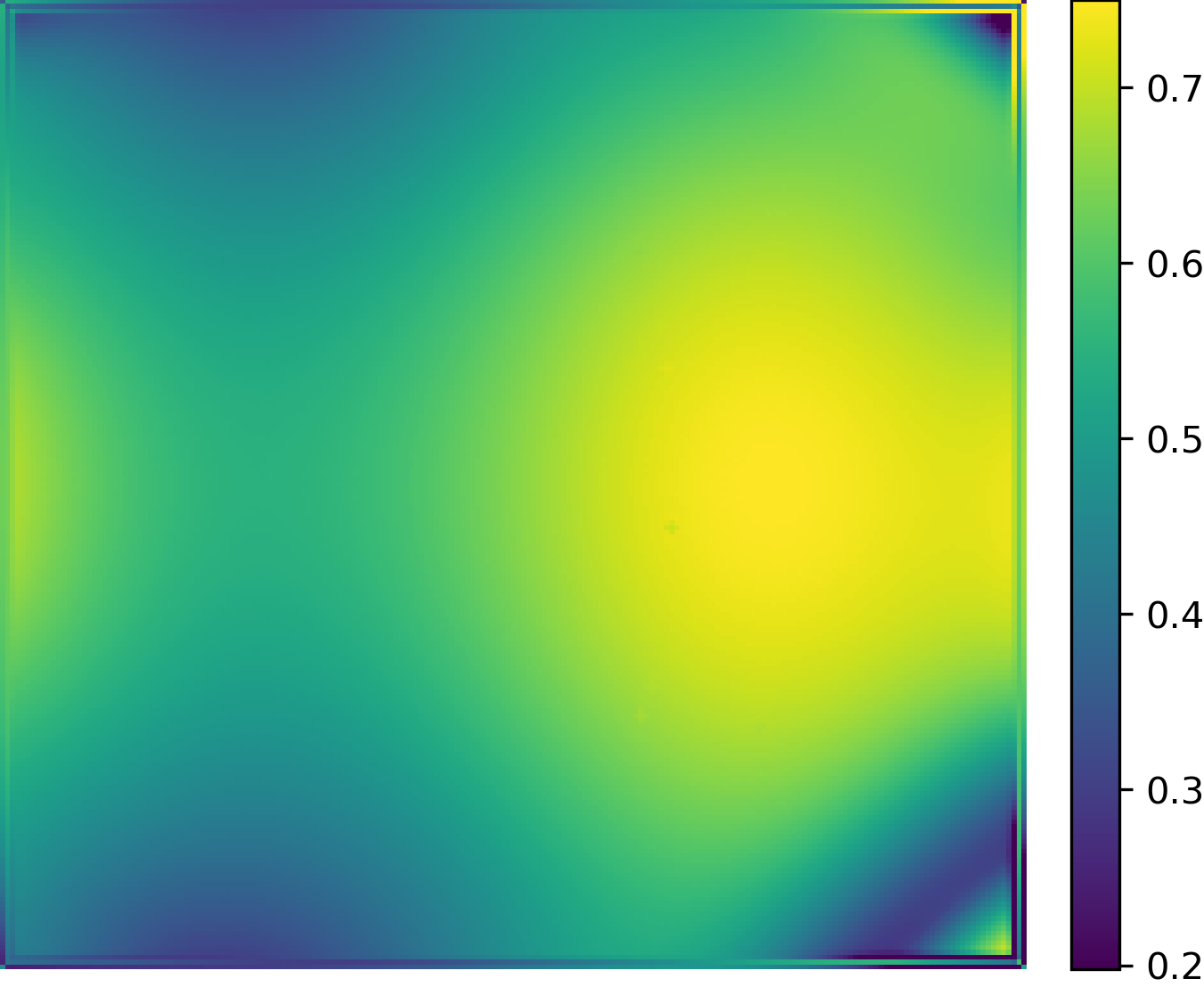}
  \caption{$\mu_{\mathrm{rec}}$ (1\% noise)}
\end{subfigure}\hfill
\begin{subfigure}{0.24\textwidth}
  \centering
  \includegraphics[width=\linewidth]{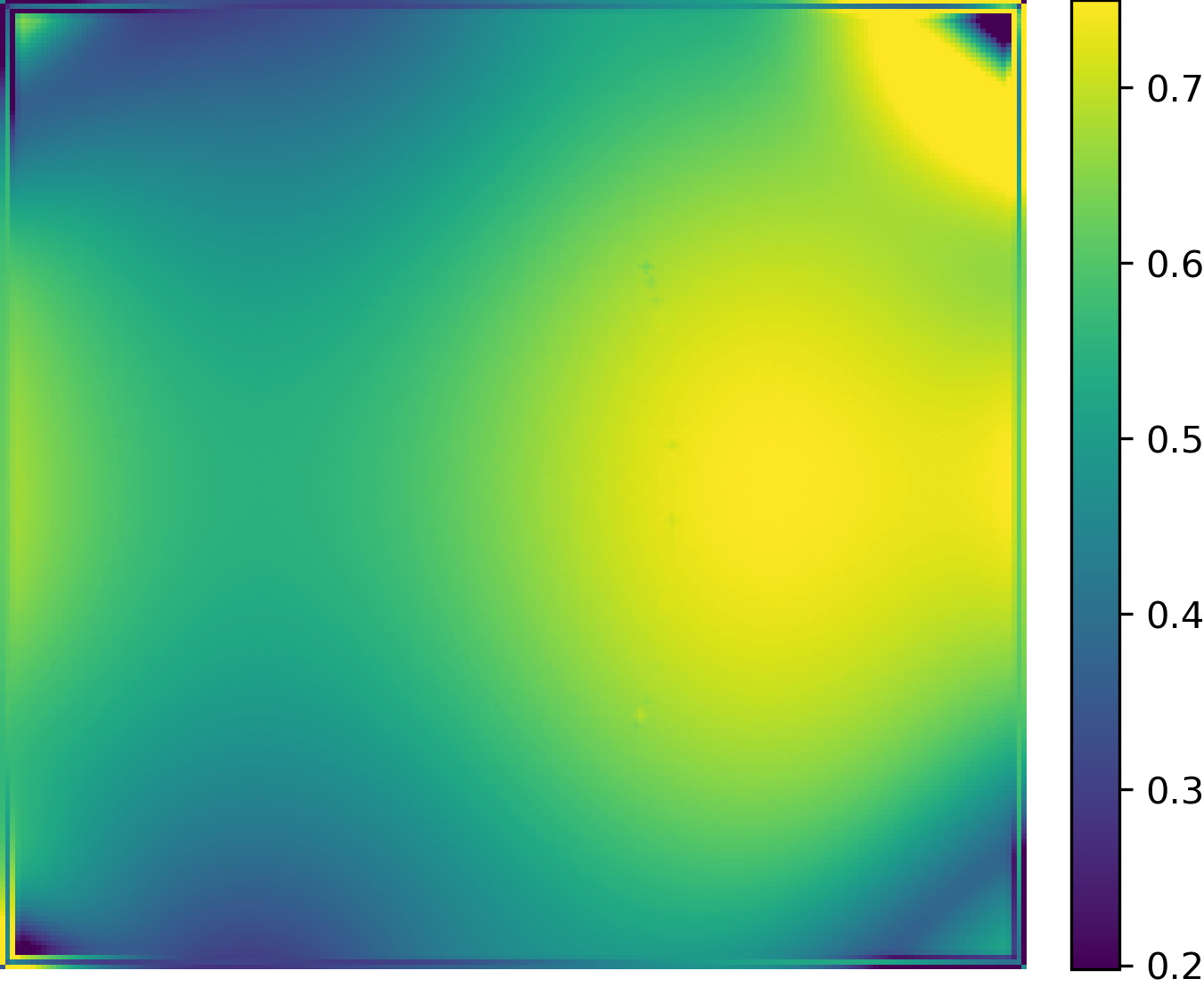}
  \caption{$\mu_{\mathrm{rec}}$ (5\% noise)}
\end{subfigure}

\vspace{2mm}

\begin{subfigure}{0.24\textwidth}
  \centering
  \includegraphics[width=\linewidth]{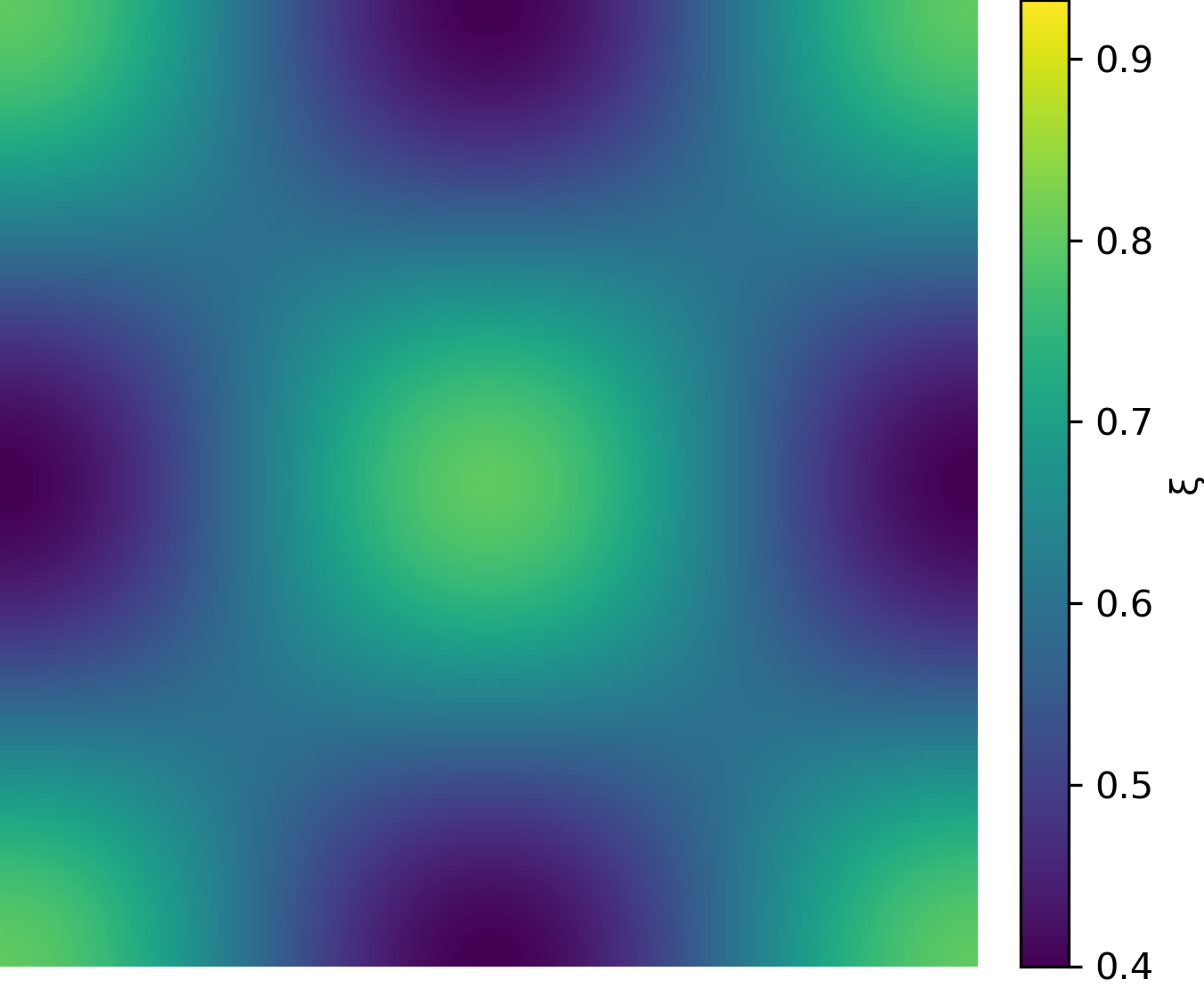}
  \caption{$\xi_{\mathrm{true}}$}
\end{subfigure}\hfill
\begin{subfigure}{0.24\textwidth}
  \centering
  \includegraphics[width=\linewidth]{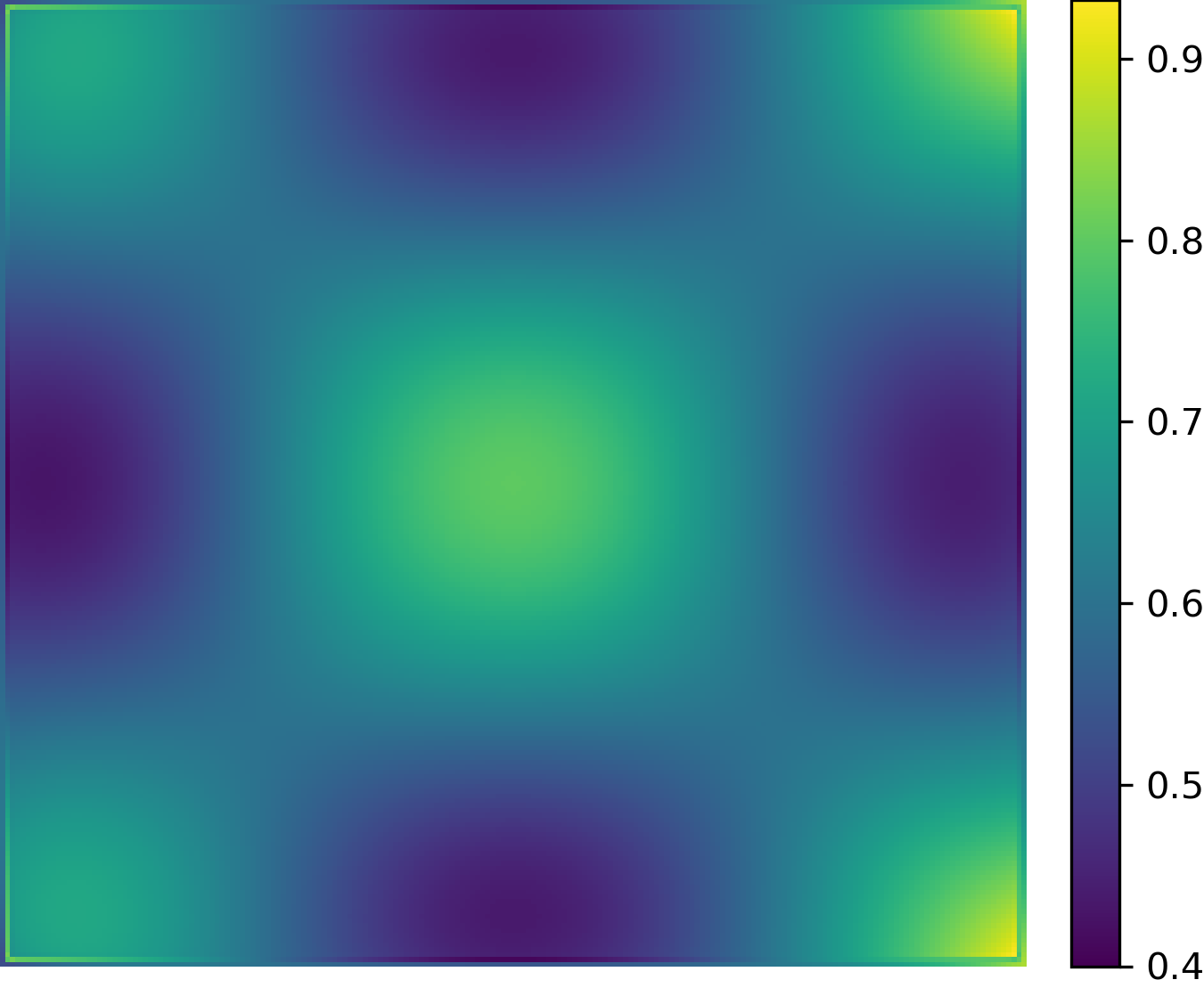}
  \caption{$\xi_{\mathrm{rec}}$}
\end{subfigure}\hfill
\begin{subfigure}{0.24\textwidth}
  \centering
  \includegraphics[width=\linewidth]{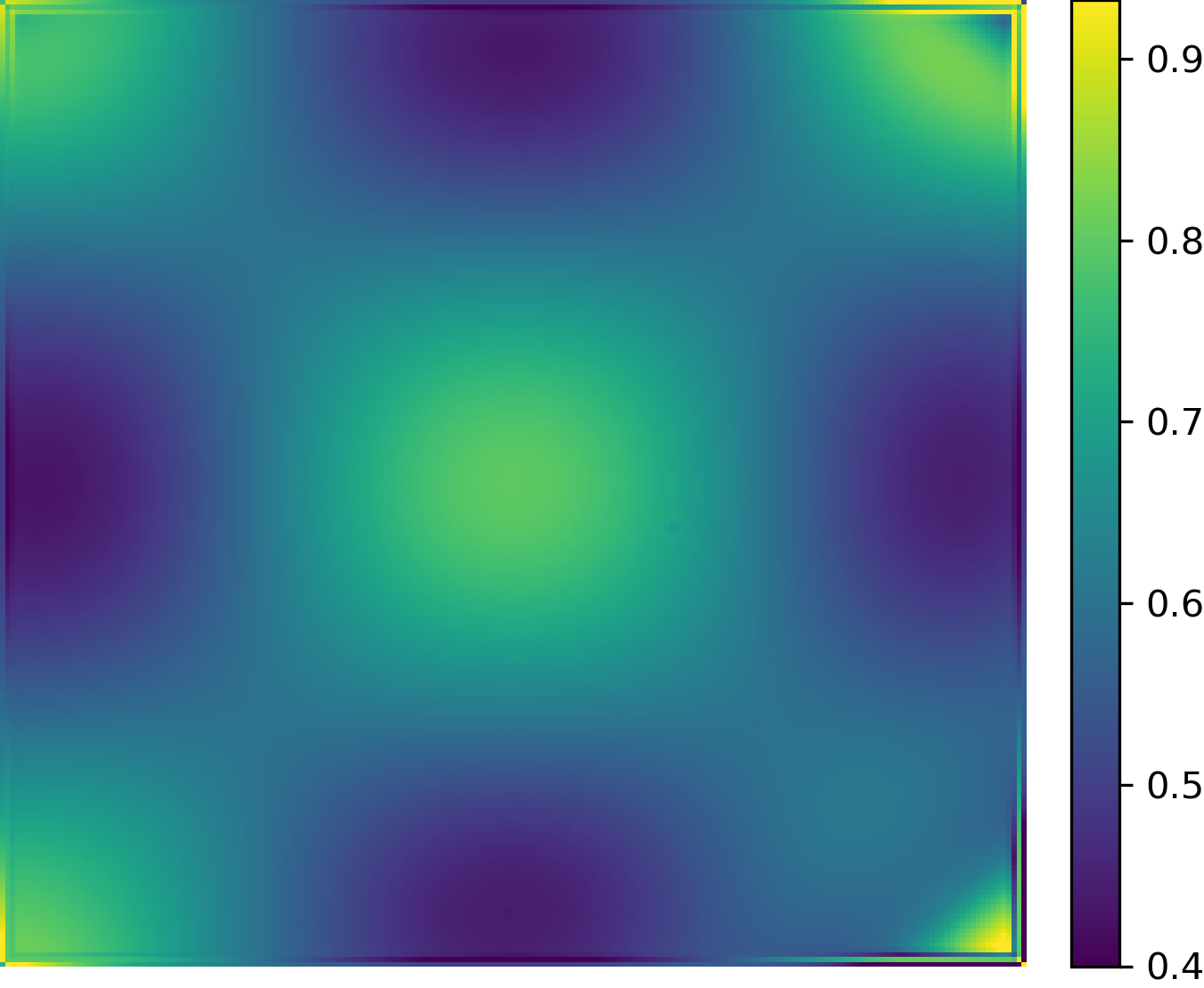}
  \caption{$\xi_{\mathrm{rec}}$ (1\% noise)}
\end{subfigure}\hfill
\begin{subfigure}{0.24\textwidth}
  \centering
  \includegraphics[width=\linewidth]{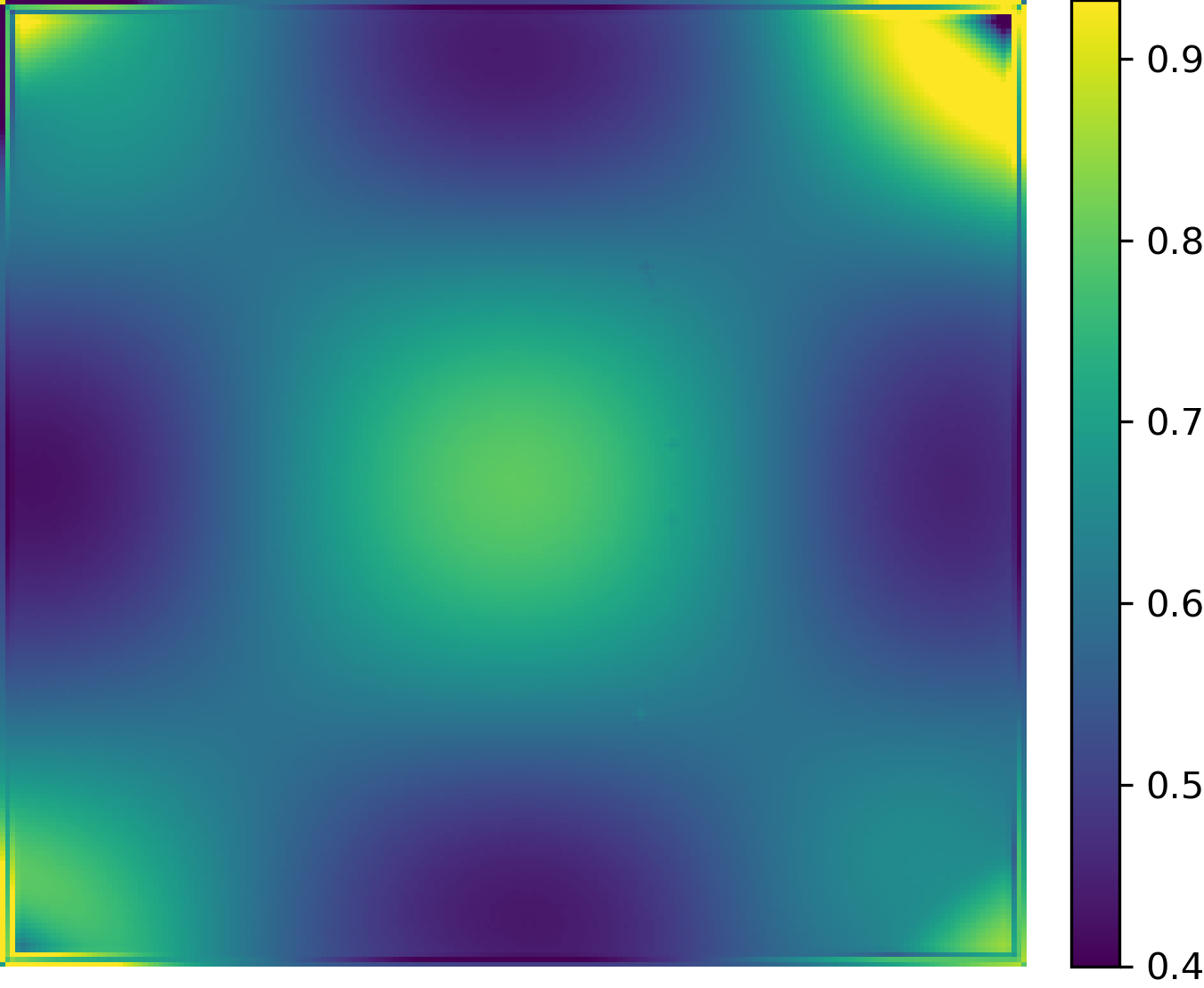}
  \caption{$\xi_{\mathrm{rec}}$ (5\% noise)}
\end{subfigure}

\caption{
Reconstruction results for the coefficients.
Top row: $\mu$. Bottom row: $\xi$.
From left to right: true coefficient, noiseless reconstruction,
reconstruction with 1\% noise, and reconstruction with 5\% noise.
}
\label{fig:coeff_recon_2}
\end{figure}

\paragraph{Experiment II.}
In our second numerical experiment, we focus on reconstructing the reaction coefficients $(\mu,\xi)$ given by:
\begin{equation}
\begin{aligned}
\mu(x,y) &= 0.4 + 0.25(1 - x^2) + 0.1 \sin(\pi y), \\
\xi(x,y) &= 0.2 + 0.6 \cos(\pi x)\cos(\pi y).
\end{aligned}
\end{equation}
The resulting reconstructions are visualized in Figure~\ref{fig:coeff_recon_2} across noise levels $\eta = 0$, $0.01$, and $0.05$ (arranged left to right), all started with the initial condition $\rho_0(x,y) = 0.5 + 0.2 \cos(0.5\pi x) + 0.15 \cos(0.5\pi y)$. The method continues to demonstrate high accuracy, yielding relative $L^2$ errors for $(\mu,\xi)$ of $(4.49,\,3.62)\times10^{-2}$, $(8.33,\,6.47)\times10^{-2}$, and $(2.44,\,1.48)\times10^{-1}$, respectively.

\paragraph{Experiment III.}
The third experiment illustrates a case of reconstruction failure. In this setup, $\mu$ is defined as a constant background superimposed with four Gaussian sources of varying amplitudes, locations, and orientations, while $\xi$ is given by $0.6 + 0.4 \cos(\pi x) \cos(\pi y)$. The observed failure is primarily attributed to the reconstruction's sensitivity to the initial guess. For this specific trial, only datasets $d_0, d_1$, and $d_2$ were utilized, with both $\mu^{(0)}$ and $\xi^{(0)}$ initialized as $1.0 + 0.3 (x^2 + y^2)$. 

\begin{figure}[!htb]
\centering
\begin{subfigure}{0.24\textwidth}
  \centering
  \includegraphics[width=\linewidth]{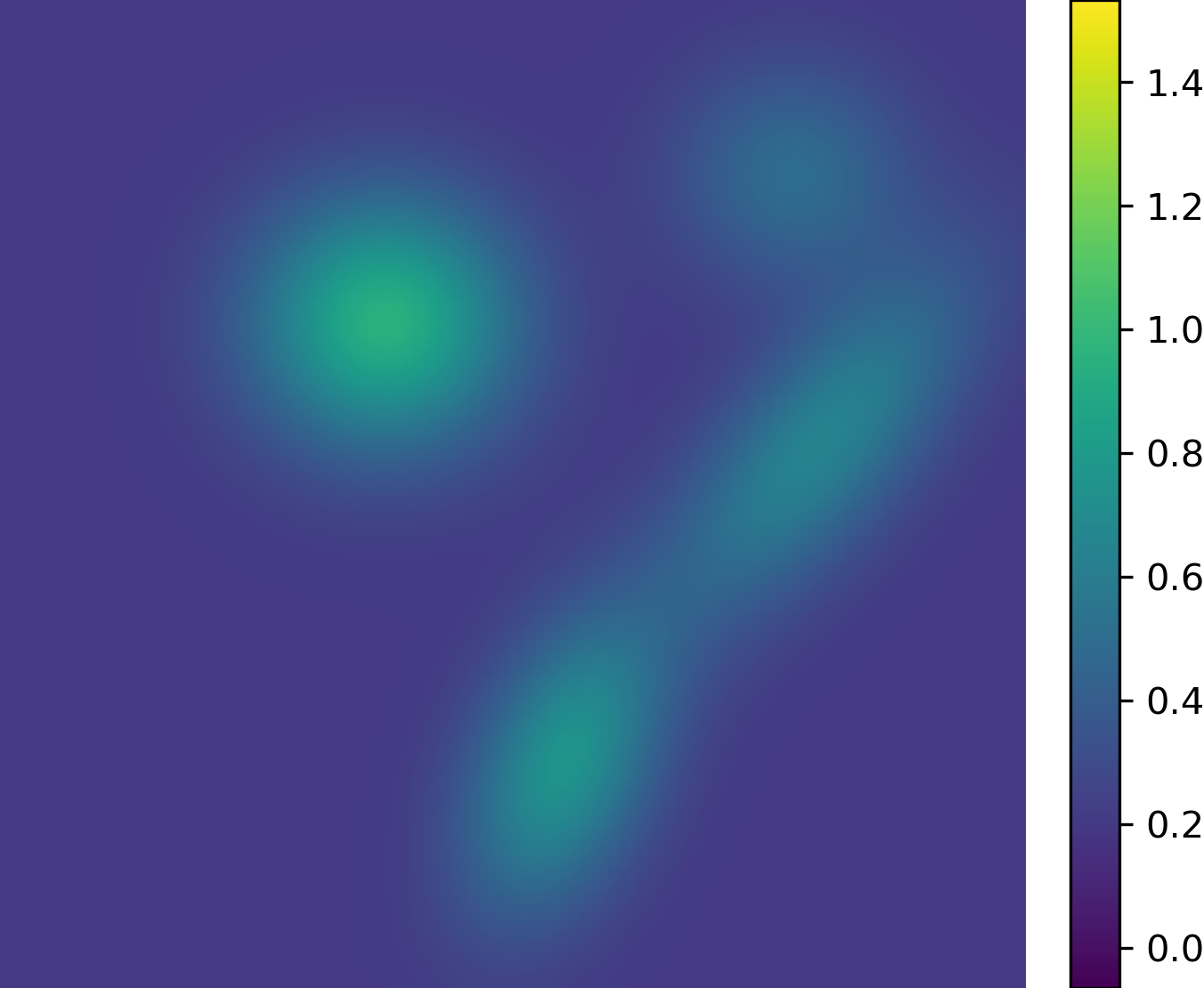}
  \caption{$\mu_{\mathrm{true}}$}
\end{subfigure}\hfill
\begin{subfigure}{0.24\textwidth}
  \centering
  \includegraphics[width=\linewidth]{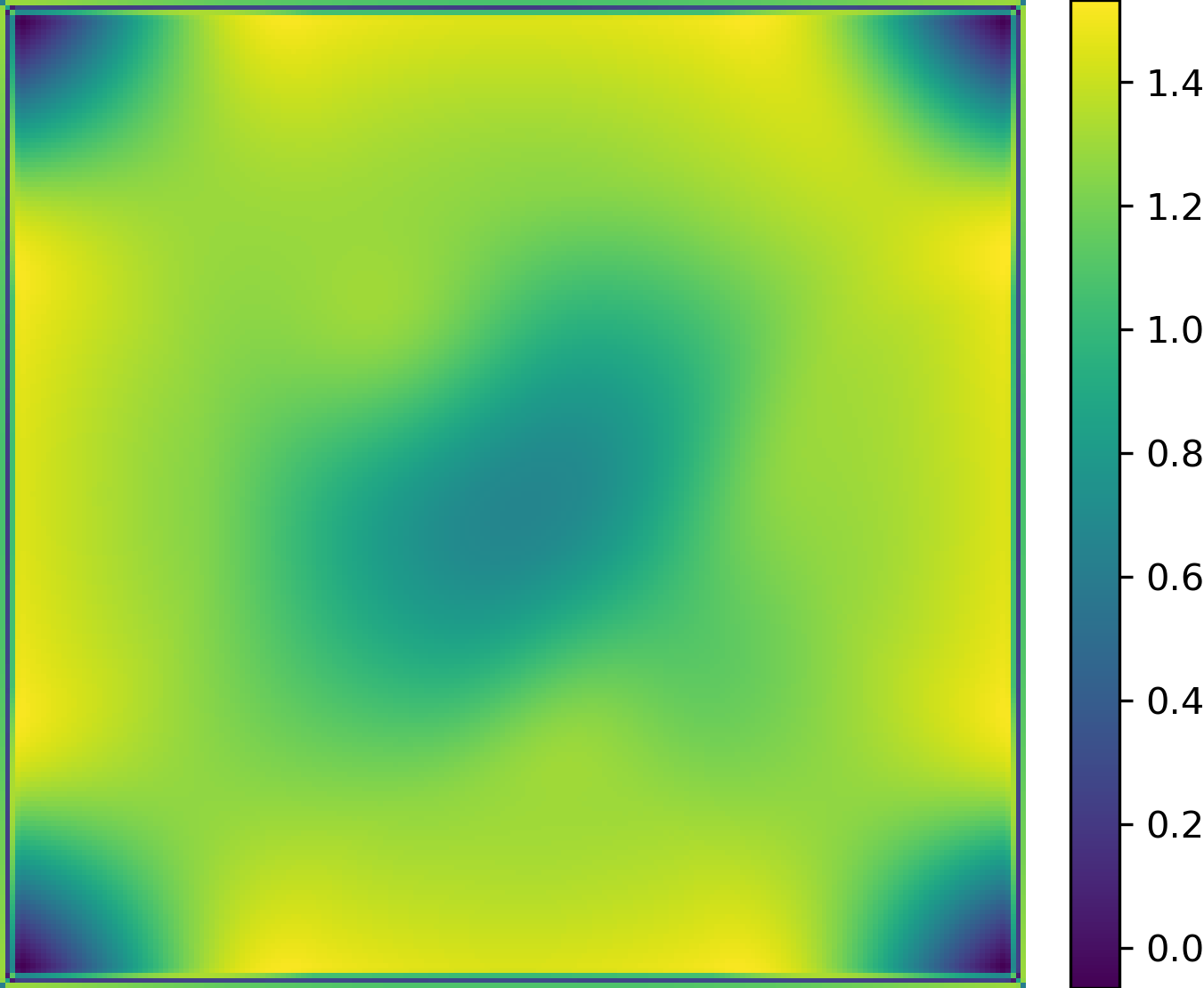}
  \caption{$\mu_{\mathrm{rec}}$}
\end{subfigure}\hfill
\begin{subfigure}{0.24\textwidth}
  \centering
  \includegraphics[width=\linewidth]{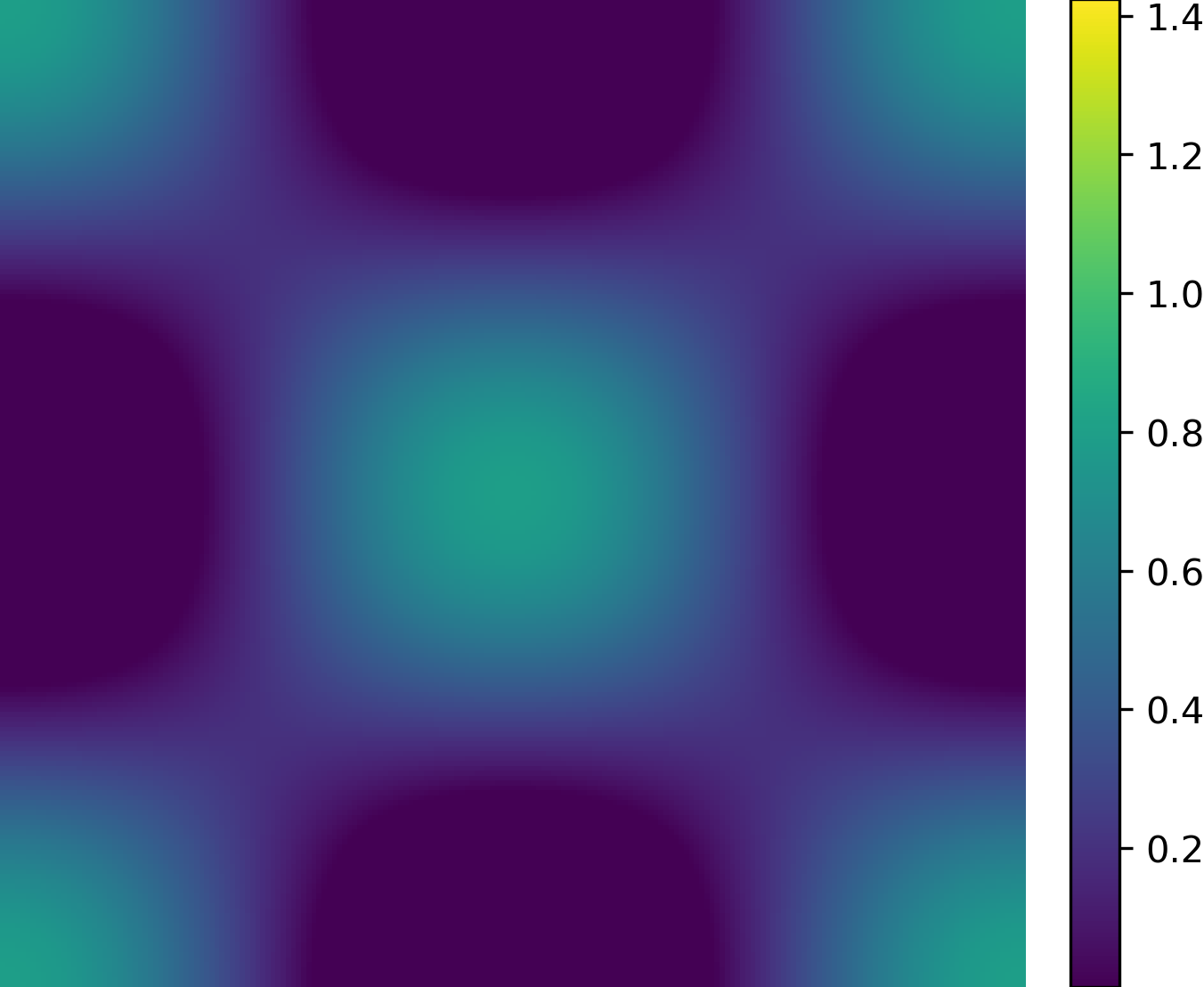}
  \caption{$\xi_{\mathrm{true}}$}
\end{subfigure}\hfill
\begin{subfigure}{0.24\textwidth}
  \centering
  \includegraphics[width=\linewidth]{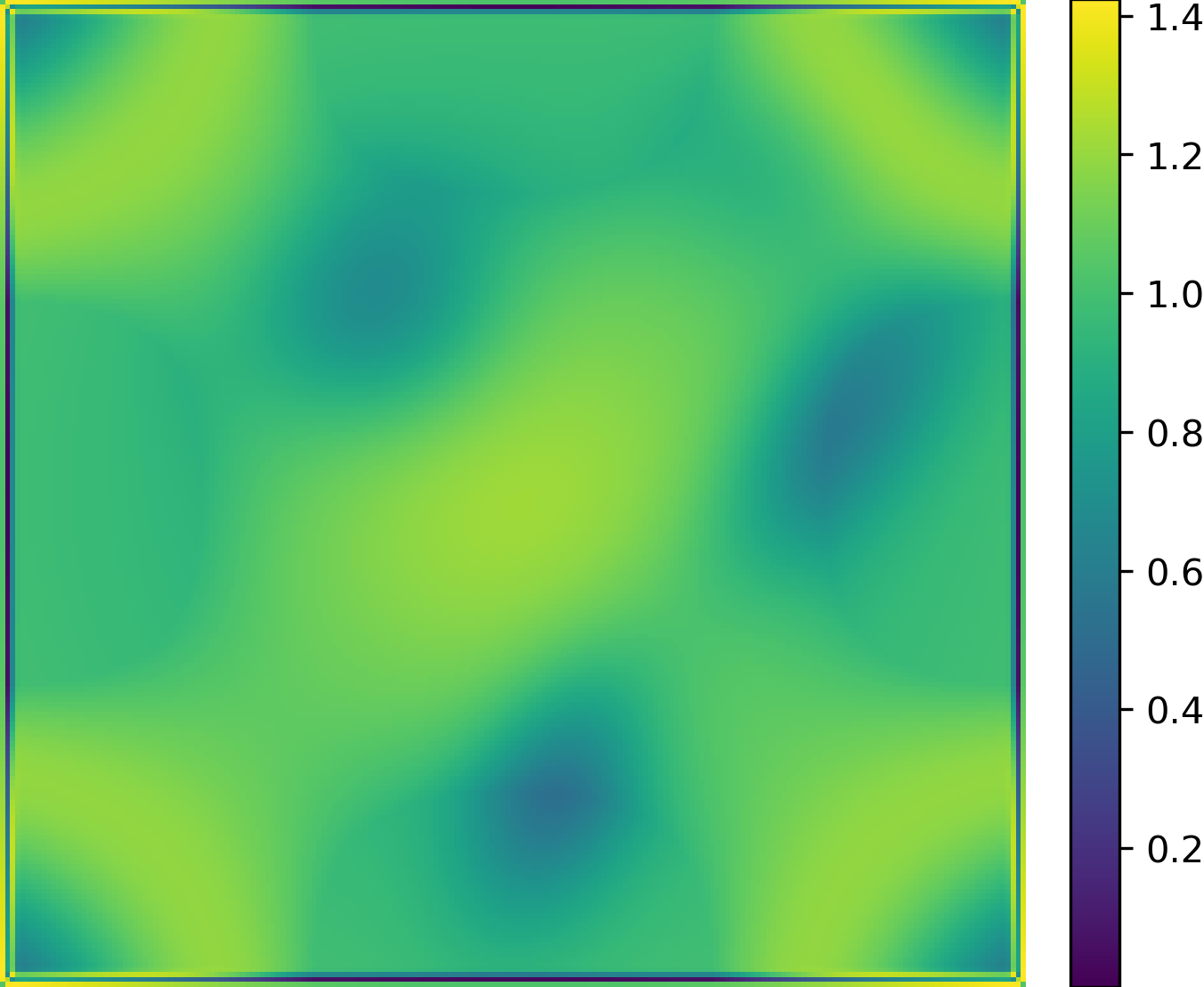}
  \caption{$\xi_{\mathrm{rec}}$}
\end{subfigure}

\caption{
Reconstruction failure for the coefficients.
From left to right: ground truth of $\mu$, noiseless reconstruction of $\mu$,
ground truth of $\xi$, and reconstruction of $\xi$.
}
\label{fig:coeff_recon_3}
\end{figure}

As shown in Figure \ref{fig:coeff_recon_3}, the reconstructed profiles of $\mu$ and $\xi$ are strongly coupled, failing to separate the distinct features of each parameter. Intuitively, this arises from the properties of the sensitivity matrix, whose columns are not guaranteed to be independent. For instance, the snapshots $d_1$ and $d_2$ become nearly linearly dependent when the sampling times $t_1$ and $t_2$ are close, leading to an ill-conditioned inverse problem.

\begin{figure}[!htb]
\centering
\begin{subfigure}{0.24\textwidth}
  \centering
  \includegraphics[width=\linewidth]{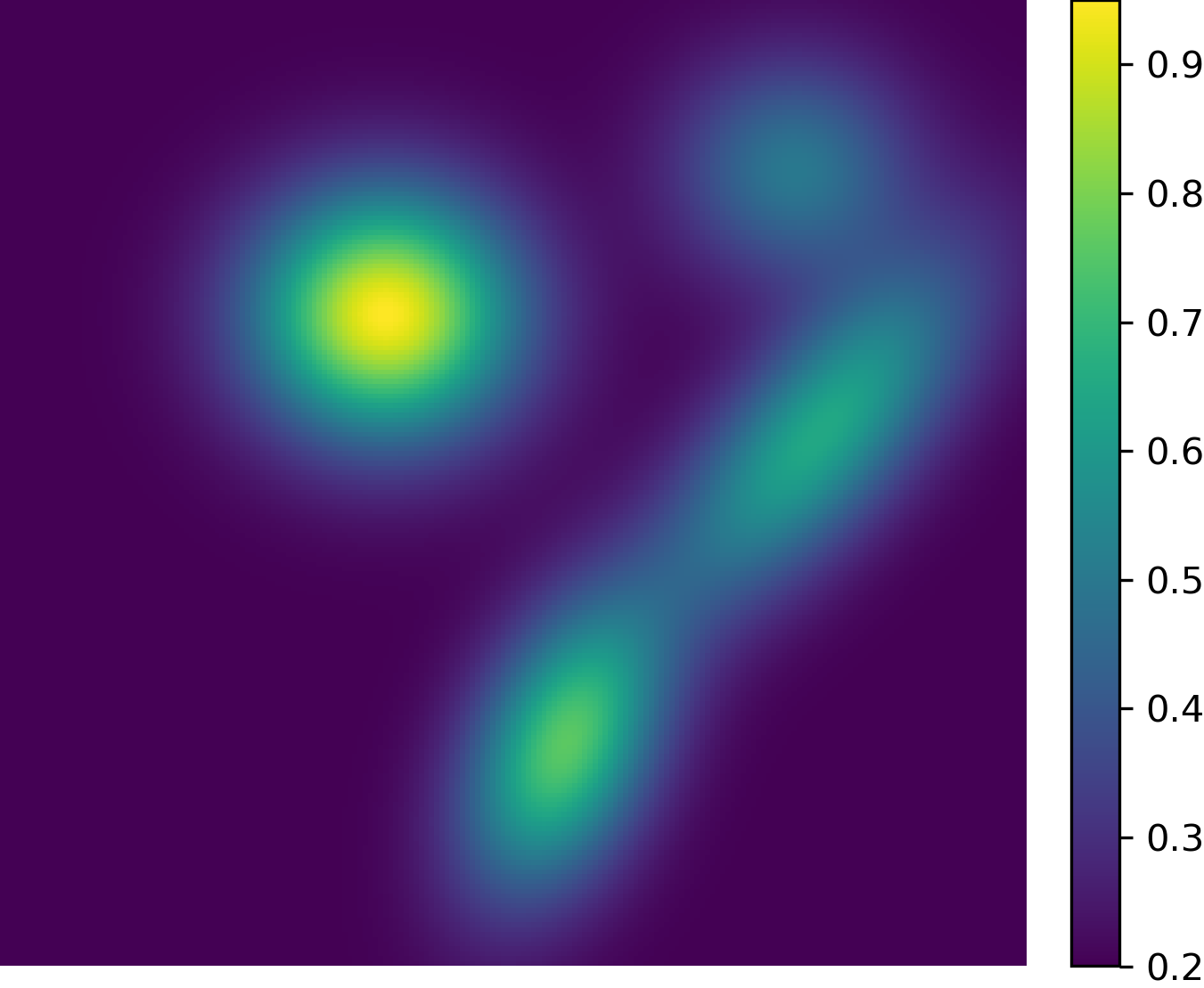}
  \caption{Ground truth}
\end{subfigure}\hfill
\begin{subfigure}{0.24\textwidth}
  \centering
  \includegraphics[width=\linewidth]{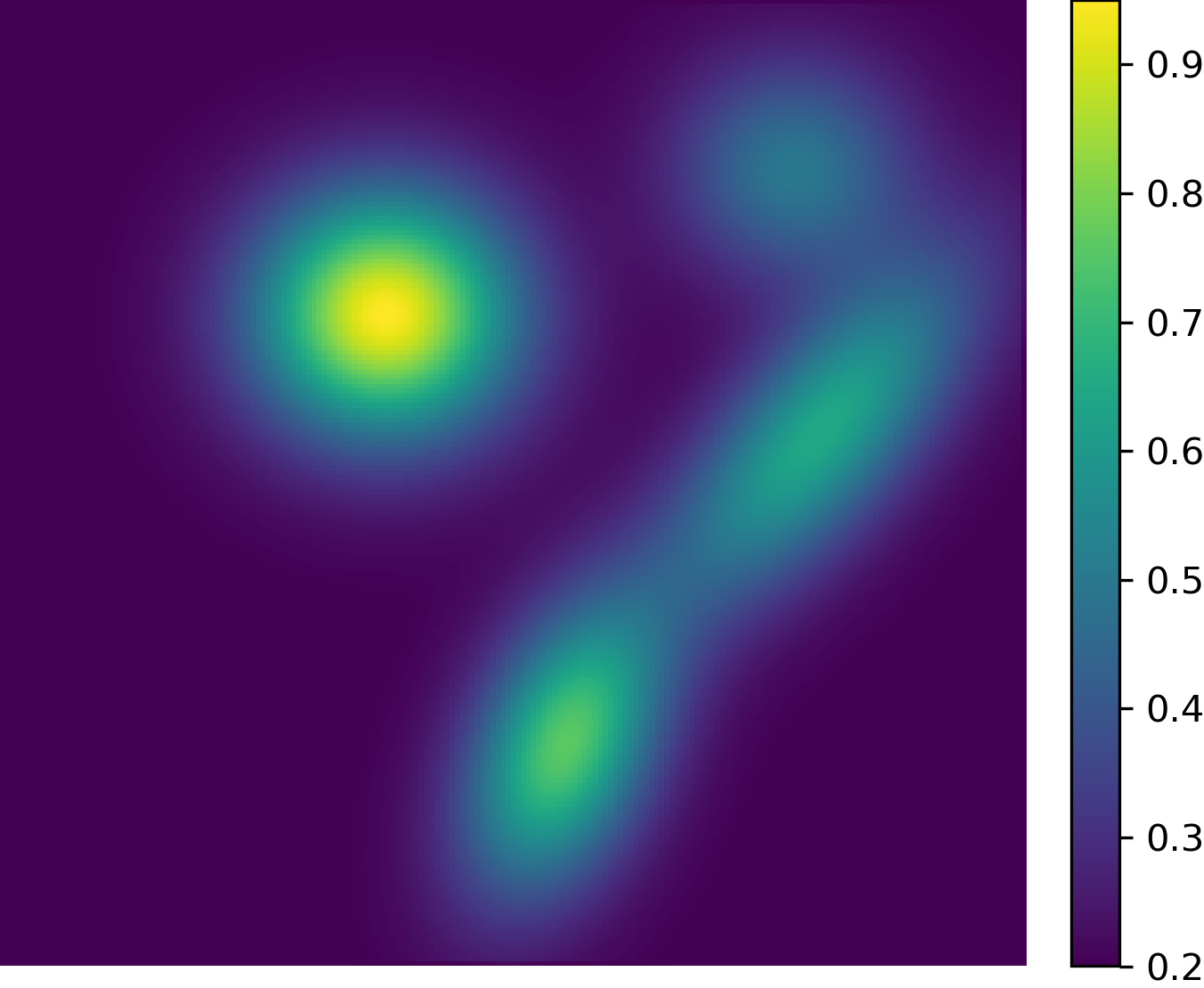}
  \caption{Noise free}
\end{subfigure}\hfill
\begin{subfigure}{0.24\textwidth}
  \centering
  \includegraphics[width=\linewidth]{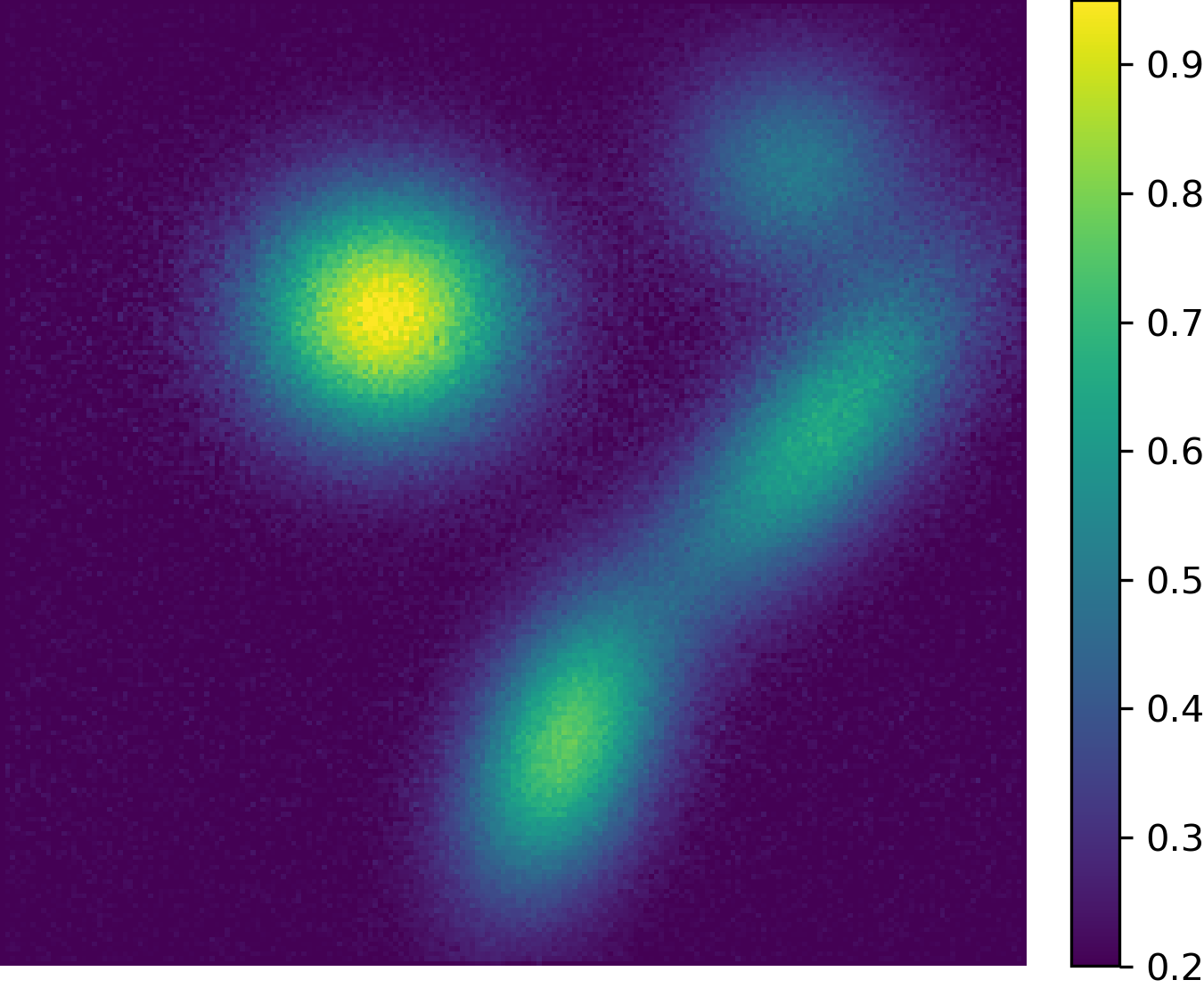}
  \caption{$5\%$ Noise}
\end{subfigure}\hfill
\begin{subfigure}{0.24\textwidth}
  \centering
  \includegraphics[width=\linewidth]{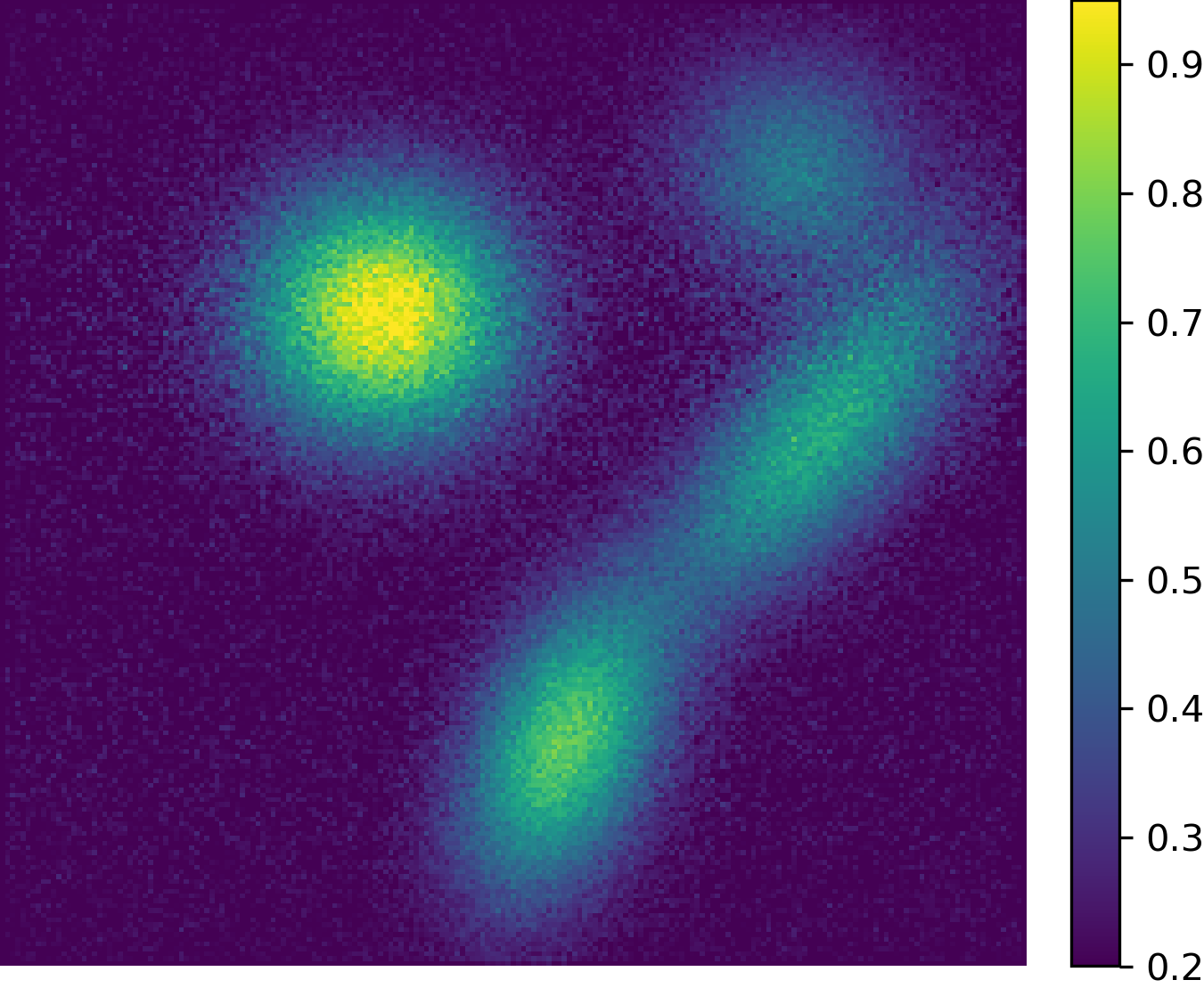}
  \caption{$10\%$ Noise}
\end{subfigure}

\caption{
Reconstruction of the initial condition given by a superposition of four Gaussian sources. 
From left to right, we show the ground truth $\rho_0$, the reconstruction using noiseless data, and the reconstructions obtained with $5\%$ and $10\%$ noise.
}
\label{fig:coeff_recon_4}
\end{figure}

\paragraph{Experiment IV.} In this experiment, we consider the reconstruction of an initial condition consisting of a superposition of four Gaussian sources with varying amplitudes, locations, and orientations. 
We performed simulations using noisy data with noise-to-signal ratios (NSR) $\eta = 0$, $\eta = 0.05$, and $\eta = 0.1$. 
The proposed method achieves high-quality reconstructions, with relative $L^2$ errors equal to $(9.22,\,10.49,\,13.84)\times10^{-2}$, respectively. 
The corresponding reconstruction results are displayed in Figure~\ref{fig:coeff_recon_4}.

\begin{figure}[!htb]
\centering
\begin{subfigure}{0.24\textwidth}
  \centering
  \includegraphics[width=\linewidth]{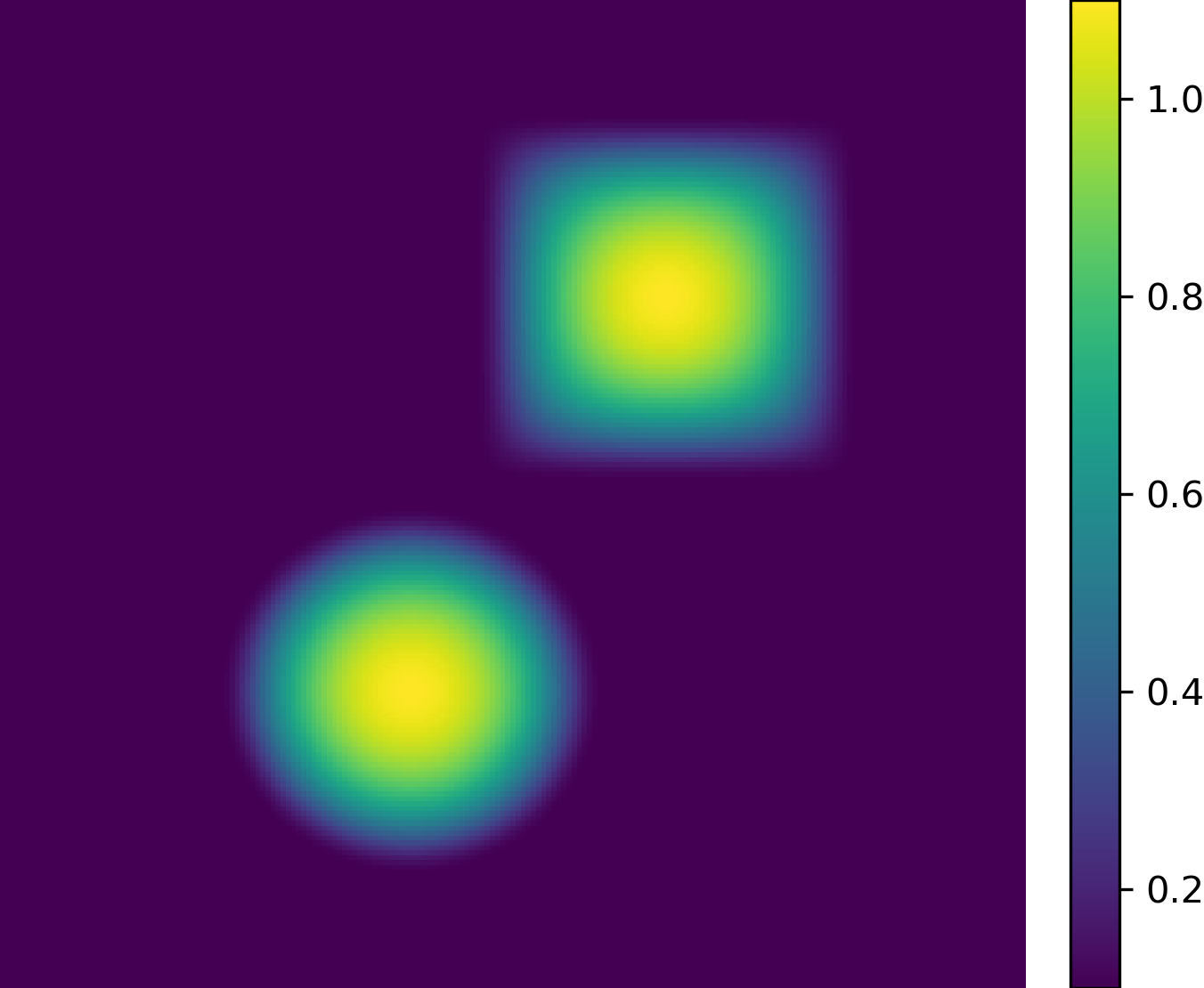}
  \caption{Ground truth}
\end{subfigure}\hfill
\begin{subfigure}{0.24\textwidth}
  \centering
  \includegraphics[width=\linewidth]{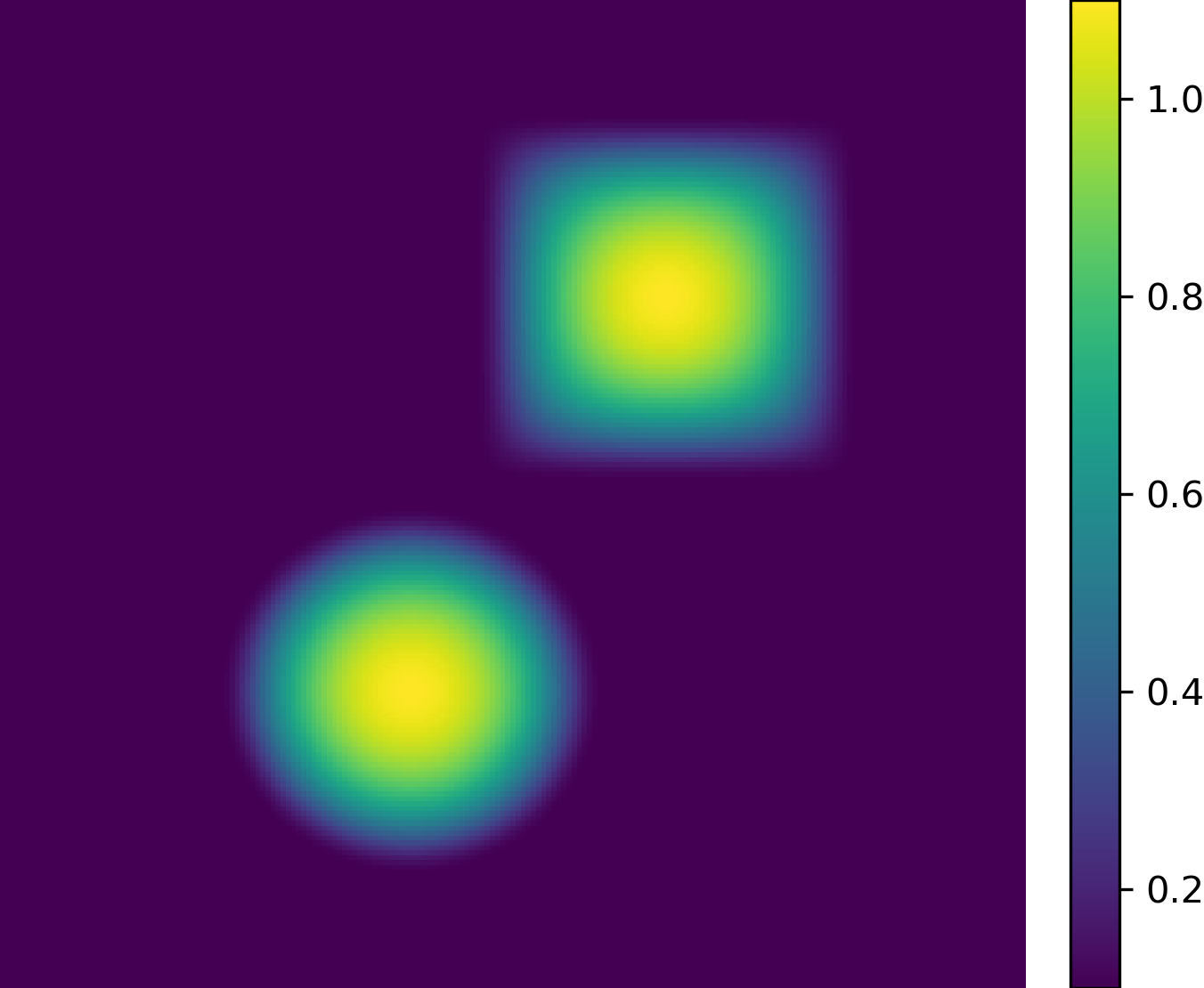}
  \caption{Noise free}
\end{subfigure}\hfill
\begin{subfigure}{0.24\textwidth}
  \centering
  \includegraphics[width=\linewidth]{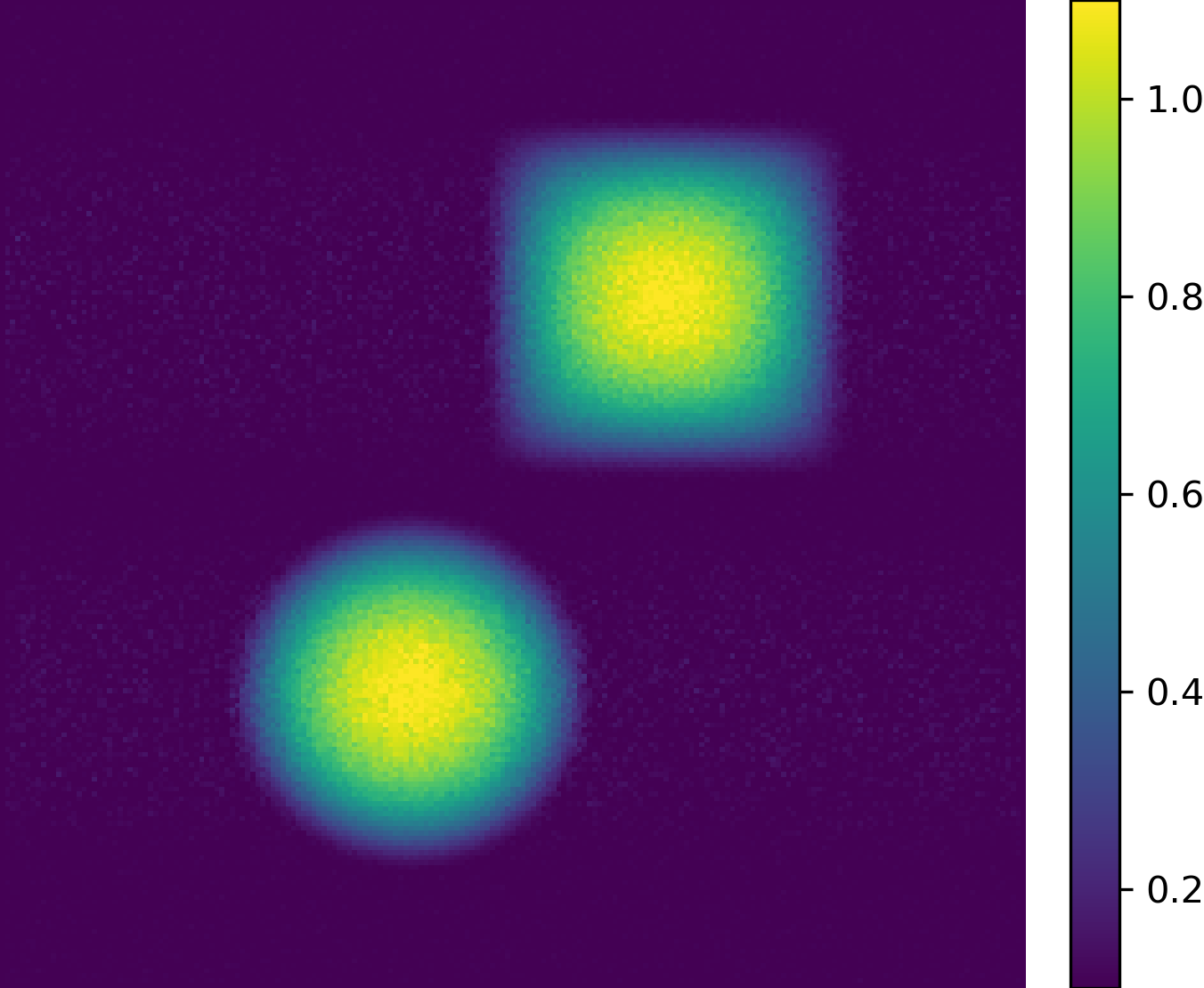}
  \caption{$5\%$ Noise}
\end{subfigure}\hfill
\begin{subfigure}{0.24\textwidth}
  \centering
  \includegraphics[width=\linewidth]{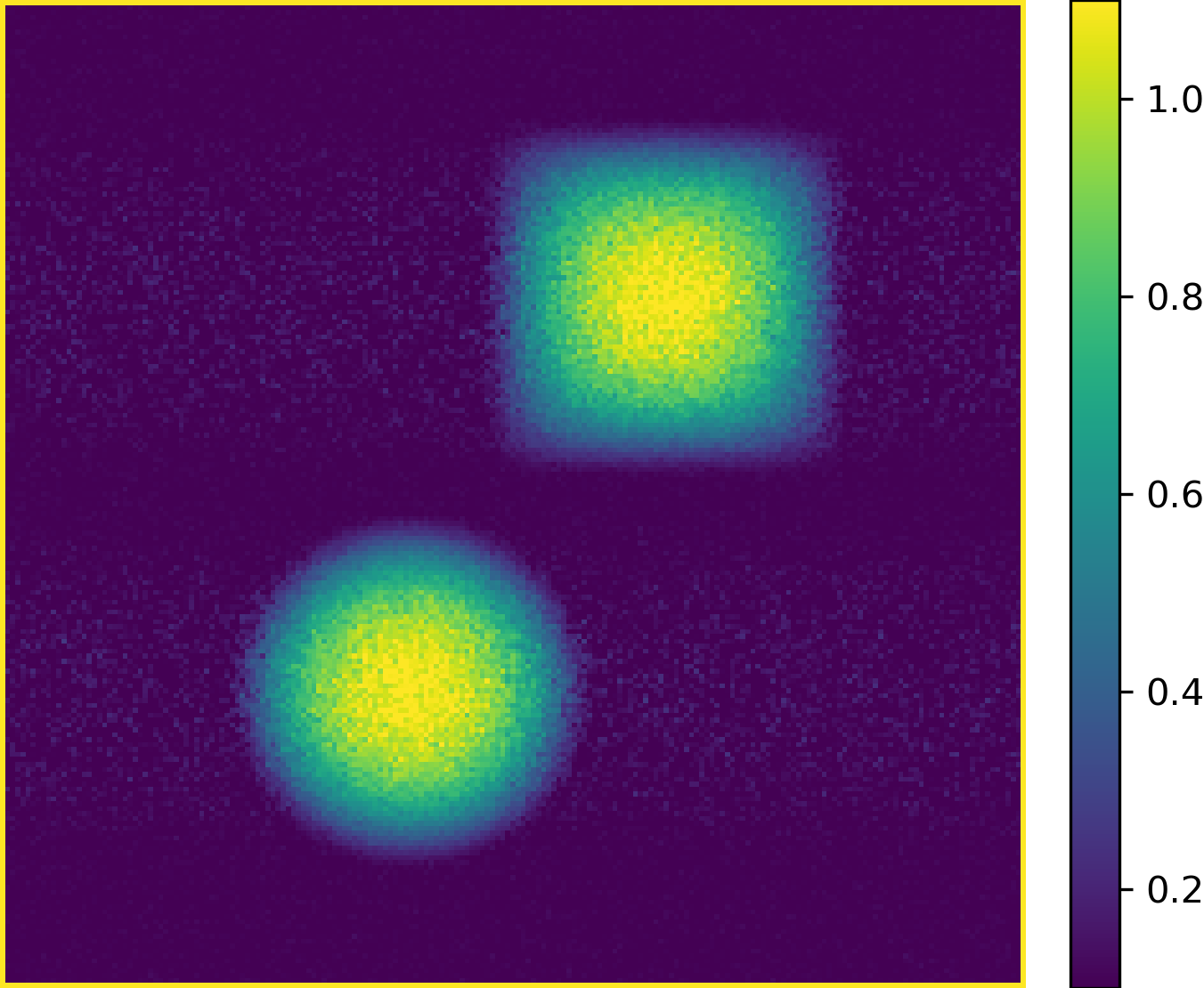}
  \caption{$10\%$ Noise}
\end{subfigure}

\caption{
Reconstruction of the initial condition consisting of a circle and a square.
From left to right: ground truth of $\rho_0$, noiseless reconstruction, reconstruction with $5\%$ noise, and reconstruction with $10\%$ noise.
}
\label{fig:coeff_recon_5}
\end{figure}

\paragraph{Experiment V.} This experiment is devoted to the reconstruction of an initial condition consisting of a circle and a square. 
Figure~\ref{fig:coeff_recon_5} presents the reconstruction results obtained with NSR $\eta = 0$, $\eta = 0.05$, and $\eta = 0.1$. 
The corresponding relative $L^2$ errors are $(9.94,\,11.50,\,57.61)\times10^{-2}$, respectively.

\begin{figure}[!htb]
\centering
\begin{subfigure}{0.24\textwidth}
  \centering
  \includegraphics[width=\linewidth]{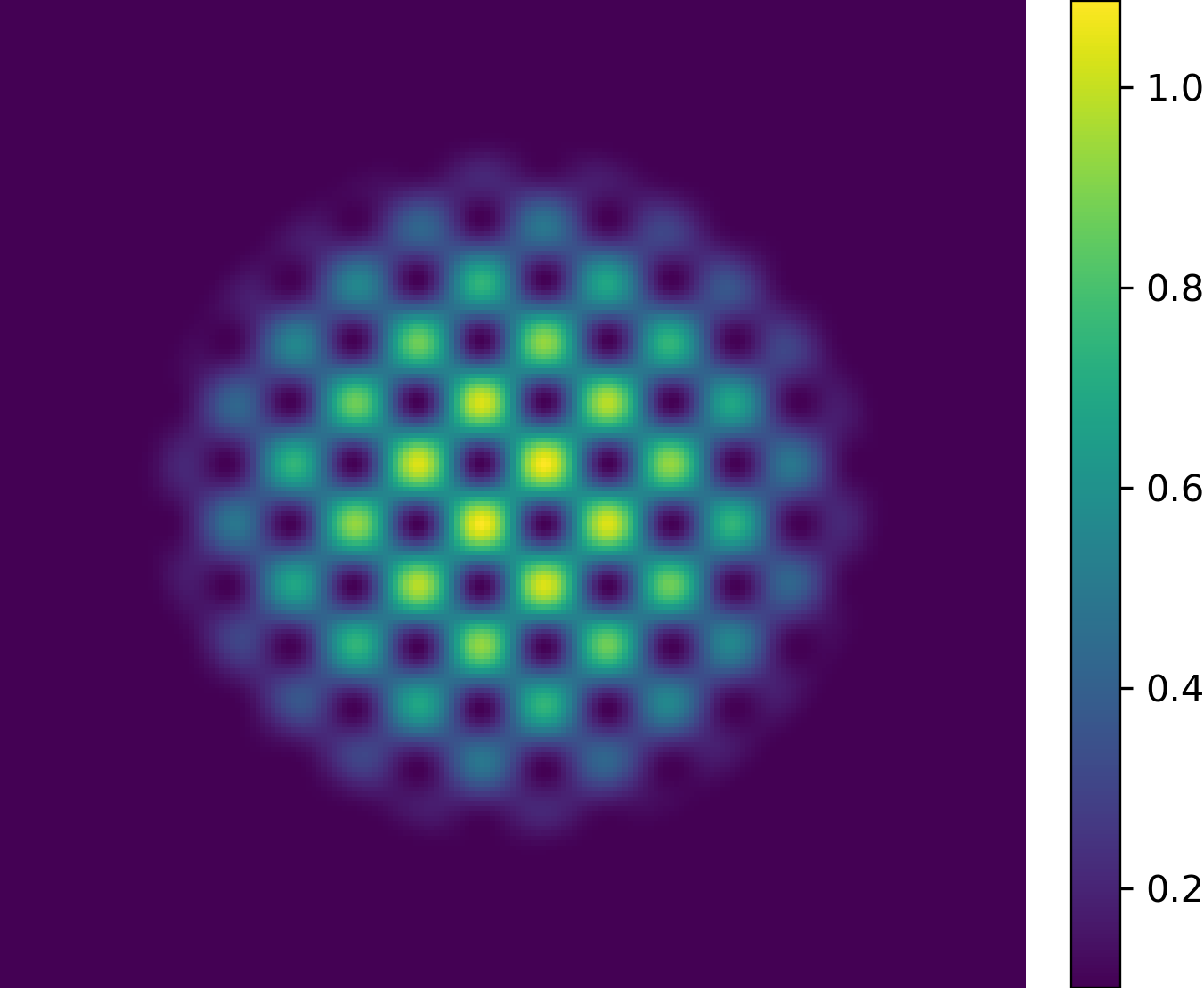}
  \caption{Ground truth}
\end{subfigure}\hfill
\begin{subfigure}{0.24\textwidth}
  \centering
  \includegraphics[width=\linewidth]{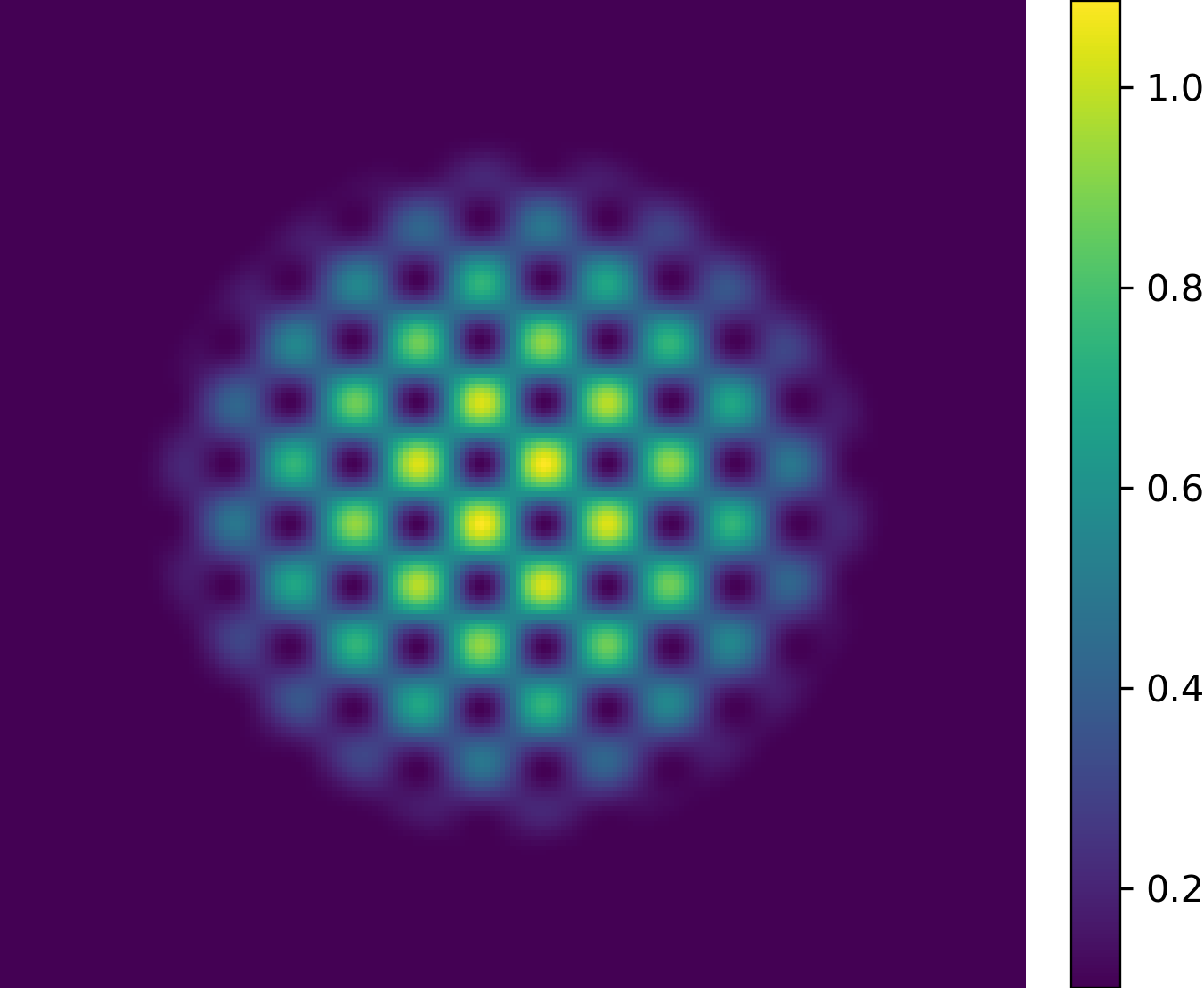}
  \caption{Noise free}
\end{subfigure}\hfill
\begin{subfigure}{0.24\textwidth}
  \centering
  \includegraphics[width=\linewidth]{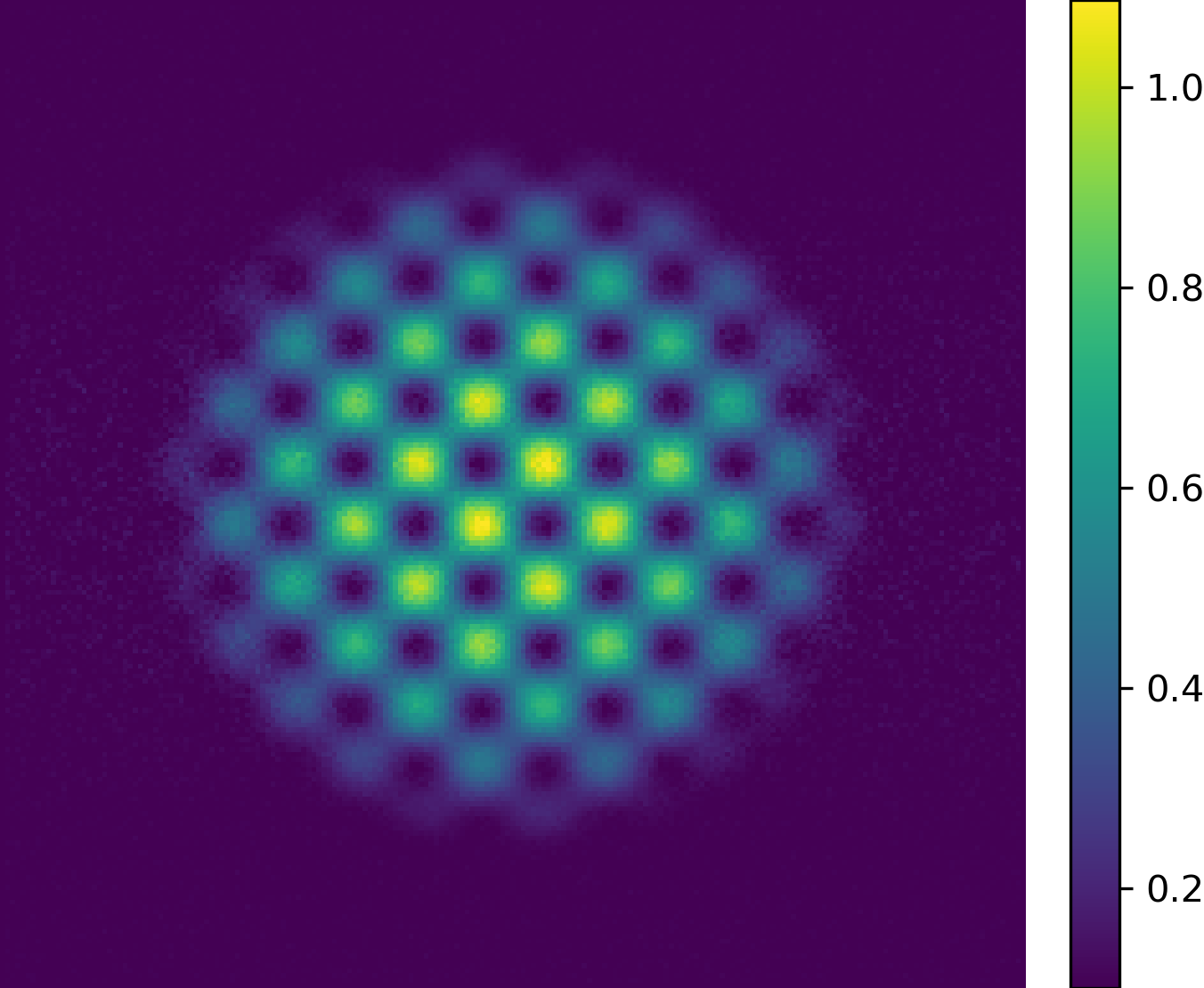}
  \caption{$5\%$ Noise}
\end{subfigure}\hfill
\begin{subfigure}{0.24\textwidth}
  \centering
  \includegraphics[width=\linewidth]{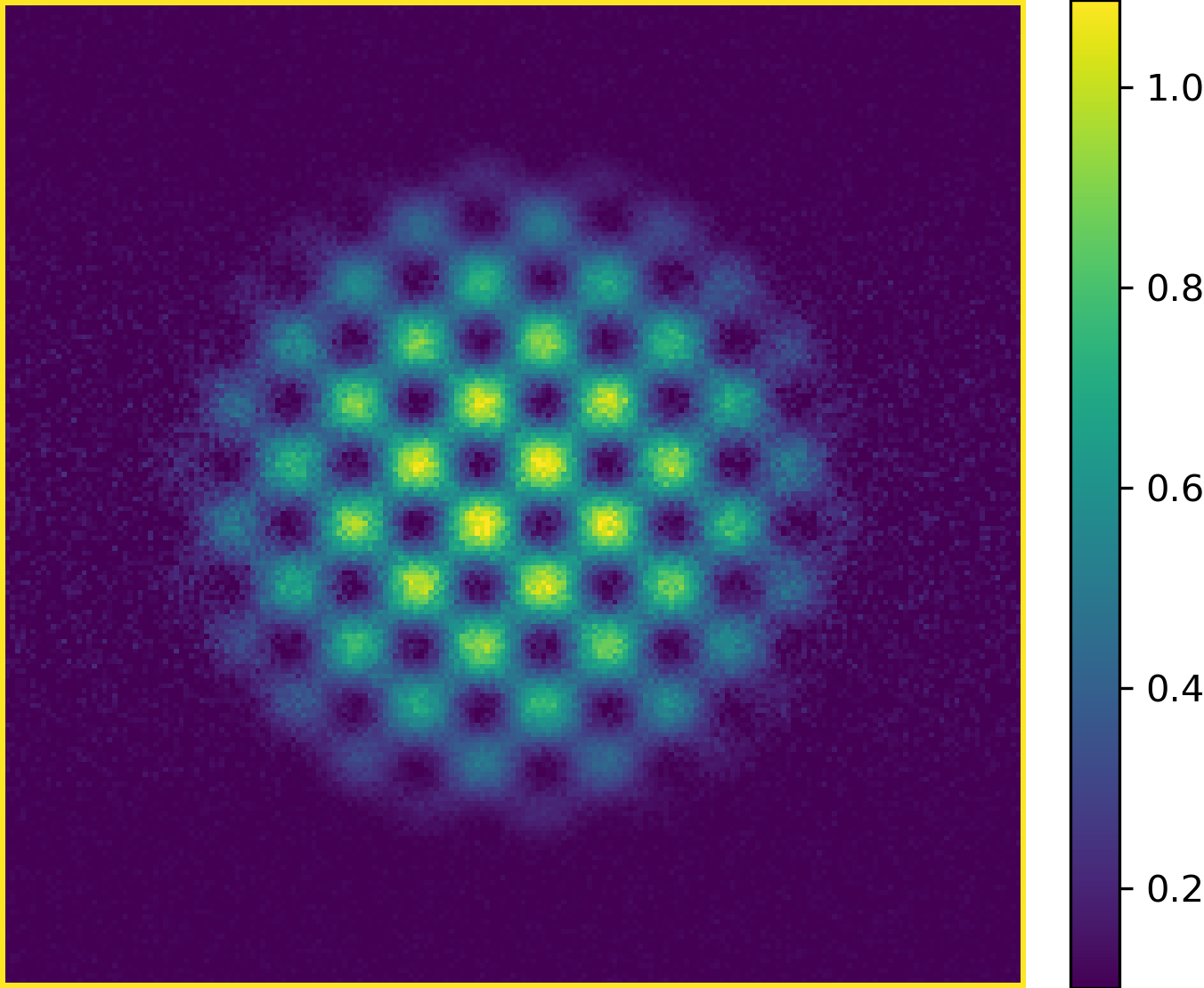}
  \caption{$10\%$ Noise}
\end{subfigure}

\caption{
Reconstruction of the initial condition of a localized wave packet.
From left to right: ground truth of $\rho_0$, noiseless reconstruction, reconstruction with $5\%$ noise, and reconstruction with $10\%$ noise.
}
\label{fig:coeff_recon_6}
\end{figure}

\paragraph{Experiment VI.} In the final experiment, we consider the reconstruction of a localized wave packet. 
As shown in Figure~\ref{fig:coeff_recon_6}, the relative $L^2$ reconstruction errors are noticeably larger. 
Specifically, for noise levels of $0\%$, $5\%$, and $10\%$, the corresponding errors are $1.24\times10^{-1}$, $1.31\times10^{-1}$, and $6.73\times10^{-1}$, respectively. 
This reduction in reconstruction accuracy is primarily caused by the high-frequency oscillatory nature of the wave packet, which amplifies the ill-posedness of the inverse problem.

The above examples validate the feasibility of the proposed algorithm and corroborate the theoretical analysis. We further observe that, in practice, the inversion of the initial condition is numerically more robust than the simultaneous inversion of the two reaction coefficients. While this may at first appear to contradict the stability hierarchy of Theorem~\ref{stabi}, where the coefficients enjoy Lipschitz stability and the initial condition only logarithmic stability, the coefficient inversion is hampered by the strong coupling between $\mu$ and $\xi$, which renders the simultaneous recovery of both intrinsically ill-conditioned. In particular, during the decoupling procedure, the reconstruction errors introduced at each stage propagate through the forward and adjoint PDE solves and accumulate across iterations, thereby amplifying the instability of the coefficient inversion.

For the inversion of the initial condition, we adopt a smooth Gaussian function as the initial guess. We further perform numerical experiments to investigate the influence of the measurement time on the reconstruction accuracy. The results indicate that measurements taken closer to the initial time yield more accurate reconstructions, as expected from the smoothing property of the forward parabolic operator, although the overall sensitivity to the measurement time remains moderate.

Moreover, using the two-stage strategy, we are able to recover both the coefficients and the initial condition from the same measurement data. In this approach, the initial condition is reconstructed using the previously recovered coefficients. Within a controllable noise regime, the errors in the coefficient reconstruction do not noticeably degrade the stability of the initial-condition inversion, which further demonstrates the robustness of the proposed algorithm.

\section{Concluding remarks}
\label{SEC:Conclusion}

In this paper, we studied an inverse problem for a nonlinear, density-dependent reaction-diffusion model of cell invasion and brain tumor growth, in which the local proliferation rate $\mu(x)$, the competition (saturation) coefficient $\xi(x)$, and the unknown initial condition $\rho_0(x)$ are reconstructed simultaneously from interior measurements on a subdomain together with a few full-domain snapshots. On the theoretical side, we established, by means of Carleman estimates, a global uniqueness result and a Lipschitz-type stability estimate for the reaction coefficients, together with a weaker logarithmic stability estimate for the initial condition that reflects the severe ill-posedness inherent in recovering the initial state. On the computational side, we introduced a time-shift strategy that treats an interior snapshot as a pseudo-initial condition, thereby decoupling the strongly coupled problem into a two-stage scheme: the reaction coefficients are recovered first, without prior knowledge of the true initial state, and the initial condition is reconstructed afterwards with the coefficients held fixed. Both stages were formulated as PDE-constrained optimization problems and solved by an adjoint-state method. A series of two-dimensional numerical experiments confirmed the feasibility, accuracy, and noise robustness of the proposed approach, and also delineated its limitations for strongly coupled coefficient pairs and highly oscillatory profiles.

Several aspects of the present work suggest opportunities for refinement. First, our analysis relies on a non-degeneracy assumption on the diffusion coefficient, namely $0<\underline\alpha\le\gamma(\rho)\le\overline\alpha<\infty$, together with the non-degeneracy condition~\eqref{ndeg} on $\partial_t\rho$. Relaxing these assumptions to cover the degenerate, porous-medium-type regime $\gamma(\rho)=\beta(x)\rho^\kappa$, which is more faithful to the sharp invasion fronts observed in practice, is an important and challenging extension. Second, the diffusion coefficient $\gamma$ was assumed known throughout; recovering it simultaneously with the reaction coefficients and the initial state would broaden the applicability of the framework, at the cost of additional ill-posedness. Third, as Experiment~III illustrates, the simultaneous recovery of $\mu$ and $\xi$ is intrinsically ill-conditioned when the snapshot times are too close, or the initial guess is poor. A careful study of optimal experimental design, in particular the placement of the measurement subdomain $\omega$ and the choice of snapshot times $t_0,t_1,t_2$, would help mitigate this coupling and reduce the extra data requirement incurred by the decoupling strategy. Replacing the deterministic, regularization-based inversion by a Bayesian formulation would, in addition, provide principled uncertainty quantification for the reconstructed quantities.
 
\section*{Acknowledgments}

This work is partially supported by the National Science Foundation through grants DMS-1937254 and DMS-2309802, and partially by the Gordon \& Betty Moore Foundation Award GBMF12801 (doi.org/10.37807/GBMF12801). 

{\small
\bibliography{BIB-REN,BIB-IP-Tumor-Growth}
\bibliographystyle{siam} 
}


\end{document}